\documentclass{article} % default font size is 10pt

\newcommand{\dom}{\fn{dom}}
\newcommand{\cod}{\fn{cod}}
\newcommand{\sSet}{{\cat{sSet}}}

\usepackage{latexsym}
\usepackage{amsfonts}
\usepackage{amssymb}
\usepackage{amsmath}
\usepackage{amscd}
\usepackage{eucal}
\usepackage[all, knot]{xy}
\xyoption{arc}
\xyoption{rotate}
\usepackage{amsthm}
\usepackage{pict2e}
\usepackage{pstricks}
\usepackage{pst-node}
\usepackage{verbatim} % allows comments
\usepackage{graphics}
\usepackage{graphicx}
\usepackage{color}   
\usepackage{framed}

\usepackage{bbm}

\usepackage{marvosym}  %not in macros-beamer-pst
\usepackage{stmaryrd}
\usepackage{textcomp}
\usepackage{wasysym}

\psset{unit=1mm,linewidth=0.5pt,
arrowlength=1.1, arrowinset=0.7,arrowsize=4.5pt,
doublesep=1pt,
labelsep=2pt,nodesep=2pt}
  
\def\1{\ensuremath{\mathbbm{1}}}%
\def\2{\ensuremath{\mathbbm{2}}}%

\newcommand{\N}{\ensuremath{\mathbb N}} %xx
\newcommand{\E}{\ensuremath{\bb{E}}}

\newcommand{\cV}{{\cl{V}}}

\newcommand{\bL}{\bb{L}}
\newcommand{\bM}{\bb{M}}

\newcommand{\bR}{\bb{R}}

\newcommand{\bw}{\ensuremath{\mathbf{w}}}

\newcommand{\ve}{\varepsilon}

\newcommand{\nt}{natural transformation}
\newcommand{\nts}{natural transformations}

\newcommand{\iso}{\cong} 
 
\newcommand{\lra}{\ensuremath{\longrightarrow}} %xx

\renewcommand{\:}{\colon}

\newcommand{\fn}[1]{\ensuremath{\mbox{\bfseries {\upshape {#1}\hspace{1.3pt}}}}}

\newcommand{\cat}[1]{\ensuremath{\textrm{\bfseries {\upshape {#1}}}}}

\newcommand{\Set}{{\cat{Set}}}
\newcommand{\Cat}{{\cat{Cat}}}

\newcommand{\Prof}{{\cat{Prof}}}

\newcommand{\Ladj}{\ensuremath{\bL\cat{Adj}}}
\newcommand{\Radj}{\ensuremath{\bR\cat{Adj}}}

\newcommand{\MADJ}{\ensuremath{\bM\cat{Adj}}}

\newcommand{\Mcat}{{\cat{Mcat}}}
\newcommand{\mcat}{{\cat{Mcat}}}
\newcommand{\CMcat}{{\cat{CMcat}}}
\newcommand{\cmcat}{{\cat{CMcat}}}

\newcommand{\tspan}{{\cat{$T$-Span}}}
\newcommand{\MAdj}{{\cat{MAdj}}}

\newcommand{\ha}{\ensuremath{\hat{a}}}

\newcommand{\shomblank}{\ensuremath{[\hh{2pt}\suscore\hh{1pt},\hh{1pt} \suscore\hh{2pt}]}}

\newcommand{\cl}[1]{\ensuremath{\mathcal {#1}}}
\newcommand{\bb}[1]{\ensuremath{\mathbb {#1}}}
\newcommand{\ed}{\end{document}}
\newcommand{\map}[1]{\ensuremath{\stackrel{{#1}}{\lra}}}

\newcommand{\demph}[1]{{\bfseries #1}}

\newcommand{\ttimes}{\hh{-2pt} \times \hh{-2pt}}

\newcommand{\op}{^{\hspace{1pt}\textrm{\upshape{op}}}}
\newcommand{\opb}{^\bullet}

\newcommand{\uscore}{\ensuremath{\underline{\hspace{0.6em}}}}
\newcommand{\suscore}{\ensuremath{\underline{\hspace{0.5em}}}}

\newcommand{\vv}[1]{\vspace*{#1}}
\newcommand{\hh}[1]{\hspace*{#1}}

\newcommand{\ladj}{\ensuremath{\dashv}}

\newcommand{\adj}{\dashv}

\newcommand{\tra}{\psset{unit=0.1cm,nodesep=0pt} \pspicture(8,0)
\pcline{->}(1,1.1)(7,1.1) \endpspicture}

\newcommand{\updown}{\psset{unit=0.1cm,nodesep=0pt} \pspicture(0,2)(0,2)
\pcline{<->}(1,0)(1,4) \endpspicture}

\newcommand{\updo}{\updown}

\newcommand{\tmapsto}{\psset{unit=0.1cm,nodesep=0pt} \pspicture(8,0) %shorter arrow
\pcline{|->}(1,1.2)(7,1.2) \endpspicture}

\newcommand{\ltmapsto}{\psset{unit=0.1cm,nodesep=0pt} \pspicture(14,0)
\pcline{|->}(2.5,1.2)(11.5,1.2) \endpspicture}

\newcommand{\tramap}[1]{\psset{unit=0.1cm,nodesep=0pt,labelsep=2pt} \pspicture(8,4)
\pcline{->}(1,1.1)(7,1.1)\naput{\ensuremath{\scriptstyle{#1}}} \endpspicture}

\newcommand{\tmap}{\tramap}

\newcommand{\ltramap}[1]{\psset{unit=0.1cm,nodesep=0pt} \pspicture(17,4)
\pcline{->}(2.5,1.1)(14.5,1.1)\naput{\ensuremath{\scriptstyle{#1}}} \endpspicture}

\newcommand{\ltmap}{\ltramap}

\newcommand{\prof}{\psset{unit=0.1cm,nodesep=0pt} \pspicture(10,0)
\pcline[labelsep=2pt]{->}(1,1.1)(9,1.1)
\pcline[tbarsize=3pt 1]{|-}(4.4,1.1)(8.5,1.1) \endpspicture}

\newcommand{\pra}{\prof}

\newcommand{\pmap}[1]{\psset{unit=0.1cm,nodesep=0pt} \pspicture(10,5)
\pcline{->}(1,1.1)(9,1.1)\naput[labelsep=2pt]{\ensuremath{\scriptstyle{#1}}}
\pcline[tbarsize=3pt 1]{|-}(4.4,1.1)(8.5,1.1) \endpspicture}

\newcommand{\cra}{\psset{unit=0.1cm,nodesep=0pt} \pspicture(10,0)
\pcline{->}(1,1.1)(9,1.1) \rput(5,1.1){\pscircle(0,0){1.6pt}}
 \endpspicture}
 
\newcommand{\cmap}[1]{\psset{unit=0.1cm,nodesep=0pt,labelsep=2pt} \pspicture(10,4)
\pcline{->}(1,1.1)(9,1.1)\naput[labelsep=2pt]{\ensuremath{\scriptstyle{#1}}} \rput(5,1.1){\pscircle(0,0){1.6pt}}
 \endpspicture}

\newcommand{\noi}{\noindent}
\newcommand{\dend}{\ensuremath{\hfill \maltese}\end{mydefinition}}

\newcommand{\igc}[1]{\begin{tabular}[c]{c}\includegraphics[width={#1}]}
\newcommand{\igt}[1]{\begin{tabular}[t]{c}\\[-24pt] \includegraphics[width={#1}]}
\newcommand{\ei}{\end{tabular}}

\newtheorem{theorem}{Theorem}[section]

\newtheorem{proposition}[theorem]{Proposition}

\newtheorem{lemma}[theorem]{Lemma}

\newtheorem{definition}[theorem]{Definition}

\newtheorem{example}[theorem]{Example}
\newtheorem{note}[theorem]{Note}
\newtheorem{nonexample}[theorem]{Non-example}
\newtheorem{examples}[theorem]{Examples}
\newtheorem{remarks}[theorem]{Remarks}
\newtheorem{exercise}[theorem]{Exercise}
\newtheorem{question}[theorem]{Question}
\newtheorem{questions}[theorem]{Questions}
\newtheorem{algorithm}[theorem]{Algorithm}
\newtheorem{method}[theorem]{Method}

\newtheorem{thm}[theorem]{Theorem}

\newtheorem{prop}[theorem]{Proposition}
\newtheorem{cor}[theorem]{Corollary}

\newenvironment{mydefinition}{\begin{definition}\upshape} {\end{definition}}
\newenvironment{myremark}{\begin{remark}\upshape} {\end{remark}}

\newenvironment{myexample}{\begin{example}\upshape} {\end{example}}

\newenvironment{myexamples}{\begin{examples}\upshape} {\end{examples}}

{ \end{sf}\end{framed}\end{minipage}
\end{center}}

\newtheoremstyle{example}{\topsep}{\topsep}%
     {}%         Body font
     {}%         Indent amount (empty = no indent, \parindent = para indent)
     {\bfseries}% Thm head font
     {.}%        Punctuation after thm head
     {8pt}%     Space after thm head (\newline = linebreak)
     {\thmname{#1}\thmnumber{ #2}\thmnote{ #3}}%         Thm head spec

   \theoremstyle{example}

\newtheorem{remark}[theorem]{Remark}

\newenvironment{prf}{\vspace{1ex}\begin{sloppypar}{\noindent\upshape
{\bfseries Proof.}}} {{\hspace*{\fill}
$\Box$}\end{sloppypar}\vspace{2ex}}

\newenvironment{prfof}[1]{\vspace{1ex}\begin{sloppypar}{\noindent
\upshape{\bfseries Proof of {#1}. }}} {{\hspace*{\fill}
$\Box$}\end{sloppypar}\vspace{2ex}}

\newcommand{\glob}{
\xy
(-5,0)*+{.}="1";
(5,0)*+{.}="2";
{\ar@/^1pc/^{} "1";"2"};
{\ar@/_1pc/_{} "1";"2"};
{\ar@{=>}^{} (0,2)*{};(0,-2)*{}} ;
\endxy}

\newcommand{\psinv}[6]{\xy
(#1,0)*+{#3}="x";
(#2,0)*+{#5}="y";
{\ar@<.7ex>^{#4} "x"; "y"};
{\ar@<.7ex>^{#6} "y"; "x"};
\endxy
}

\psset{unit=0.1cm,labelsep=0pt,nodesep=3pt}

\begin{document}

\title{Multivariable adjunctions and mates}

\author{Eugenia Cheng, Nick Gurski and Emily Riehl\\Department of Mathematics, University of Sheffield \\E-mail: e.cheng@sheffield.ac.uk, nick.gurski@sheffield.ac.uk\\and\\Department of Mathematics, Harvard University\\E-mail: eriehl@math.harvard.edu}

\maketitle

\begin{abstract}
We present the notion of ``cyclic double multicategory'', as a structure in which to organise multivariable adjunctions and mates.
The classic example of a 2-variable adjunction is the hom/tensor/cotensor trio of functors; we generalise this situation to $n+1$ functors of $n$ variables.  Furthermore, we generalise the mates correspondence, which enables us to pass between natural transformations involving left adjoints to those involving right adjoints.  While the standard mates correspondence is described using an isomorphism of double categories, the multivariable version requires the framework of ``double multicategories''.  Moreover, we show that the analogous isomorphisms of double multicategories give a cyclic action on the multimaps, yielding the notion of ``cyclic double multicategory''.  The work is motivated by and applied to Riehl's approach to algebraic monoidal model categories.
\end{abstract}

%\newpage

\setcounter{tocdepth}{2}
\tableofcontents

%%%%%%%%%%%%%%%%%%%%%%%%%%%%%%%%%%%%%%%%%%%%%%%%%%%%%%%%%%%%%%%%%%%%%
%%%%%%%%%%%%%%%%%%%%%%%%%%%%%%%%%%%%%%%%%%%%%%%%%%%%%%%%%%%%%%%%%%%%%

\section*{Introduction}

In this paper we introduce the notion of ``cyclic double multicategory'', as a structure in which to organise a multivariable generalisation of adjunctions and mates.

Standard adjunctions are to be thought of as having one variable; the most common example of a 2-variable adjunction is the tensor/hom/cotensor trio of functors, which has the following features:

\begin{itemize}
 \item it involves 3 functors of 2 variables,
\item each pair of functors is related by a 1-variable adjunction if we fix a variable, and
\item some care is required over dualities of source and target categories.
\end{itemize}
This example is given in detail in Section~\ref{homtensor}.

We seek to generalise this situation to $n+1$ functors of $n$ variables.  Further more, we would like to generalise the ``mates correspondence'', which enables us neatly to pass between natural transformations involving left adjoints to those involving right adjoints.

The mates correspondence is elegantly described using the framework of double categories, via a double category given as follows:
\begin{itemize}
 \item 0-cells are categories,
\item horizontal 1-cells are functors,
\item vertical 1-cells are adjunctions (pointing in a fixed chosen direction e.g. the direction of the left adjoint), and
\item 2-cells are natural transformations.
\end{itemize}
Even after we have fixed the direction of the 1-cells, there is a choice for the 2-cells---we could still take the \nts\ to live in the squares involving either the left or right adjoints.  This produces \emph{a priori} two different double categories for each choice of 1-cell direction, but the mates correspondence says precisely that there is an isomorphism of double categories between them.

For multivariable adjunctions we need our vertical 1-cells to have a string of objects as its source and a single object as its target, as for the maps in a multicategory.  (A functor of $n$ variables can be thought of as a functor of 1 variable where the variable happens to be a string of $n$ objects, but then these functors will not compose.)

However, instead of the left/right of a 1-variable adjunction, an $n$-variable adjunction has a cyclic symmetry.  For example a 2-variable adjunction involves functors
\[\begin{array}{c}
   A \times B \map{F} C\op \\
B \times C \map{G} A\op \\
C \times A \map{H} B\op.
  \end{array}\]
(Note the duality that arises as a category ``cycles'' between the source and target.)  Hence our $n$-variable adjunctions are the multimaps of a ``cyclic'' multicategory---a multicategory equipped with additional structure in the form of
\begin{itemize}
 \item an involution (such as $(\ )\op$), and
\item a cyclic action on homsets, invoking the involution appropriately.
\end{itemize}
For the $n$-variable mates correspondence we are interested in a correspondence of natural transformations such as below (for the 2-variable example):
\[\psset{unit=0.08cm,labelsep=0pt,nodesep=3pt}
\pspicture(45,22)

% a1 a2 up right
% b1 b2

\rput(0,20){\rnode{a1}{$A\times B$}} % top left
\rput(25,20){\rnode{a2}{$A' \times B'$}} % top right

\rput(0,0){\rnode{b1}{$C\op$}}   % bottom left
\rput(25,0){\rnode{b2}{${C'}\op$}}  % bottom right

\psset{nodesep=3pt,labelsep=2pt,arrows=->}
\ncline{a1}{a2}\naput{{\scriptsize $$}} % top
\ncline{b1}{b2}\nbput{{\scriptsize $$}} % bottom
\ncline{a1}{b1}\nbput{{\scriptsize $F$}} % left
\ncline{a2}{b2}\naput{{\scriptsize $F'$}} % right

\psset{labelsep=1.5pt}
\pnode(9.5,7){c1}
\pnode(15.5,13){c2}
\ncline[doubleline=true]{c1}{c2} \nbput[npos=0.4]{{\scriptsize $$}}
\endpspicture
\pspicture(45,22)

% a1 a2 up right
% b1 b2

\rput(0,20){\rnode{a1}{$B\times C$}} % top left
\rput(25,20){\rnode{a2}{$B' \times C'$}} % top right

\rput(0,0){\rnode{b1}{$A\op$}}   % bottom left
\rput(25,0){\rnode{b2}{${A'}\op$}}  % bottom right

\psset{nodesep=3pt,labelsep=2pt,arrows=->}
\ncline{a1}{a2}\naput{{\scriptsize $$}} % top
\ncline{b1}{b2}\nbput{{\scriptsize $$}} % bottom
\ncline{a1}{b1}\nbput{{\scriptsize $G$}} % left
\ncline{a2}{b2}\naput{{\scriptsize $G'$}} % right

\psset{labelsep=1.5pt}
\pnode(9.5,7){c1}
\pnode(15.5,13){c2}
\ncline[doubleline=true]{c1}{c2} \nbput[npos=0.4]{{\scriptsize $$}}
\endpspicture
\pspicture(25,22)

% a1 a2 up right
% b1 b2

\rput(0,20){\rnode{a1}{$C\times A$}} % top left
\rput(25,20){\rnode{a2}{$C' \times A'$}} % top right

\rput(0,0){\rnode{b1}{$B\op$}}   % bottom left
\rput(25,0){\rnode{b2}{${B'}\op$}}  % bottom right

\psset{nodesep=3pt,labelsep=2pt,arrows=->}
\ncline{a1}{a2}\naput{{\scriptsize $$}} % top
\ncline{b1}{b2}\nbput{{\scriptsize $$}} % bottom
\ncline{a1}{b1}\nbput{{\scriptsize $H$}} % left
\ncline{a2}{b2}\naput{{\scriptsize $H'$}} % right

\psset{labelsep=1.5pt}
\pnode(9.5,7){c1}
\pnode(15.5,13){c2}
\ncline[doubleline=true]{c1}{c2} \nbput[npos=0.4]{{\scriptsize $$}}
\endpspicture
\]
and \nts\ of this form are thus the 2-cells in our ``cyclic double multicategory''.  Recall that a double category is a category object in \Cat; similarly a cyclic double multicategory is a category object in the category of cyclic multicategories.  The cyclic action will be the $n$-variable mates correspondence.

The motivation for this work is the third author's work on algebraic monoidal model categories.  In the theory of algebraic model categories \cite{rie1} the double category framework for 1-variable adjunctions and mates plays a crucial role.  For the monoidal version \cite{rie2}, multivariable adjunctions and mates are needed.

This paper is organised as follows.  In Section~\ref{mates} we recall the standard theory of mates.  In Section~\ref{adj} we define multivariable adjunctions and the multivariable mates correspondence.  In Section~\ref{multicat} we give the definition of cyclic double multicategory, building up gradually through multicategories, cyclic multicategories and double multicategories.  We show that multivariable adjunctions form a cyclic double multicategory.  In Section~\ref{modelcat} we describe the application to algebraic monoidal model categories.

Our notion of cyclic multicategory is non-symmetric and thus generalises the notion of (non-symmetric) cyclic operad given in \cite{bb1}; symmetric cyclic operads are defined in \cite{gk1} and a multicategory version is mentioned in \cite{jk2}.  Our notion of double multicategory is not the same as the notion of fc-multicategory (introduced by Leinster in \cite{lei8} and renamed ``virtual double category'' by Cruttwell and Shulman in \cite{cs1}); fc-multicategories do not involve vertical 1-cells of higher arities.

\subsubsection*{Notation}

Throughout this paper we will write $A\opb$ for $A\op$.  Also, for $n$-variable adjunctions and cyclic multicategories, we will need to use subscripts cyclically.  Thus we will index objects by $0, \ldots, n$ with lists taken cyclically, mod $n+1$.  For example we will frequently use the string $a_{i+1}, \ldots, a_{i-1}$ which means
\[a_{i+1}, a_{i+2}, \ldots, a_n, a_0, a_1, \ldots, a_{i-1}.\]

\subsubsection*{Acknowledgements}

The third author was supported by both the NSF Graduate Research Fellowship Program and an NSF Mathematical Sciences Postdoctoral Research Fellowship.  She also had useful conversations with Anthony Elmendorf, Daniel Sch\"{a}ppi, and Dominic Verity.  This work was catalysed by a visit partially sponsored by the University of Sheffield Mathematical Sciences Research Centre.  

%%%%%%%%%%%%%%%%%%%%%%%%%%%%%%%%%%%%%%%%%%%%%%%%%%%%%%%%%%%%%%%%%%%%%
%%%%%%%%%%%%%%%%%%%%%%%%%%%%%%%%%%%%%%%%%%%%%%%%%%%%%%%%%%%%%%%%%%%%%

\section{Mates}\label{mates}

In this section we describe the situation we will be generalising.  Suppose we have the following categories, functors and adjunctions
\[
\psset{labelsep=2pt}
\pspicture(20,20)

\rput[b](0,0){
\pspicture(8,20)
\rput(3,20){\rnode{A}{$A$}}
\rput(3,0){\rnode{B}{$B$}}
\psset{nodesep=3pt,arrows=->, labelsep=2pt}
\ncline[offset=-4pt]{A}{B}\naput[labelsep=1.6pt]{$\scriptstyle\adj$}\nbput{{\scriptsize $F$}}
\ncline[offset=-5pt]{B}{A}\nbput{{\scriptsize $G$}}
\endpspicture}

\rput[b](20,0){
\pspicture(8,20)
\rput(3,20){\rnode{A}{\hh{2pt}$A'$}}
\rput(3,0){\rnode{B}{\hh{1pt}$B'$}}
\psset{nodesep=3pt,arrows=->, labelsep=2pt}
\ncline[offset=-4pt]{A}{B}\naput[labelsep=1.6pt]{$\scriptstyle\adj$}\nbput{{\scriptsize $F'$}}
\ncline[offset=-5pt]{B}{A}\nbput{{\scriptsize $G'$}}
\endpspicture}
\endpspicture
\]
with unit and counit $(\eta, \varepsilon)$ and $(\eta', \varepsilon')$ respectively.  Then given functors $S$ and $T$ and a natural transformation $\alpha$ as shown
\[\psset{unit=0.08cm,labelsep=0pt,nodesep=3pt}
\pspicture(20,22)

% tl tr down left
% bl br

\rput(0,20){\rnode{tl}{$A$}} % top left
\rput(20,20){\rnode{tr}{$A'$}} % top right

\rput(0,0){\rnode{bl}{$B$}}   % bottom left
\rput(20,0){\rnode{br}{$B'$}}  % bottom right

\psset{nodesep=3pt,labelsep=2pt,arrows=->}
\ncline{tl}{tr}\naput{{\scriptsize $S$}} % top
\ncline{bl}{br}\nbput{{\scriptsize $T$}} % bottom
\ncline{tl}{bl}\nbput{{\scriptsize $F$}} % left
\ncline{tr}{br}\naput{{\scriptsize $F'$}} % right

\psset{labelsep=1.5pt}

\pnode(13,13){a3}
\pnode(7,7){b3}
\ncline[doubleline=true]{a3}{b3} \nbput[npos=0.4]{{\scriptsize $\alpha$}}

\endpspicture\]
its \demph{mate} $\bar{\alpha}$ is the \nt\
\[\psset{unit=0.08cm,labelsep=0pt,nodesep=3pt}
\pspicture(20,22)

% tl tr down right
% bl br

\rput(0,20){\rnode{tl}{$A$}} % top left
\rput(20,20){\rnode{tr}{$A'$}} % top right

\rput(0,0){\rnode{bl}{$B$}}   % bottom left
\rput(20,0){\rnode{br}{$B'$}}  % bottom right

\psset{nodesep=3pt,labelsep=2pt,arrows=->}
\ncline{tl}{tr}\naput{{\scriptsize $S$}} % top
\ncline{bl}{br}\nbput{{\scriptsize $T$}} % bottom
\ncline{bl}{tl}\naput{{\scriptsize $G$}} % left
\ncline{br}{tr}\nbput{{\scriptsize $G'$}} % right

\psset{labelsep=1.5pt}

\pnode(7,13){a3}
\pnode(13,7){b3}
\ncline[doubleline=true]{a3}{b3} \nbput[npos=0.4]{{\scriptsize $\bar{\alpha}$}}

\endpspicture\]
obtained as the following composite

\[\psset{unit=0.08cm,labelsep=1pt,nodesep=1pt}
\pspicture(-20,-10)(40,22)

% a0 a1 a2
%    b1 b2 b3

\rput(-20,20){\rnode{a0}{$B$}} % 
\rput(0,20){\rnode{a1}{$A$}} % top left
\rput(20,20){\rnode{a2}{$A'$}} % top right

\rput(0,0){\rnode{b1}{$B$}}   % bottom left
\rput(20,0){\rnode{b2}{$B'$}}  % bottom right
\rput(40,0){\rnode{b3}{$A'.$}} % 

\psset{nodesep=2pt,labelsep=2pt,arrows=->}
\ncline{a0}{a1}\naput{{\scriptsize $G$}} % top
\ncline{a0}{b1}\nbput{{\scriptsize $1$}} % top

\ncline{a1}{a2}\naput{{\scriptsize $S$}} % top
\ncline{b1}{b2}\nbput{{\scriptsize $T$}} % bottom
\ncline{a1}{b1}\nbput[labelsep=1pt]{{\scriptsize $F$}} % left
\ncline{a2}{b2}\naput[labelsep=1pt]{{\scriptsize $F'$}} % right

\ncline{a2}{b3}\naput{{\scriptsize $1$}} % top
\ncline{b2}{b3}\nbput{{\scriptsize $G'$}} % top

\psset{doubleline=true,labelsep=0pt, nodesep=3pt}
\pnode(13,13){a3}
\pnode(7,7){b3}
\ncline{a3}{b3} \nbput[npos=0.4]{{\scriptsize $\alpha$}}

\rput(-9,11)
{\pcline(5,5)(0,0) \nbput[npos=0.4]{{\scriptsize $\varepsilon$}}}
\rput(24,3)
{\pcline(5,5)(0,0)\naput[npos=0.4,labelsep=-1pt]{{\scriptsize $\eta'$}}}

\endpspicture\]
Conversely we can start with 
\[\psset{unit=0.08cm,labelsep=0pt,nodesep=3pt}
\pspicture(20,22)

% tl tr down right
% bl br

\rput(0,20){\rnode{tl}{$A$}} % top left
\rput(20,20){\rnode{tr}{$A'$}} % top right

\rput(0,0){\rnode{bl}{$B$}}   % bottom left
\rput(20,0){\rnode{br}{$B'$}}  % bottom right

\psset{nodesep=3pt,labelsep=2pt,arrows=->}
\ncline{tl}{tr}\naput{{\scriptsize $S$}} % top
\ncline{bl}{br}\nbput{{\scriptsize $T$}} % bottom
\ncline{bl}{tl}\naput{{\scriptsize $G$}} % left
\ncline{br}{tr}\nbput{{\scriptsize $G'$}} % right

\psset{labelsep=1.5pt}

\pnode(7,13){a3}
\pnode(13,7){b3}
\ncline[doubleline=true]{a3}{b3} \nbput[npos=0.4]{{\scriptsize $\beta$}}

\endpspicture\]
and obtain the mate
\[\psset{unit=0.08cm,labelsep=0pt,nodesep=3pt}
\pspicture(20,22)

% tl tr down left
% bl br

\rput(0,20){\rnode{tl}{$A$}} % top left
\rput(20,20){\rnode{tr}{$A'$}} % top right

\rput(0,0){\rnode{bl}{$B$}}   % bottom left
\rput(20,0){\rnode{br}{$B'$}}  % bottom right

\psset{nodesep=3pt,labelsep=2pt,arrows=->}
\ncline{tl}{tr}\naput{{\scriptsize $S$}} % top
\ncline{bl}{br}\nbput{{\scriptsize $T$}} % bottom
\ncline{tl}{bl}\nbput{{\scriptsize $F$}} % left
\ncline{tr}{br}\naput{{\scriptsize $F'$}} % right

\psset{labelsep=1.5pt}

\pnode(13,13){a3}
\pnode(7,7){b3}
\ncline[doubleline=true]{a3}{b3} \naput[labelsep=0pt,npos=0.4]{{\scriptsize $\bar{\beta}$}}

\endpspicture\]
as the composite 

\[\psset{unit=0.08cm,labelsep=1pt,nodesep=1pt}
\pspicture(-20,-5)(40,22)

%     a1 a2 b3
%  a0 b1 b2 

\rput(-20,0){\rnode{a0}{$A$}} % 
\rput(0,20){\rnode{a1}{$A$}} % top left
\rput(20,20){\rnode{a2}{$A'$}} % top right

\rput(0,0){\rnode{b1}{$B$}}   % bottom left
\rput(20,0){\rnode{b2}{$B'.$}}  % bottom right
\rput(40,20){\rnode{b3}{$B'$}} % top left

\psset{nodesep=2pt,labelsep=2pt,arrows=->}
\ncline{a0}{a1}\naput{{\scriptsize $1$}} % top
\ncline{a0}{b1}\nbput{{\scriptsize $F$}} % top

\ncline{a1}{a2}\naput{{\scriptsize $S$}} % top
\ncline{b1}{b2}\nbput{{\scriptsize $T$}} % bottom
\ncline{b1}{a1}\naput[labelsep=1pt]{{\scriptsize $G$}} % left
\ncline{b2}{a2}\nbput[labelsep=1pt]{{\scriptsize $G'$}} % right

\ncline{a2}{b3}\naput{{\scriptsize $F'$}} % top
\ncline{b2}{b3}\nbput{{\scriptsize $1$}} % top

\psset{doubleline=true,labelsep=0pt, nodesep=3pt}
\pnode(7,13){a3}
\pnode(13,7){b3}
\ncline{a3}{b3} \nbput[npos=0.4]{{\scriptsize $\beta$}}

\rput(-9,3)
{\pcline(0,5)(5,0) \nbput[npos=0.4]{{\scriptsize $\varepsilon$}}}
\rput(24,11)
{\pcline(0,5)(5,0)\naput[npos=0.4,labelsep=-1pt]{{\scriptsize $\eta'$}}}

\endpspicture\]

\noi By triangle identities these processes of ``conjugation'' are inverse to one another.  Furthermore, the correspondence respects both horizontal and vertical composition in the following sense.  Given adjunctions
\[
\psset{labelsep=2pt}
\pspicture(40,20)

\rput[b](0,0){
\pspicture(8,20)
\rput(3,18){\rnode{A}{$A_1$}}
\rput(3,0){\rnode{B}{$B_1$}}
\psset{nodesep=3pt,arrows=->, labelsep=2pt}
\ncline[offset=-4pt]{A}{B}\naput[labelsep=1.6pt]{$\scriptstyle\adj$}\nbput{{\scriptsize $F_1$}}
\ncline[offset=-5pt]{B}{A}\nbput{{\scriptsize $G_1$}}
\endpspicture}

\rput[b](20,0){
\pspicture(8,20)
\rput(3,18){\rnode{A}{$A_2$}}
\rput(3,0){\rnode{B}{$B_2$}}
\psset{nodesep=3pt,arrows=->,labelsep=2pt}
\ncline[offset=-4pt]{A}{B}\naput[labelsep=1.6pt]{$\scriptstyle\adj$}\nbput{{\scriptsize $F_2$}}
\ncline[offset=-5pt]{B}{A}\nbput{{\scriptsize $G_2$}}
\endpspicture}

\rput[b](40,0){
\pspicture(8,20)
\rput(3,18){\rnode{A}{$A_3$}}
\rput(3,0){\rnode{B}{$B_3$}}
\psset{nodesep=3pt,arrows=->, labelsep=2pt}
\ncline[offset=-4pt]{A}{B}\naput[labelsep=1.6pt]{$\scriptstyle\adj$}\nbput{{\scriptsize $F_3$}}
\ncline[offset=-5pt]{B}{A}\nbput{{\scriptsize $G_3$}}
\endpspicture}

\endpspicture
\]
and \nts\
\[\psset{unit=0.08cm,labelsep=0pt,nodesep=3pt}
\pspicture(20,22)

% tl tr
% bl br

\rput(0,20){\rnode{tl}{$A_1$}} % top left
\rput(20,20){\rnode{tr}{$A_2$}} % top mid
\rput(40,20){\rnode{trr}{$A_3$}} % top right

\rput(0,0){\rnode{bl}{$B_1$}}   % bottom left
\rput(20,0){\rnode{br}{$B_2$}}  % bottom mid
\rput(40,0){\rnode{brr}{$B_3$}}  % bottom right

\psset{nodesep=3pt,labelsep=2pt,arrows=->}
\ncline{tl}{tr}\naput{{\scriptsize $S_1$}} % top
\ncline{bl}{br}\nbput{{\scriptsize $T_1$}} % bottom
\ncline{tr}{trr}\naput{{\scriptsize $S_2$}} % top
\ncline{br}{brr}\nbput{{\scriptsize $T_2$}} % bottom

\ncline{tl}{bl}\nbput{{\scriptsize $F_1$}} % left
\ncline{tr}{br}\naput[labelsep=1pt]{{\scriptsize $F_2$}} % right
\ncline{trr}{brr}\naput{{\scriptsize $F_3$}} % right

\psset{labelsep=1.5pt}

\pnode(13,13){a3}
\pnode(7,7){b3}
\ncline[doubleline=true]{a3}{b3} \naput[npos=0.4,labelsep=0pt]{{\scriptsize $\alpha_1$}}

\pnode(33,13){a3}
\pnode(27,7){b3}
\ncline[doubleline=true]{a3}{b3} \naput[npos=0.4,labelsep=0pt]{{\scriptsize $\alpha_2$}}

\endpspicture\]
we have
\[\overline{\alpha_2 \ast \alpha_1} = \overline{\alpha_2} \ast \overline{\alpha_1}\]
where $\ast$ is to be interpreted with the appropriate whiskering, so in fact the honest equality is
\[\overline{T_2 \alpha_1 \circ \alpha_2S_1} = \overline{\alpha_2}T_1 \circ S_2 \overline{\alpha_1}.\]
For ``vertical'' composition, given adjunctions
\[
\psset{labelsep=2pt}
\pspicture(40,40)

\rput[b](0,0){
\pspicture(8,20)
\rput(3,40){\rnode{A}{$A_1$}}
\rput(3,20){\rnode{B}{$B_1$}}
\rput(3, 0){\rnode{C}{$C_1$}}

\psset{nodesep=3pt,arrows=->, labelsep=2pt}

\ncline[offset=-4pt]{A}{B}\naput[labelsep=1.6pt] {$\scriptstyle\adj$}\nbput{{\scriptsize $F_1$}} 
\ncline[offset=-5pt]{B}{A}\nbput{{\scriptsize $G_1$}}

\ncline[offset=-4pt]{B}{C}\naput[labelsep=1.6pt] {$\scriptstyle\adj$}\nbput{{\scriptsize $H_1$}} 
\ncline[offset=-5pt]{C}{B}\nbput{{\scriptsize $K_1$}}

\endpspicture}

\rput[b](20,0){
\pspicture(8,20)
\rput(3,40){\rnode{A}{$A_2$}}
\rput(3,20){\rnode{B}{$B_2$}}
\rput(3, 0){\rnode{C}{$C_2$}}

\psset{nodesep=3pt,arrows=->, labelsep=2pt}

\ncline[offset=-4pt]{A}{B}\naput[labelsep=1.6pt] {$\scriptstyle\adj$}\nbput{{\scriptsize $F_2$}} 
\ncline[offset=-5pt]{B}{A}\nbput{{\scriptsize $G_2$}}

\ncline[offset=-4pt]{B}{C}\naput[labelsep=1.6pt] {$\scriptstyle\adj$}\nbput{{\scriptsize $H_2$}} 
\ncline[offset=-5pt]{C}{B}\nbput{{\scriptsize $K_2$}}

\endpspicture}

\endpspicture
\]
and \nts\
\[\psset{unit=0.08cm,labelsep=0pt,nodesep=3pt}
\pspicture(20,42)

% a1 a2
% b1 b2
% c1 c2

\rput(0,40){\rnode{a1}{$A_1$}} % top left
\rput(20,40){\rnode{a2}{$A_2$}} % top mid

\rput(0,20){\rnode{b1}{$B_1$}}   % bottom left
\rput(20,20){\rnode{b2}{$B_2$}}  % bottom mid

\rput(0,0){\rnode{c1}{$C_1$}}   % bottom left
\rput(20,0){\rnode{c2}{$C_2$}}  % bottom mid

\psset{nodesep=3pt,labelsep=2pt,arrows=->}
\ncline{a1}{a2}\naput{{\scriptsize $S$}} % top
\ncline{b1}{b2}\nbput{{\scriptsize $T$}} % bottom
\ncline{c1}{c2}\nbput{{\scriptsize $U$}} % bottom

\ncline{a1}{b1}\nbput{{\scriptsize $F_1$}} % left
\ncline{b1}{c1}\nbput{{\scriptsize $H_1$}} % left

\ncline{a2}{b2}\naput{{\scriptsize $F_2$}} % right
\ncline{b2}{c2}\naput{{\scriptsize $H_2$}} % right

\psset{labelsep=1.5pt}

\pnode(13,33){a3}
\pnode(7,27){b3}
\ncline[doubleline=true]{a3}{b3} \naput[npos=0.4,labelsep=0pt]{{\scriptsize $\alpha_1$}}

\pnode(13,13){a3}
\pnode(7,7){b3}
\ncline[doubleline=true]{a3}{b3} \naput[npos=0.4,labelsep=0pt]{{\scriptsize $\alpha_2$}}

\endpspicture\]
we have 
\[\overline{\alpha_2 \circ \alpha_1} = \overline{\alpha_2} \circ \overline{\alpha_1}\]
which actually means
\[\overline{\alpha_2 F_1 \circ H_2 \alpha_1} = G_2 \overline{\alpha_2} \circ \overline{\alpha_1} K_1.\]
Both of these facts are easily checked using 2-pasting diagrams and triangle identities.  

This situation is conveniently formalised using double categories.  In the following definition we have chosen the direction of the vertical 1-cells to correspond to the direction of the left adjoints.  

\begin{mydefinition}
 We define two double categories \Ladj\ and \Radj\ with the same 0- and 1-cells, but different 2-cells.  In both cases the 0-cells are categories, the horizontal 1-cells are functors, and a vertical 1-cell $A \tra B$ is an adjunction
\[
\psset{labelsep=2pt}
\pspicture(20,8)
\rput(0,3){\rnode{A}{$A$}}
\rput(20,3){\rnode{B}{$B.$}}
\psset{nodesep=2pt,arrows=->, labelsep=2pt}
\ncline[offset=5pt]{A}{B}\nbput{$\scriptstyle\perp$}\naput{{\scriptsize $F$}}
\ncline[offset=5pt]{B}{A}\naput{{\scriptsize $G$}}
\endpspicture
\]
A 2-cell 
\[\psset{unit=0.1cm,labelsep=0pt,nodesep=3pt}
\pspicture(30,20)

% a1 a2
% b1 b2

\rput(0,20){\rnode{a1}{$A_1$}} % top left
\rput(20,20){\rnode{a2}{$A_2$}} % top right

\rput(0,0){\rnode{b1}{$B_1$}}   % bottom left
\rput(20,0){\rnode{b2}{$B_2$}}  % bottom right

\psset{nodesep=3pt,labelsep=2pt,arrows=->}
\ncline{a1}{a2}\naput{{\scriptsize $S$}} % top
\ncline{b1}{b2}\nbput{{\scriptsize $T$}} % bottom

\psset{nodesep=3pt,arrows=->, labelsep=2pt}
\ncline[offset=-4pt]{a1}{b1}
  \naput[labelsep=1.6pt]{$\scriptstyle\adj$}
  \nbput{{\scriptsize $F_1$}}
\ncline[offset=-5pt]{b1}{a1}\nbput{{\scriptsize $G_1$}}

\ncline[offset=-4pt]{a2}{b2}
  \naput[labelsep=1.6pt]{$\scriptstyle\adj$}
  \nbput{{\scriptsize $F_2$}}
\ncline[offset=-5pt]{b2}{a2}\nbput{{\scriptsize $G_2$}}

\pnode(10,14){a3}
\pnode(10,6){b3}
\ncline[doubleline=true]{a3}{b3}

\endpspicture\]
is given in each case as follows.

\begin{itemize}
 \item In \Ladj\ such a 2-cell is a \nt\ 
\[\psset{unit=0.08cm,labelsep=0pt,nodesep=3pt}
\pspicture(20,22)

% tl tr down left
% bl br

\rput(0,20){\rnode{tl}{$A_1$}} % top left
\rput(20,20){\rnode{tr}{$A_2$}} % top right

\rput(0,0){\rnode{bl}{$B_1$}}   % bottom left
\rput(20,0){\rnode{br}{$B_2.$}}  % bottom right

\psset{nodesep=3pt,labelsep=2pt,arrows=->}
\ncline{tl}{tr}\naput{{\scriptsize $S$}} % top
\ncline{bl}{br}\nbput{{\scriptsize $T$}} % bottom
\ncline{tl}{bl}\nbput{{\scriptsize $F_1$}} % left
\ncline{tr}{br}\naput{{\scriptsize $F_2$}} % right

\psset{labelsep=1.5pt}

\pnode(13,13){a3}
\pnode(7,7){b3}
\ncline[doubleline=true]{a3}{b3} \nbput[npos=0.4]{{\scriptsize $\alpha$}}

\endpspicture\]

 \item In \Radj\ such a 2-cell is a \nt\ 
\[\psset{unit=0.08cm,labelsep=0pt,nodesep=3pt}
\pspicture(20,22)

% tl tr down right
% bl br

\rput(0,20){\rnode{tl}{$A_1$}} % top left
\rput(20,20){\rnode{tr}{$A_2$}} % top right

\rput(0,0){\rnode{bl}{$B_1$}}   % bottom left
\rput(20,0){\rnode{br}{$B_2.$}}  % bottom right

\psset{nodesep=3pt,labelsep=2pt,arrows=->}
\ncline{tl}{tr}\naput{{\scriptsize $S$}} % top
\ncline{bl}{br}\nbput{{\scriptsize $T$}} % bottom
\ncline{bl}{tl}\naput{{\scriptsize $G_1$}} % left
\ncline{br}{tr}\nbput{{\scriptsize $G_2$}} % right

\psset{labelsep=1.5pt}

\pnode(7,13){a3}
\pnode(13,7){b3}
\ncline[doubleline=true]{a3}{b3} \nbput[npos=0.4]{{\scriptsize $\alpha$}}

\endpspicture\]

\end{itemize}

\end{mydefinition}

\begin{theorem}{\bfseries\upshape \cite[Proposition 2.2]{ks1}}
 
There is an isomorphism of double categories
\[ \Ladj \iso \Radj\]
which is the identity on 0- and 1-cells (horizontal and vertical); on 2-cells it is given by taking mates.
\end{theorem}

We now look at this from a slightly different point of view that seems a little contrived here, but leads to a natural framework for the $n$-variable generalisation.  The idea is to notice that an adjunction
\[
\psset{labelsep=2pt}
\pspicture(20,8)
\rput(0,3){\rnode{A}{$A$}}
\rput(20,3){\rnode{B}{$B$}}
\psset{nodesep=2pt,arrows=->, labelsep=2pt}
\ncline[offset=5pt]{A}{B}\nbput{$\scriptstyle\perp$}\naput{{\scriptsize $F$}}
\ncline[offset=5pt]{B}{A}\naput{{\scriptsize $G$}}
\endpspicture
\]
is equivalently an adjunction
\[
\psset{labelsep=2pt}
\pspicture(20,8)
\rput(0,3){\rnode{A}{$B\opb$}}
\rput(20,3){\rnode{B}{$A\opb.$}}
\psset{nodesep=2pt,arrows=->, labelsep=2pt}
\ncline[offset=5pt]{A}{B}\nbput{$\scriptstyle\perp$}\naput{{\scriptsize $G\opb$}}
\ncline[offset=5pt]{B}{A}\naput{{\scriptsize $F\opb$}}
\endpspicture
\]
Now, we could deal with this by introducing yet another pair of double categories $\Ladj_R$ and $\Radj_R$ as above but whose vertical 1-cells point in the direction of the \emph{right} adjoints; the 2-cell directions must also be changed accordingly.  We would then get isomorphisms of double categories
\[\begin{array}{c}
   (\ )\opb \: \Ladj \lra \Radj_R \\
   (\ )\opb \: \Radj \lra \Ladj_R. 
  \end{array}\]
However, we can actually express all this structure using one single version of the above four isomorphic double categories, as follows.

Given a 2-cell in \Ladj, that is, a \nt\
\[\psset{unit=0.08cm,labelsep=0pt,nodesep=3pt}
\pspicture(20,22)

% tl tr down left
% bl br

\rput(0,20){\rnode{tl}{$A_1$}} % top left
\rput(20,20){\rnode{tr}{$A_2$}} % top right

\rput(0,0){\rnode{bl}{$B_1$}}   % bottom left
\rput(20,0){\rnode{br}{$B_2$}}  % bottom right

\psset{nodesep=3pt,labelsep=2pt,arrows=->}
\ncline{tl}{tr}\naput{{\scriptsize $S$}} % top
\ncline{bl}{br}\nbput{{\scriptsize $T$}} % bottom
\ncline{tl}{bl}\nbput{{\scriptsize $F_1$}} % left
\ncline{tr}{br}\naput{{\scriptsize $F_2$}} % right

\psset{labelsep=1.5pt}

\pnode(13,13){a3}
\pnode(7,7){b3}
\ncline[doubleline=true]{a3}{b3} \nbput[npos=0.4]{{\scriptsize $\alpha$}}

\endpspicture\]
its mate
\[\psset{unit=0.08cm,labelsep=0pt,nodesep=3pt}
\pspicture(20,22)

% tl tr down right
% bl br

\rput(0,20){\rnode{tl}{$A_1$}} % top left
\rput(20,20){\rnode{tr}{$A_2$}} % top right

\rput(0,0){\rnode{bl}{$B_1$}}   % bottom left
\rput(20,0){\rnode{br}{$B_2$}}  % bottom right

\psset{nodesep=3pt,labelsep=2pt,arrows=->}
\ncline{tl}{tr}\naput{{\scriptsize $S$}} % top
\ncline{bl}{br}\nbput{{\scriptsize $T$}} % bottom
\ncline{bl}{tl}\naput{{\scriptsize $G_1$}} % left
\ncline{br}{tr}\nbput{{\scriptsize $G_2$}} % right

\psset{labelsep=1.5pt}

\pnode(7,13){a3}
\pnode(13,7){b3}
\ncline[doubleline=true]{a3}{b3} \nbput[npos=0.4]{{\scriptsize $\bar{\alpha}$}}

\endpspicture\]
is not \emph{a priori} a 2-cell of \Ladj\ as its source and target involve \emph{right} adjoints $G_1$ and $G_2$.  However, it can be dualised to give
\[\psset{unit=0.08cm,labelsep=0pt,nodesep=3pt}
\pspicture(20,22)

% tl tr up left
% bl br

\rput(0,20){\rnode{tl}{$A_1\opb$}} % top left
\rput(20,20){\rnode{tr}{$A_2\opb$}} % top right

\rput(0,0){\rnode{bl}{$B_1\opb$}}   % bottom left
\rput(20,0){\rnode{br}{$B_2\opb$}}  % bottom right

\psset{nodesep=3pt,labelsep=2pt,arrows=->}
\ncline{tl}{tr}\naput{{\scriptsize $S\opb$}} % top
\ncline{bl}{br}\nbput{{\scriptsize $T\opb$}} % bottom
\ncline{bl}{tl}\naput{{\scriptsize $G_1\opb$}} % left
\ncline{br}{tr}\nbput{{\scriptsize $G_2\opb$}} % right

\psset{labelsep=1.5pt}

\pnode(13,7){a3}
\pnode(7,13){b3}
\ncline[doubleline=true]{a3}{b3} \nbput[npos=0.4]{{\scriptsize $\bar{\alpha}\opb$}}

\endpspicture\]
where we must reverse the 2-cell direction as the target category has been dualised.   Thus, turning the diagram round so that the left adjoints point downwards, we have 

\[\psset{unit=0.08cm,labelsep=0pt,nodesep=3pt}
\pspicture(0,-6)(20,24)

\rput(0,0){\rnode{tl}{$A_1\opb$}} % top left
\rput(20,00){\rnode{tr}{$A_2\opb$}} % top right

\rput(0,20){\rnode{bl}{$B_1\opb$}}   % bottom left
\rput(20,20){\rnode{br}{$B_2\opb$}}  % bottom right

\psset{nodesep=3pt,labelsep=2pt,arrows=->}
\ncline{tl}{tr}\nbput{{\scriptsize $S\opb$}} % top
\ncline{bl}{br}\naput{{\scriptsize $T\opb$}} % bottom
\ncline{bl}{tl}\nbput{{\scriptsize $G_1\opb$}} % left
\ncline{br}{tr}\naput{{\scriptsize $G_2\opb$}} % right

\psset{labelsep=0pt}

\pnode(13,13){a3}
\pnode(7,7){b3}
\ncline[doubleline=true]{a3}{b3} \nbput[npos=0.4,labelsep=-1pt]{{\scriptsize $\bar{\alpha}\opb$}}

\endpspicture\]
which \emph{is} a 2-cell of \Ladj\ as $G_1\opb$ and $G_2\opb$ are \emph{left} adjoints.  Thus the mates correspondence actually gives us some extra structure on \Ladj\ in the form of isomorphisms:
\begin{itemize}
 \item $\Ladj_v(A,B) \iso \Ladj_v(B\opb, A\opb)$, and
\item $\Ladj_2(S,T) \iso \Ladj_2(T\opb, S\opb)$.
\end{itemize}
Here the subscript $v$ indicates the hom-set of vertical 1-cells, and the subscript $2$ indicates the hom-set of 2-cells with respect to their horizontal 1-cell boundaries.

We will see that these isomorphisms are the beginning of a cyclic structure: the 1-ary part.  The situation has a slightly different flavour, without technically being different, if we put it in the language of ``mutual left adjoints''.

\begin{mydefinition}
 Consider functors
\[\begin{array}{ccc}
   A \tmap{F} B\opb & \mbox{so} & A \opb \tmap{F\opb} B\\
B \tmap{G} A\opb & & B\opb \tmap{G\opb} A.
  \end{array}\]
A \demph{mutual left adjunction} of $F$ and $G$ is an adjunction
\[F\opb \ladj G\]
or equivalently
\[G\opb \ladj F.\]
Note that this is given by isomorphisms
\[B(Fa, b) \iso A(Gb,a)\]
natural in $a$ and $b$.  If we started with
\[\begin{array}{cl}
   A\opb \tmap{F} B & \mbox{and}\\
B\opb \tmap{G} A
  \end{array}\]
then the adjunctions
\[ F\opb \ladj G \mbox{\ \ or \ \ } G\opb \ladj F\]
as above would be given by isomorphisms
\[B(b,Fa) \iso A(a,Gb),\]
which is called a \demph{mutual right adjoint}.  

Note that the unit and counit for a mutual left adjoint as above have components
\[\eta_a  \:  GFa \tra a \in A \ \makebox[0pt][l]{\ and}\]
\[\varepsilon_b \:  FGb \tra b \in B,\]
whereas for a mutual right adjoint the components are
\[\eta_a  \:  a \tra GFa \in A \ \makebox[0pt][l]{\ and}\]
\[\varepsilon_b \:  b \tra FGb \in B.\]

\end{mydefinition}

\begin{myremark}\label{epsilonremark}
The unit and counit given above are for an adjunction
\[ F\opb \ladj G\]
whereas for the (equivalent) adjunction
\[ G\opb \ladj F\]
the unit and counit are the other way round, that is,
\[\ve_a  \:  GFa \tra a \in A \ \makebox[0pt][l]{\ and}\]
\[\eta_b \:  FGb \tra b \in B.\]
In the spirit of symmetry, we will refer to both of these natural transformations as $\ve$; it will be clear from the source and target for which adjunction this is actually a counit.  

\end{myremark}

\begin{myremark}
Note that some ambiguity can arise when we are dealing with multivariable functors that are covariant in some variables and contravariant in others.  In this case we will be careful to refer to the \emph{actual} adjunctions or specify the hom-set isomorphisms. 
\end{myremark}

We can now express the mates correspondence for mutual left adjunctions.  Given a mutual left adjunction between
\[\begin{array}{cl}
   A \tmap{F} B\opb & \mbox{and}\\
B \tmap{G} A\opb
  \end{array}\]
the mates correspondence together with duality as above gives us a correspondence between \nts\

\[\psset{unit=0.08cm,labelsep=0pt,nodesep=3pt}
\pspicture(0,0)(20,25)

% a1 a2 up right
% b1 b2

\rput(0,20){\rnode{a1}{$A$}} % top left
\rput(20,20){\rnode{a2}{$A'$}} % top right

\rput(0,0){\rnode{b1}{$B\opb$}}   % bottom left
\rput(20,0){\rnode{b2}{$B\opb$}}  % bottom right

\psset{nodesep=3pt,labelsep=2pt,arrows=->}
\ncline{a1}{a2}\naput{{\scriptsize $S$}} % top
\ncline{b1}{b2}\nbput{{\scriptsize $T\opb$}} % bottom
\ncline{a1}{b1}\nbput{{\scriptsize $F$}} % left
\ncline{a2}{b2}\naput{{\scriptsize $F'$}} % right

\psset{labelsep=1.5pt}
\pnode(7,7){c1}
\pnode(13,13){c2}
\ncline[doubleline=true]{c1}{c2} \nbput[npos=0.4]{{\scriptsize $$}}

\rput(33,10){and}

\endpspicture
\hh{3em}\hh{3em}
\pspicture(0,0)(20,22)

% a1 a2 up right
% b1 b2

\rput(0,20){\rnode{a1}{$B$}} % top left
\rput(20,20){\rnode{a2}{${B'}\opb$}} % top right

\rput(0,0){\rnode{b1}{$A\opb$}}   % bottom left
\rput(20,0){\rnode{b2}{${A'}\opb.$}}  % bottom right

\psset{nodesep=3pt,labelsep=2pt,arrows=->}
\ncline{a1}{a2}\naput{{\scriptsize $T$}} % top
\ncline{b1}{b2}\nbput{{\scriptsize $S\opb$}} % bottom
\ncline{a1}{b1}\nbput{{\scriptsize $G$}} % left
\ncline{a2}{b2}\naput{{\scriptsize $G'$}} % right

\psset{labelsep=1.5pt}
\pnode(7,7){c1}
\pnode(13,13){c2}
\ncline[doubleline=true]{c1}{c2} \nbput[npos=0.4]{{\scriptsize $$}}

\endpspicture
\]

\vspace{2em}\noi This is obtained from the ordinary mates correspondence by taking some appropriate duals. This is the $n=1$ part of the $n$-variable case, in which we look at \nts\

\[\psset{unit=0.08cm,labelsep=0pt,nodesep=3pt}
\pspicture(0,-5)(70,22)

% a1 a2 up right
% b1 b2

\rput(0,20){\rnode{a1}{$A_1 \times \cdots \times A_n$}} % top left
\rput(70,20){\rnode{a2}{$A_1' \times \cdots \times A_n'$}} % top right

\rput(0,0){\rnode{b1}{$A_{0}\opb$}}   % bottom left
\rput(70,0){\rnode{b2}{${A_{0}'}\opb$}}  % bottom right

\psset{nodesep=3pt,labelsep=2pt,arrows=->}
\ncline{a1}{a2}\naput{{\scriptsize $S_1 \times \cdots \times S_n$}} % top
\ncline{b1}{b2}\nbput{{\scriptsize $S_{0}\opb$}} % bottom
\ncline{a1}{b1}\nbput{{\scriptsize $$}} % left
\ncline{a2}{b2}\naput{{\scriptsize $$}} % right

\psset{labelsep=1.5pt}
\pnode(32,7){c1}
\pnode(38,13){c2}
\ncline[doubleline=true]{c1}{c2} \nbput[npos=0.4]{{\scriptsize $$}}

\endpspicture\]
and
\[\psset{unit=0.08cm,labelsep=0pt,nodesep=3pt}
\pspicture(0,-5)(70,24)

% a1 a2 up right
% b1 b2

\rput(0,20){\rnode{a1}{$A_2 \times \cdots \times A_{0}$}} % top left
\rput(70,20){\rnode{a2}{$A_2' \times \cdots \times A_{0}'$}} % top right

\rput(0,0){\rnode{b1}{$A_1\opb$}}   % bottom left
\rput(70,0){\rnode{b2}{${A_1'}\opb$}}  % bottom right

\psset{nodesep=3pt,labelsep=2pt,arrows=->}
\ncline{a1}{a2}\naput{{\scriptsize $S_2 \times \cdots \times S_{0}$}} % top
\ncline{b1}{b2}\nbput{{\scriptsize $S_1\opb$}} % bottom
\ncline{a1}{b1}\nbput{{\scriptsize $$}} % left
\ncline{a2}{b2}\naput{{\scriptsize $$}} % right

\psset{labelsep=1.5pt}
\pnode(32,7){c1}
\pnode(38,13){c2}
\ncline[doubleline=true]{c1}{c2} \nbput[npos=0.4]{{\scriptsize $$}}

\endpspicture
\]
and every cyclic variant.

%%%%%%%%%%%%%%%%%%%%%%%%%%%%%%%%%%%%%%%%%%%%%%%%%%%%%%%%%%%%%%%%%%%%%
%%%%%%%%%%%%%%%%%%%%%%%%%%%%%%%%%%%%%%%%%%%%%%%%%%%%%%%%%%%%%%%%%%%%%
%%%%%%%%%%%%%%%%%%%%%%%%%%%%%%%%%%%%%%%%%%%%%%%%%%%%%%%%%%%%%%%%%%%%%
%%%%%%%%%%%%%%%%%%%%%%%%%%%%%%%%%%%%%%%%%%%%%%%%%%%%%%%%%%%%%%%%%%%%%
%%%%%%%%%%%%%%%%%%%%%%%%%%%%%%%%%%%%%%%%%%%%%%%%%%%%%%%%%%%%%%%%%%%%%
%%%%%%%%%%%%%%%%%%%%%%%%%%%%%%%%%%%%%%%%%%%%%%%%%%%%%%%%%%%%%%%%%%%%%
%%%%%%%%%%%%%%%%%%%%%%%%%%%%%%%%%%%%%%%%%%%%%%%%%%%%%%%%%%%%%%%%%%%%%
%%%%%%%%%%%%%%%%%%%%%%%%%%%%%%%%%%%%%%%%%%%%%%%%%%%%%%%%%%%%%%%%%%%%%
%%%%%%%%%%%%%%%%%%%%%%%%%%%%%%%%%%%%%%%%%%%%%%%%%%%%%%%%%%%%%%%%%%%%%
%%%%%%%%%%%%%%%%%%%%%%%%%%%%%%%%%%%%%%%%%%%%%%%%%%%%%%%%%%%%%%%%%%%%%
%%%%%%%%%%%%%%%%%%%%%%%%%%%%%%%%%%%%%%%%%%%%%%%%%%%%%%%%%%%%%%%%%%%%%
%%%%%%%%%%%%%%%%%%%%%%%%%%%%%%%%%%%%%%%%%%%%%%%%%%%%%%%%%%%%%%%%%%%%%
%%%%%%%%%%%%%%%%%%%%%%%%%%%%%%%%%%%%%%%%%%%%%%%%%%%%%%%%%%%%%%%%%%%%%
%%%%%%%%%%%%%%%%%%%%%%%%%%%%%%%%%%%%%%%%%%%%%%%%%%%%%%%%%%%%%%%%%%%%%
%%%%%%%%%%%%%%%%%%%%%%%%%%%%%%%%%%%%%%%%%%%%%%%%%%%%%%%%%%%%%%%%%%%%%
%%%%%%%%%%%%%%%%%%%%%%%%%%%%%%%%%%%%%%%%%%%%%%%%%%%%%%%%%%%%%%%%%%%%%
%%%%%%%%%%%%%%%%%%%%%%%%%%%%%%%%%%%%%%%%%%%%%%%%%%%%%%%%%%%%%%%%%%%%%
%%%%%%%%%%%%%%%%%%%%%%%%%%%%%%%%%%%%%%%%%%%%%%%%%%%%%%%%%%%%%%%%%%%%%

%%%%%%%%%%%%%%%%%%%%%%%%%%%%%%%%%%%%%%%%%%%%%%%%%%%%%%%%%%%%%%%%%%%%%
%%%%%%%%%%%%%%%%%%%%%%%%%%%%%%%%%%%%%%%%%%%%%%%%%%%%%%%%%%%%%%%%%%%%%

\section{Multivariable adjunctions}\label{adj}

In this section we define multivariable adjunctions.  The basic idea is that for an ``$n$-variable adjunction'' we have $n+1$ categories $A_0, \cdots, A_n$ and $n+1$ multifunctors, each of which has one of the $A\opb_i$ as its target, and the product of the other $n$ categories as its source.  These multifunctors can all be restricted to functors with a single category as their source, by fixing an object in each of the other categories.  For every pair $i \neq j$ there is a pair of contravariant functors obtained in this way involving $A_i$ and $A_j$.  These should be in a specified adjunction; moreover, of course, all these adjunctions should be coherent in an appropriate way.

We first give the definition of this structure, and then immediately prove  Theorem~\ref{twotwo} giving a more ``economical'' characterisation, in which \emph{a priori} we specify only one multifunctor, and a family of 1-variable adjoints for it.  Using standard results about parametrised representability, these 1-variable adjoints then extend uniquely to $n$-variable multifunctors with the required structure.  It is the characterisation in Theorem~\ref{twotwo} that we will use in the rest of the work.

% Note that the left/right terminology these adjunctions is potentially somewhat more confusing than in 1-variable adjunctions; here again it will turn out to be clearest to think in terms of Theorem~\ref{twotwo} rather than the definition itself.  
% 

% 
% \subsection{Mutual left adjunctions}
% 
% We first need to fix some notation for ``mutual left adjoints''.
% 
% \begin{definition}
% Consider functors
% 
% \[\begin{array}{rclcrcl}
%  C & \tmap{F} & D\opb & \mbox{ so } & C{\opb} & \tmap{F{\opb}} & D \\
%  D & \tmap{G} & C{\opb} & \mbox{ so } & D{\opb} & \tmap{G{\opb}} & C 
% \end{array}\]
% %
% A \demph{mutual left adjunction} of $F$ and $G$ is an adjunction $F\opb \ladj G$ or equivalently $G\opb \ladj F$. That is, isomorphisms
% \[\begin{array}{lccc}
% &   D(F\opb X, Y) & \cong & C\opb(X, GY) \\
% \mbox{i.e.} & D(FX, Y) & \cong & C(GY, X) \\
% \mbox{or} & D\opb(Y, FX) & \cong & C(G\opb Y, X) 
%   \end{array}\]
% natural in $X$ and $Y$.
% 
% Similarly a \demph{mutual right adjunction} is given by 
% \[D(X, FY) \cong C(X, GY)\]
% natural in $X$ and $Y$.
% 
% \end{definition}
% 
% Note that it is more usual to write contravariant functors as $C\opb \tra D$ but for $n$-variable adjunctions the above convention is more convenient.
% 
% **I originally drew some naturality squares here as I was afraid of being confused.**

\subsection{Definition of multivariable adjunctions}

\begin{mydefinition}
 Let $n \in \N$.  An \demph{$n$-variable (mutual) left adjunction} is given by the following data and axioms.

\begin{itemize}
 \item Categories $A_0, \ldots, A_{n}$.
\item Functors 
\[\begin{array}{lcl}
   A_1 \times A_2 \times \cdots \times A_{n-1} \times A_n & \tmap{F_0} & A_{0}\opb \\[3pt]
   A_2 \times A_3 \times \cdots \times A_n \times A_{0} & \tmap{F_1} & A_{1}\opb \\[3pt]
%    A_3 \times A_4 \times \cdots \times A_{0} \times A_1 & \tmap{F_2} & A_{2}\opb \\[3pt]
%    A_4 \times \cdots \times A_{0} \times A_1 \times A_2 & \tmap{F_3} & A_{3}\opb \\
 & \vdots & \\[3pt]
   A_{i+1} \times  \cdots \times A_{i-1} & \tmap{F_i} & A_{i}\opb \\
 & \vdots & \\[3pt]
   A_{0}  \times \cdots \times A_{n-1} & \tmap{F_{n}} & A_{n}\opb. 
  \end{array}\]

Here the subscripts are all to be taken mod $n+1$. Where possible, we will adopt the convention that the subscript on a multifunctor matches the subscript of its target category.

\item  For all $0 \leq i \leq n$, and for all $a_{i+1} \in A_{i+1}, \ldots, a_{i-2} \in A_{i-2}$ a mutual left adjunction between

\[
\psset{unit=0.1cm,labelsep=2pt,nodesep=3pt}
\pspicture(40,12)
\rput[r](0,9){\rnode{a1}{$A_i$}}  % named node, with something placed there
\rput[l](40,9){\rnode{a2}{$A_{i-1}\opb$}}  % named node, with something placed there
\rput[r](0,1){\rnode{b1}{$A_{i-1}$}}  % named node, with something placed there
\rput[l](40,1){\rnode{b2}{$A_i\opb$}}  % named node, with something placed there

\ncline{->}{a1}{a2} \naput{{\scriptsize $F_{i-1}(\uscore, a_{i+1}, a_{i+2}, \ldots, a_{i-2})$}}
\ncline{->}{b1}{b2} \naput{{\scriptsize $F_{i}(a_{i+1}, a_{i+2}, \ldots, a_{i-2}, \uscore)$}}

\endpspicture
\]
thus isomorphisms
\[A_{i-1}\big(F_{i-1}(a_i, \ldots, a_{i-2}), a_{i-1}\big) \cong A_i\big(F_{i}(a_{i+1}, \ldots, a_{i-1}), a_i\big)
\]
natural in $a_{i-1}$ and $a_{i}$.  If we use the shorthand $\hat{a}_i$ for the sequence $a_{i+1}, \ldots, a_{i-1}$, this isomorphism takes the appealing form
\[A_{i-1}\big(F_{i-1}(\hat{a}_{i-1}), a_{i-1}\big) \iso A_i\big(F_i(\hat{a}_i), a_i\big).\]

\end{itemize}

The following axioms must be satisfied:

\begin{itemize}

\item the above isomorphisms must additionally be natural in all variables, and 

\item the ``cycle'' of isomorphisms commutes:
\[
\psset{unit=0.1cm,labelsep=1pt,nodesep=3pt}
\pspicture(0,-5)(60,35)
\rput(0,15){\rnode{an}{$A_n\big(F_n(\hat{a}_n), a_n\big)$}}  
\rput(10,30){\rnode{a0}{$A_0\big(F_0(\hat{a}_0), a_0\big)$}}  
\rput(50,30){\rnode{a1}{$A_1\big(F_1(\hat{a}_1), a_1\big)$}}  
\rput(60,15){\rnode{a2}{$A_2\big(F_2(\hat{a}_2), a_2\big)$}} 
\rput(50,0){\rnode{a3}{$A_3\big(F_3(\hat{a}_3), a_3\big)$}}  

\rput(10,0){\rnode{a4}{\rotatebox{125}{$\ldots$}}}  
\rput(45,-4){\rotatebox{52}{$\ldots$}}

\ncline{->}{an}{a0} \naput{{\scriptsize $\sim$}}
\ncline{->}{a0}{a1} \naput{{\scriptsize $\sim$}}
\ncline{->}{a1}{a2} \naput{{\scriptsize $\sim$}}
\ncline{->}{a2}{a3} \naput{{\scriptsize $\sim$}}

\ncline{->}{a4}{an} \naput{{\scriptsize $\sim$}}
% \ncline{->}{a}{a} \naput{{\scriptsize $\sim$}}

\endpspicture
\]

\end{itemize}

We say that the functor $F_0$ is \demph{equipped with $n$-variable left adjoints} $F_1, \cdots, F_n$. This terminology makes more sense in the light of the following theorem.
\end{mydefinition}

\begin{theorem}\label{theoremA}\label{twotwo}
 The following description precisely corresponds to an $n$-variable left adjunction.

\begin{itemize}
 \item categories $A_0, \ldots, A_{n}$
\item a functor  $A_1 \times \cdots \times A_n  \tmap{F_0}  A_{0}\opb$

\item for all $0,j, k$ distinct, and for all $a_j \in A_j$, a mutual left adjoint for the functor
\[F_0(a_1, \ldots, a_{k-1}, \uscore, a_{k+1}, \ldots, a_n) \:  A_k \tra A_{0}\opb\]

\end{itemize}

\end{theorem}

\begin{myremark}
Note we say that $F$ is equipped with $n$-variable \demph{left} adjoints if each of its 1-variable restrictions has a \demph{left} adjoint. $F$ is equipped with $n$-variable \demph{right} adjoints if each of its 1-variable restrictions has a \demph{right} adjoint.
\end{myremark}

To prove this we use the following result of Mac Lane \cite[IV.7, Theorem 3]{mac1}.

\begin{theorem}\label{theoremB}
 Given categories $A,B,C$, a functor $F \:  A \times B \tra C\opb$, and for all $b \in B$ a mutual left adjoint $G(b,\uscore): C \tra A\opb$ for the functor
\[F(\uscore, b)  \:  A \tra C\opb\]
i.e. isomorphisms
\begin{equation}\label{one}C(F(a,b), c) \iso A(G(b,c), a)\end{equation}
natural in $a$ and $c$, there is a unique way to extend the functors
\[G(b,\uscore) \:  C \tra A\opb\]
to a single functor
\[G \:  B \times C \tra A\opb\]
such that the isomorphism (\ref{one}) is also natural in $b$.
\end{theorem}

This is a standard result about parametrised representability; we give a 2-categorical expression of Mac Lane's proof, as this will be useful later.

\begin{prf}
 We write $1 \tmap{b} B$ for the functor picking out the object $b \in B$.  The hypothesis of the theorem then says that for each such $b$ we have a right adjoint for the composite
\[
% a1 a2 a3
\psset{unit=0.1cm,labelsep=2pt,nodesep=2pt}
\pspicture(30,5)
\rput(-5,0){\rnode{a1}{$A\opb$}}  % named node, with something placed there
\rput(15,0){\rnode{a2}{$A\opb \times B\opb$}}  % named node, with something placed there
\rput(35,0){\rnode{a3}{$C$}}  % named node, with something placed there
\ncline{->}{a1}{a2} \naput{{\scriptsize $1 \times b\opb$}}
\ncline{->}{a2}{a3} \naput{{\scriptsize $F\opb$}}
\endpspicture
\]
which we call
\[
% b1  b3
\psset{unit=0.1cm,labelsep=2pt,nodesep=2pt}
\pspicture(30,5)

\rput(0,0){\rnode{b1}{$C$}}  % named node, with something placed there
\rput(30,00){\rnode{b3}{$A\opb$}}  % named node, with something placed there

\ncline{->}{b1}{b3} \naput{{\scriptsize $G(b,\uscore)$}}
\endpspicture
\]
with unit and counit

% a1 a2
%    b2
%    c2 c3

\[
\psset{unit=0.1cm,labelsep=2pt,nodesep=2pt}
\pspicture(0,-5)(60,40)
\rput(0,30){\rnode{a1}{$C$}}  % named node, with something placed there
\rput(30,30){\rnode{a2}{$A\opb$}}  % named node, with something placed there
\rput(30,15){\rnode{b2}{$A\opb \times B\opb$}}  % named node, with something placed there
\rput(30,0){\rnode{c2}{$C$}}  % named node, with something placed there
\rput(60,0){\rnode{c3}{$A\opb.$}}  % named node, with something placed there

\ncline{->}{a1}{a2} \naput{{\scriptsize $G(b,\uscore)$}}
\ncline{->}{a1}{c2} \nbput{{\scriptsize $1$}}
\ncline[nodesepA=3pt]{->}{a2}{b2} \nbput{{\scriptsize $1 \times b\opb$}}
\ncline{->}{a2}{c3} \naput{{\scriptsize $1$}}
\ncline[nodesepA=3pt]{->}{b2}{c2} \naput{{\scriptsize $F\opb$}}
\ncline{->}{c2}{c3} \nbput{{\scriptsize $G(b,\uscore)$}}

\rput(19,24){\rnode{d1}{$$}}
\rput(14,19){\rnode{d2}{$$}}
\rput(44,9){\rnode{e1}{$$}}
\rput(39,4){\rnode{e2}{$$}}

\psset{labelsep=1.5pt}

\ncline[doubleline=true, labelsep=1pt]
{->}{d1}{d2} \nbput{{\scriptsize $\varepsilon_b$}}
\ncline[doubleline=true]
{->}{e1}{e2} \naput{{\scriptsize $\eta_b$}}

\endpspicture
\]
Now, extending the individual functors
\[G(b, \uscore) \: C \tra A\opb\]
to a functor
\[G \: B \times C \tra A\opb\]
%  
% \[\begin{array}{ccccc}
% G(b, \uscore) &:& C & \tra & A\opb \\
% G &:& B \times C & \tra & A\opb
% \end{array}\]
consists of giving, for each morphism $b_1 \tmap{f} b_2$ in $B$, a natural transformation
\[
\psset{unit=0.1cm,labelsep=2pt,nodesep=3pt}
\pspicture(20,20)

\rput(0,10){\rnode{a1}{$C$}}  % named node, with something placed there
\rput(20,10){\rnode{a2}{$A\opb$}}  % named node, with something placed there

\pcline[doubleline=true]{->}(10,13)(10,7)

\ncarc[arcangle=45]{->}{a1}{a2}\naput{{\scriptsize $G(b_1, \uscore)$}}
\ncarc[arcangle=-45]{->}{a1}{a2}\nbput{{\scriptsize $G(b_2, \uscore)$}}

\endpspicture
\]
and checking functoriality.  The \nt\ is given as the mate of 
\[
\psset{unit=0.1cm,labelsep=2pt,nodesep=2pt}
\pspicture(50,20)

\rput(5,10){\rnode{a1}{$A\opb$}}  % named node, with something placed there
\rput(30,10){\rnode{a2}{$A\opb \times B\opb$}}  % named node, with something placed there
\rput(50,10){\rnode{a3}{$C$}}

\pcline[doubleline=true]{->}(13,13)(13,7) \naput{{\scriptsize $1 \ttimes f\opb$}}

\ncarc[arcangle=45]{->}{a1}{a2}\naput{{\scriptsize $1 \ttimes b_2\opb$}}
\ncarc[arcangle=-45]{->}{a1}{a2}\nbput{{\scriptsize $1\ttimes b_1\opb$}}
\ncline{->}{a2}{a3} \naput{{\scriptsize $F$}}

\endpspicture
\]
that is

% a1 a2
%    b2
%    c2 c3

\[
\psset{unit=0.1cm,labelsep=2pt,nodesep=2pt}
\pspicture(0,-5)(60,35)
\rput(0,30){\rnode{a1}{$C$}}  % named node, with something placed there
\rput(30,30){\rnode{a2}{$A\opb$}}  % named node, with something placed there
\rput(30,13){\rnode{b2}{$A\opb \times B\opb$}}  % named node, with something placed there
\rput(30,0){\rnode{c2}{$C$}}  % named node, with something placed there
\rput(60,0){\rnode{c3}{$A\opb.$}}  % named node, with something placed there

\ncline{->}{a1}{a2} \naput{{\scriptsize $G(b_1,\uscore)$}}
\ncline{->}{a1}{c2} \nbput{{\scriptsize $1$}}
%\ncline[nodesepA=3pt]{->}{a2}{b2} \nbput{{\scriptsize $1 \times b\opb$}}
\ncline{->}{a2}{c3} \naput{{\scriptsize $1$}}
\ncline[nodesepA=3pt]{->}{b2}{c2} \naput[labelsep=1pt]{{\scriptsize $F\opb$}}
\ncline{->}{c2}{c3} \nbput{{\scriptsize $G(b_2,\uscore)$}}

\rput(19,24){\rnode{d1}{$$}}
\rput(14,19){\rnode{d2}{$$}}
\rput(44,9){\rnode{e1}{$$}}
\rput(39,4){\rnode{e2}{$$}}

\ncline[doubleline=true, labelsep=1pt]
{->}{d1}{d2} \nbput{{\scriptsize $\varepsilon_{b_1}$}}
\ncline[doubleline=true]
{->}{e1}{e2} \naput{{\scriptsize $\eta_{b_2}$}}

\ncarc[arcangle=38]{->}{a2}{b2}
\ncarc[arcangle=-38]{->}{a2}{b2}
\pcline[doubleline=true]
{->}(32,23)(27,23) \naput[npos=0.3]{{\scriptsize $1 \ttimes f\opb$}}

\endpspicture
\]

Functoriality then follows from the functoriality of the mates correspondence.  Now we further need that the isomorphism
\[C(F(a,b),c) \tmap{\sim} A(G(b,c), a)\]
is natural in $b$.  A diagram chase shows that this is equivalent to the following diagram commuting for all $f \:  b_1 \tra b_2$ in B:
% a1 a2
% a3 a4

\[
\psset{unit=0.1cm,labelsep=3pt,nodesep=3pt}
\pspicture(0,-5)(40,25)

%%%%%%%%%% top

\rput(0,21){\rnode{a1}{$G(b_2, F(a,b_1))$}}  % top left
\rput(35,21){\rnode{a2}{$G(b_1, F(a,b_1))$}}  % top right
\rput(0,1){\rnode{a3}{$G(b_2, F(a,b_2))$}}  % bottom left
\rput(35,1){\rnode{a4}{$a$}}  % bottom right

\ncline{->}{a1}{a2} \naput{{\scriptsize $G(f,1)$}} % top
\ncline{->}{a3}{a4} \nbput{{ $\eta_{b_2, a}$}} % bottom
\ncline{->}{a1}{a3} \nbput{{\scriptsize $G(1, F(1,f))$}} % left
\ncline{->}{a2}{a4} \naput{{ $\eta_{b_1, a}$}} % right

\endpspicture
\]
or dually an analogous diagram involving $\varepsilon$'s:
% a1 a2
% a3 a4

\[
\psset{unit=0.1cm,labelsep=3pt,nodesep=3pt}
\pspicture(0,-5)(40,25)

%%%%%%%%%% top

\rput(0,21){\rnode{a1}{$F(b_2, G(a,b_1))$}}  % top left
\rput(35,21){\rnode{a2}{$F(b_1, G(a,b_1)$}}  % top right
\rput(0,2){\rnode{a3}{$F(b_2, G(a,b_2))$}}  % bottom left
\rput(35,2){\rnode{a4}{$a$.}}  % bottom right

\ncline{->}{a1}{a2} \naput{{\scriptsize $F(f,1)$}} % top
\ncline{->}{a3}{a4} \nbput{{ $\varepsilon_{b_2, a}$}} % bottom
\ncline{->}{a1}{a3} \nbput{{\scriptsize $F(1, G(1,f))$}} % left
\ncline{->}{a2}{a4} \naput{{ $\varepsilon_{b_1, a}$}} % right

\endpspicture
\]
2-categorically this is 

% a1 a2 a3
%       b3
%       c3

\begin{equation}\label{rhino}
\psset{unit=0.1cm,labelsep=2pt,nodesep=2pt}
\pspicture(0,-10)(80,40)

\rput[b](6,0){
\pspicture(40,35)
\rput(0,33){\rnode{a1}{$C$}}
\rput(20,33){\rnode{a2}{$B \times C$}}
\rput(40,33){\rnode{a3}{$A\opb$}}
\rput(40,15){\rnode{b3}{$A\opb \times B \opb$}}
\rput(40,0){\rnode{c3}{$C$}}

\ncline{->}{a2}{a3} \naput{{\scriptsize $G$}}
\ncline{->}{a3}{b3} \naput{{\scriptsize $1\ttimes b_2\opb$}}
\ncline{->}{b3}{c3} \naput{{\scriptsize $F\opb$}}

%\ncarc[arcangleA=290,arcangleB=120]{->}{a1}{c3} \nbput{{\scriptsize $1$}}

\nccurve[angleA=300,angleB=160]{->}{a1}{c3} \nbput{{\scriptsize $1$}}

\psset{labelsep=1pt}
\ncarc[arcangle=35,nodesepA=3pt]{->}{a1}{a2} \naput{{\scriptsize$b_1 \hh{-2pt} \times \hh{-2pt} 1$}}
\ncarc[arcangle=-35,nodesepA=3pt]{->}{a1}{a2} \nbput{{\scriptsize$b_2 \hh{-2pt} \times \hh{-2pt}  1$}}

\psset{doubleline=true}
\rput{0}(21,17){\pcline{->}(4,4)(0,0) \naput[labelsep=1pt]{{\scriptsize $\varepsilon_{b_2}$}}}
\rput{0}(7,31){\pcline{->}(0,4)(0,0) \naput{{\scriptsize $f \hh{-2pt} \times \hh{-2pt}  1$}}}

\endpspicture}

\rput(41,18){$=$}

\rput[b](65,0){
\pspicture(40,35)
\rput(0,33){\rnode{a1}{$C$}}
\rput(20,33){\rnode{a2}{$B \times C$}}
\rput(40,33){\rnode{a3}{$A\opb$}}
\rput(40,15){\rnode{b3}{$A\opb \times B \opb$}}
\rput(40,0){\rnode{c3}{$C$.}}

\ncline{->}{a1}{a2} \naput{{\scriptsize $b_1 \ttimes 1$}}
\ncline{->}{a2}{a3} \naput{{\scriptsize $G$}}
%\ncline{->}{a3}{b3} \naput{{\scriptsize $1\ttimes b_2\opb$}}
\ncline{->}{b3}{c3} \naput{{\scriptsize $F\opb$}}

%\ncarc[arcangleA=290,arcangleB=120]{->}{a1}{c3} \nbput{{\scriptsize $1$}}

\nccurve[angleA=300,angleB=160]{->}{a1}{c3} \nbput{{\scriptsize $1$}}

\psset{labelsep=1pt}
\ncarc[arcangle=44,nodesepA=3pt]{->}{a3}{b3} \naput{{\scriptsize$1 \hh{-2pt} \times \hh{-2pt} b_2\opb$}}
\ncarc[arcangle=-44,nodesepA=3pt]{->}{a3}{b3} \nbput{{\scriptsize$1 \hh{-2pt} \times \hh{-2pt}  b_1\opb$}}

\psset{doubleline=true}
\rput{0}(21,17){\pcline{->}(4,4)(0,0) \naput[labelsep=1pt]{{\scriptsize $\varepsilon_{b_1}$}}}
\rput[c](37,25){\pcline{->}(5,0)(0,0) \naput[npos=0.3]{{\scriptsize $1 \ttimes f\opb$}}}

\endpspicture
}

\endpspicture
\end{equation}

Now by our definition we have 

\begin{equation}\label{rhino2}
\psset{unit=0.1cm,labelsep=2pt,nodesep=2pt}
\pspicture(70,40)

%%% right hand side

\rput[b](71,12){
\pspicture(50,20)

\rput(5,10){\rnode{a1}{$C$}}  % named node, with something placed there
\rput(30,10){\rnode{a2}{$B \times C$}}  % named node, with something placed there
\rput(50,10){\rnode{a3}{$A\opb$}}

\pcline[doubleline=true]{->}(13,13)(13,7) \naput{{\scriptsize $f \ttimes 1$}}

\ncarc[arcangle=45]{->}{a1}{a2}\naput{{\scriptsize $b_1 \ttimes 1$}}
\ncarc[arcangle=-45]{->}{a1}{a2}\nbput{{\scriptsize $b_2 \ttimes 1$}}
\ncline{->}{a2}{a3} \naput{{\scriptsize $G$}}

\endpspicture}

\rput(38,22){$=$}

%%%%%%%%%%%%%%%%%%%%%%%%%%%%%%%%%%%%%%%%%%%%%%%%%%%%%%%%%%%%%%%%%%%
% left hand side 

% a1 a2
%    b2
%    c2 c3

\rput[b](0,0){
\psset{unit=0.1cm,labelsep=2pt,nodesep=2pt}
\pspicture(0,-5)(60,35)
\rput(0,30){\rnode{a1}{$C$}}  % named node, with something placed there
\rput(30,30){\rnode{a2}{$A\opb$}}  % named node, with something placed there
\rput(30,13){\rnode{b2}{$A\opb \times B\opb$}}  % named node, with something placed there
\rput(30,0){\rnode{c2}{$C$}}  % named node, with something placed there
\rput(60,0){\rnode{c3}{$A\opb$}}  % named node, with something placed there

\ncline{->}{a1}{a2} \naput{{\scriptsize $G(b_1,\uscore)$}}
\ncline{->}{a1}{c2} \nbput{{\scriptsize $1$}}
%\ncline[nodesepA=3pt]{->}{a2}{b2} \nbput{{\scriptsize $1 \times b\opb$}}
\ncline{->}{a2}{c3} \naput{{\scriptsize $1$}}
\ncline[nodesepA=3pt]{->}{b2}{c2} \naput[labelsep=1pt]{{\scriptsize $F\opb$}}
\ncline{->}{c2}{c3} \nbput{{\scriptsize $G(b_2,\uscore)$}}

\rput(19,24){\rnode{d1}{$$}}
\rput(14,19){\rnode{d2}{$$}}
\rput(44,9){\rnode{e1}{$$}}
\rput(39,4){\rnode{e2}{$$}}

\ncline[doubleline=true, labelsep=1pt]
{->}{d1}{d2} \nbput{{\scriptsize $\varepsilon_{b_1}$}}
\ncline[doubleline=true]
{->}{e1}{e2} \naput{{\scriptsize $\eta_{b_2}$}}

\ncarc[arcangle=38]{->}{a2}{b2}
\ncarc[arcangle=-38]{->}{a2}{b2}
\pcline[doubleline=true]
{->}(32,23)(27,23) \naput[npos=0.3]{{\scriptsize $1 \ttimes f\opb$}}

\endpspicture}

\endpspicture
\end{equation}
\noi since the right-hand side is the definition of $G$ on morphisms of $B$.  Then equation (\ref{rhino}) follows from a triangle identity for $\eta_{b_2}$ and $\varepsilon_{b_2}$; dually the equation for $\eta$ holds by a triangle identity for $\eta_{b_1}$ and $\varepsilon_{b_1}$. 

For uniqueness we suppose we have a functor $G$ satisfying the naturality condition as shown in diagram~(\ref{rhino}) above.  Then as above, equation (\ref{rhino2}) must hold, showing that our construction of $G$ is unique. 
\end{prf}

\begin{prfof}{Theorem~\ref{theoremA}}
First we show that the structure in the theorem gives rise to an $n$-variable left adjunction.  First we need to define for all $i \neq 0$ a functor
\[F_i \: A_{i+1} \times \cdots \times A_{i-1} \tra A_{i}\opb.\]
Now, we have for for all $a_1, \ldots, a_{i-1}, a_{i+1}, \ldots, a_n$ a left adjoint for the functor
\[F_0(a_1, \ldots, a_{i-1}, \uscore, a_{i+1}, \ldots, a_n) \: A_{i} \tra A_0\opb,\]
equivalently a right adjoint for its opposite
\[F_0\opb(a_1, \ldots, a_{i-1}, \uscore, a_{i+1}, \ldots, a_n) \: A_{i}\opb \tra A_0\]
called, say
\[F_i(a_{i+1}, \ldots, a_n, \uscore, a_1, \ldots, a_{i-1}) \: A_{0} \tra A_{i}\opb.\]
By Theorem~\ref{theoremB} it extends uniquely to a functor
\[F_i \: A_{i+1} \times \cdots \times A_{i-1} \tra A_{i}\opb\]
making the isomorphism
\[A_0\big(F_0(\ha_0), a_0 \big) \iso A_i \big(F_i(\ha_i) , a_i \big)\]
natural in every variable (where \emph{a priori} it was only natural in $a_0$ and $a_i$).  This is by putting
\[\begin{array}{rcl}
   A &=& A_{i} \\
   B &=& A_{i+1} \times \cdots \times A_{i-1} \\
   C &=& A_0
  \end{array}\]
in the theorem.  It remains to show that we have the correct adjunctions.  Now by the above hom-set isomorphism we construct the composite isomorphism
\[A_{i-1}\big(F_{i-1}(\ha_{i-1}), a_{i-1}\big) \tmap{\sim} A_{0}\big(F_{0}(\ha_{0}), a_{0}\big) \tmap{\sim} A_{i}\big(F_{i}(\ha_{i}), a_{i}\big)\]
which we already know to be natural in every variable, and by construction the cycle of isomorphisms commutes as required. 

Conversely given an $n$-variable adjunction we use the cycle of isomorphisms to specify an isomorphism
\[A_{0}\big(F_{0}(\ha_{0}), a_{0}\big) \tmap{\sim} A_{1}\big(F_{1}(\ha_{1}), a_{1}\big) \tmap{\sim} \cdots \tmap{\sim} A_{i}\big(F_{i}(\ha_{i}), a_{i}\big).\]
Then, fixing all variables except $a_{i}$ and $a_0$ we get the required adjunction. \end{prfof}

It is instructive to work through this definition for some small values of $n$.

\begin{myexample} $n=1$

A 1-variable adjunction is just an ordinary adjunction, but in the notation of the definition it is given by
\begin{itemize}
 \item categories $A_0, A_1$, 
\item functors $\begin{array}[t]{cc}
		  A_1 \tmap{F_0} A_0\opb \\
		  A_0 \tmap{F_1} A_1\opb
                \end{array}$, and
\item an adjunction $F_0\opb \ladj F_1$.

\end{itemize}
\end{myexample}

\begin{myexample}\label{cycliso} $n=2$
 
A 2-variable adjunction is given by categories, functors and adjunctions as follows:
\[\begin{array}{cc}
\hh{3em}   A \times B \map{F} C\opb \hh{3em} & F(\uscore, b)\opb \ladj G(b, \uscore) \\
   B \times C \map{G} A\opb & G(\uscore, c)\opb \ladj H(c, \uscore) \\
   C \times A \map{H} B\opb & H(\uscore, a)\opb \ladj F(a, \uscore)
  \end{array}\]
given by a ``cycle of isomorphisms''
\[
\psset{unit=0.1cm,labelsep=0pt,nodesep=3pt}
\pspicture(30,17)
\rput(0,15){$C(F(a,b), c)$}  
\rput(30,15){$A(G(b,c), a)$}
\rput(15,2){$B(H(c,a),b)$}

\rput(15,15){$\iso$}
\rput(6,8.7){\rotatebox{320}{$\iso$}}
\rput(24,9){\rotatebox{40}{$\iso$}}

\endpspicture
\]
natural in $a$, $b$ and $c$.

Theorem~\ref{theoremA} says that to specify this it is equivalent to specify the functor $F$ along with, for each $a \in A$  and $b \in B$ left adjoints for the functors $F(\uscore, b)$ and $F(a, \uscore)$, that is functors
\[G(b,\uscore) \: C \tra A\opb\]
\[H(\uscore, a) \: C \tra B\opb\]
and isomorphisms
\[\begin{array}{ccl}
   A(G(b,c),a) \iso C(F(a,b), c) & \hh{3em} & \mbox{natural in $a$ and $c$}\\
A(G(b,c), a) \iso B(H(c,a), b) && \mbox{natural in $b$ and $c$.}
  \end{array}\]
\end{myexample}

Note that the original definition has $n+1$ adjunctions specified cyclically, each involving a pair of ``numerically adjacent'' categories and naturality in all $n+1$ variables; Theorem~\ref{theoremA} specifies $n$ adjunctions, each involving $A_0$ and one other category, and natural only in 2 variables.

\begin{remark}
For $n=0$ it is useful to say that a ``0-variable adjunction'' is a functor $1 \tra A$ as these will be the 0-ary maps in our eventual multicategory structure. The fact that these compose is the following lemma.
\end{remark}

\begin{lemma}
Consider an $n$-variable adjunction as above.  Fix $0\leq k \leq n$ and $a_k \in A_k$.  Then fixing $a_k$ in each functor $F_i, i \neq k$ yields an $(n-1)$-variable adjunction in the evident way. 
\end{lemma}

Obviously we can repeat this process to fix any number of variables to restrict an multivariable adjunction to one in a smaller number of variables.  Note that apart from being a crucial component of the eventual multicategory structure, this fact is also used in the proof of the multivariable mates correspondence (Theorem~\ref{nvarmates}).

% **comment about left and right?  Actually what we have is: for each variable in turn, a left adjoint.  Alternatively we could have: for each variable in turn, a right adjoint.  I think that's the clearest way of putting it.**

\begin{proposition}\label{leftright}
An $n$-variable left adjunction of functors $F_0, \ldots, F_n$ is equivalently an $n$-variable right adjunction of $F_0\opb, \ldots, F_n\opb$. 
\end{proposition}

%%%%%%%%%%%%%%%%%%%%%%%%%%%%%%%%%%%%%%%%%%%%%%%%%%%%%%%%%%%%%%%%%%%%%
%%%%%%%%%%%%%%%%%%%%%%%%%%%%%%%%%%%%%%%%%%%%%%%%%%%%%%%%%%%%%%%%%%%%%
%%%%%%%%%%%%%%%%%%%%%%%%%%%%%%%%%%%%%%%%%%%%%%%%%%%%%%%%%%%%%%%%%%%%%
%%%%%%%%%%%%%%%%%%%%%%%%%%%%%%%%%%%%%%%%%%%%%%%%%%%%%%%%%%%%%%%%%%%%%
%%%%%%%%%%%%%%%%%%%%%%%%%%%%%%%%%%%%%%%%%%%%%%%%%%%%%%%%%%%%%%%%%%%%%
%%%%%%%%%%%%%%%%%%%%%%%%%%%%%%%%%%%%%%%%%%%%%%%%%%%%%%%%%%%%%%%%%%%%%

\subsection{A motivating example}\label{homtensor}

We begin by presenting the standard example of a 2-variable adjunction that we have generalised, the ``tensor/hom/cotensor'' adjunction.  The only slightly tricky thing is taking care of the dualities.

Let $\cV$ be a monoidal category, so we have a functor
\[\uscore \otimes \uscore \: \cV \times \cV \tra \cV.\]
Then \cV\ is biclosed if 
\begin{itemize}
 \item  $\forall b \in \cV$ the functor $\uscore \otimes b$ has a right adjoint $[\hh{1.5pt}b\hh{1pt}, \uscore\hh{2pt}]$ (``hom''), and
\item  $\forall a \in \cV$ the functor $a \otimes \uscore$ has a right adjoint $a \pitchfork \uscore$ (``cotensor'').
\end{itemize}

The first adjunction gives us isomorphisms
\[\cV(a\otimes b, c) \cong \cV(a, [b,c])\]
natural in $a$ and $c$; by parametrised representability the functor
\[ [b, \uscore] \: \cV \lra \cV\]
extends to a functor
\[ [\uscore, \uscore] \: \cV\opb \times \cV \tra \cV\]
uniquely making the isomorphisms natural in $b$ as well.

Similarly for the second adjunction we get a functor
\[ \uscore \pitchfork \uscore \: \cV\opb \times \cV \tra \cV\]
making the isomorphism                                                      
\[\cV(a \otimes b, c ) \cong \cV(b, a \pitchfork c)\]
natural in all three variables.

\vv{6pt}

More generally for categories $A,B,C$ a tensor-hom-cotensor adjunction consists of functors and adjunctions

\[\psset{labelsep=1pt}\begin{array}{cccclccc}
 \hh{1pt} A \hh{1pt} \times \hh{1pt} B & \hh{-10pt} \ltmap{\suscore \hh{1pt} \otimes \hh{1pt} \suscore\hh{3pt}} \hh{-10pt} & C &\hh{2.5em}& \forall a \in A \ & a \otimes \uscore & \ladj & a \pitchfork \uscore \\[2pt]
 B\opb \times C & \hh{-10pt} \ltmap{\shomblank\hh{3pt}} \hh{-10pt} & A &\hh{3em}& \forall b \in B & \uscore \otimes b & \ladj & [\hh{1.5pt}b\hh{1pt}, \uscore\hh{2pt}] \\[2pt]
 A\opb \times C & \hh{-10pt} \ltmap{\suscore \hh{1pt} \pitchfork \hh{1pt} \suscore\hh{3pt}} \hh{-10pt} & B &\hh{3em}& \forall c \in C & [ \hh{2pt} \uscore \hh{1pt},\hh{1pt} c \hh{2pt}]\opb & \ladj & \uscore \pitchfork c 
  \end{array}\]
and by parametrised representability it follows that the following isomorphisms are natural in all three variables:
\[ A(a,[b,c]) \cong B(b, a \pitchfork c) \cong C(a \otimes b, c).\]

For our standard framework with functors
\[\begin{array}{ccc}
   A_1 \times A_2 \tmap{F_0} A_0\opb \\
   A_2 \times A_0 \tmap{F_1} A_1\opb \\
   A_0 \times A_1 \tmap{F_2} A_2\opb 
  \end{array}\]
we can put
\[\setlength{\arraycolsep}{0.2em}
\begin{array}{rclcrcl}
   A_1 &=& A\opb & \hh{2em} & F_0 &=& (\uscore \otimes \uscore)\opb \\
   A_2 &=& B\opb     && F_1 &=& [\uscore, \uscore] \\
   A_0 &=& C && F_2 &=& \uscore \pitchfork \uscore
  \end{array}\]
and $F_0, F_1$, $F_2$ then form a 2-variable left adjunction (although $\uscore \pitchfork \uscore$ now has domain $C \times A\opb$ instead of  $A\opb \times 	C$).

This 2-variable adjunction is  is the starting point of the discussion in Section~\ref{modelcat}. 

\subsection{Composition}

Just as ordinary adjunctions can be composed (with care over directions) so can $n$-variable adjunctions, with care over directions, dualities and arities.  The only difficulty in the following theorem is the notation.  The idea is to compose $n$-variable adjunctions in the manner of multimaps in a multicategory; indeed this is what they will be in Section~\ref{multicat}.  In this section all multivariable adjunctions are left adjunctions; of course the right adjunctions follow dually.

\begin{theorem}\label{twopointseven}
Suppose we have the following multivariable left adjunctions.
\[\begin{array}{lll}
   A_{11} \times \cdots \times A_{1n_1} \tmap{ F_{10}} A_{10}\opb = B_1 \hh{1em}& \mbox{with $n_1$-variable adjoints} & F_{11}, \ldots, F_{1n_1} \\
    A_{21} \times \cdots \times A_{2n_2} \tmap{F_{20}} A_{20}\opb = B_2 & \mbox{with $n_2$-variable adjoints} & F_{21}, \ldots, F_{2n_2} \\
 \hh{5.5em} \vdots & \\
   A_{k1} \times \cdots \times A_{kn_k} \tmap{F_{k0}} A_{k0}\opb = B_k & \mbox{with $n_k$-variable adjoints} & F_{k1}, \ldots, F_{kn_k} \\[1em]
\mbox{and} \\[1em]
B_1 \times \cdots \times B_k \tmap{G_0} B_0\opb  = A_{00}\opb & \mbox{with $k$-variable adjoints} & G_1, \ldots, G_k 
  \end{array}\]
Then the composite functor 
\[G_0(F_{10}, F_{20}, \ldots, F_{k0}) \: A_{11} \times \cdots \times A_{1n_1} \times \cdots \times A_{k1} \times \cdots \times A_{kn_k} \tra B_{0}\opb\]
is canonically equipped with $(n_1 + \cdots + n_k)$-variable adjoints.  This composition makes categories and multivarible adjunctions into a multicategory.  
\end{theorem}

\begin{prf}
We write $n_1 + \cdots + n_k=m$ and call the above composite $H_{00}$.  We must construct $m$-variable adjoints for $H_{00}$, so first we need $m$ functors which we call
\[\begin{array}{l}
   H_{11}, \ldots, H_{1n_1} \\
   H_{21}, \ldots, H_{2n_2} \\
\hh{2em} \vdots \\
H_{k1}, \ldots, H_{kn_k}
  \end{array}\]
where $H_{ij}$ has target category $A_{ij}\opb$ and its source is then determined cyclically. 

As the notation is rather complex we will give one example with all variables written down and then convert to a shorthand for convenience.  We define $H_{11}$ by

$H_{11}(a_{12}, \ldots, a_{1n_1}, \ldots, a_{kn_1}, b_0) = $ \\[0.6em]
\hh{3em}$F_{11}\big(a_{12}, \ldots, a_{1n_1}, G_1(F_{20}(a_{21}, \ldots, a_{2n_2}), \ldots, F_{k0}(a_{k1}, \ldots, a_{kn_k}), b_0)\big)$

\vv{1em}\noi where each $a_{ij} \in A_{ij}$ and $b_0 \in B_0$.  We think this is clearer if we do not write the variables explicitly, giving
\[H_{11} = F_{11}\big(\ \underline{\hh{4em}} \ , G_1(F_{20}, \ldots, F_{k0}, \bullet)\big).\]
Here the long line indicates a string of variables and a dot indicates a single variable.  From the sources and targets of all the relevant functors it is unambiguous what the variables need to be, though somewhat tedious to write them out.  The remaining functors $H_{ij}$ can then be written like this:
\[\begin{array}{lcl}
  H_{12} &=& F_{12}\big(\ \underline{\hh{3.8em}} \ , G_1(F_{20}, \ldots, F_{k0}, \bullet), \bullet \big) \\
  H_{13} &=& F_{13}\big(\ \underline{\hh{3.2em}} \ , G_1(F_{20}, \ldots, F_{k0}, \bullet), \bullet, \bullet \big) \\
  H_{14} &=& F_{14}\big(\ \underline{\hh{2.6em}} \ , G_1(F_{20}, \ldots, F_{k0}, \bullet), \bullet, \bullet, \bullet \big) \\
\hh{2em} \vdots \\
  H_{1n_1} &=& F_{1n_1}\big(G_1(F_{20}, \ldots, F_{k0}, \bullet), \ \underline{\hh{4em}} \ \big) \\[8pt]
H_{21} &=& F_{21}\big(\ \underline{\hh{4em}} \ , G_2(F_{30}, \ldots, F_{k0}, \bullet, F_{10})\big) \\
H_{22} &=& F_{22}\big(\ \underline{\hh{3.2em}} \ , G_2(F_{30}, \ldots, F_{k0}, \bullet, F_{10}), \bullet\big)\\
H_{23} &=& F_{23}\big(\ \underline{\hh{2.6em}} \ , G_2(F_{30}, \ldots, F_{k0}, \bullet, F_{10}), \bullet, \bullet \big)\\
\hh{2em} \vdots \\
H_{2n_2} &=& F_{2n_2}\big(G_2(F_{30}, \ldots, F_{k0}, \bullet, F_{10}), \ \underline{\hh{4em}} \ \big)\\
\hh{2em} \vdots \\
H_{k1} &=& F_{k1}\big(\ \underline{\hh{4em}} \ , G_k(\bullet, F_{10}, \ldots, F_{k-1,0}) \big) \\
\hh{2em} \vdots \\
H_{kn_k} &=& F_{kn_k}\big(G_k(\bullet, F_{10}, \ldots, F_{k-1,0}),\ \underline{\hh{4em}} \ \big) \\
  \end{array}\]
It remains to exhibit the adjunctions required, which will take the form of the following isomorphisms.
\[\begin{array}{lcl}
   A_{00}\big( H_{00}(\ha_{00}), a_{00}\big) & \cong &  A_{11}\big( H_{11}(\ha_{11}), a_{11}\big) \\
& \cong &  A_{12}\big( H_{12}(\ha_{12}), a_{12}\big) \\
&& \hh{2em} \vdots \\
& \cong &  A_{kn_k}\big( H_{kn_k}(\ha_{kn_k}), a_{kn_k}\big) \\
  \end{array}\]
The schematic diagram in Table~\ref{bigtable} indicates which adjunctions of $F_{ij}$'s and $G_l$'s are involved with each of the adjunctions for the $H_{ij}$'s.  The vertical arrows indicate individual adjunctions.
\begin{table}[ht]
\caption{\label{bigtable} Individual adjunctions forming composite multivariable adjunctions}
\[\begin{array}{c|cccccc}
 H_{00} & G_{0} & F_{10} & F_{20} & F_{30} & \cdots & F_{k0} \\[-5pt]
        & \updo & \updo  & \\
 H_{11} & G_{1} & F_{11} & F_{20} & F_{30} & \cdots & F_{k0} \\[-5pt]
        &       & \updo  &        &        &        &        \\
 H_{13} & G_{1} & F_{13} & F_{20} & F_{30} & \cdots & F_{k0} \\[-5pt]
        &       & \updo  &        &        &        &        \\
 H_{14} & G_{1} & F_{14} & F_{20} & F_{30} & \cdots & F_{k0} \\[-5pt]
 \vdots \\
        &       & \updo  &        &        &        &        \\
 H_{1n_1} & G_{1} & F_{1n_1} & F_{20} & F_{30} & \cdots & F_{k0} \\[-5pt]
        & \updo & \updo  & \updo  &        &        &        \\
 H_{21} & G_{2} & F_{11} & F_{21} & F_{30} & \cdots & F_{k0} \\[-5pt]
        &       &        & \updo  &        &        &        \\
 H_{22} & G_{2} & F_{11} & F_{22} & F_{30} & \cdots & F_{k0} \\[-5pt]
        &       &        & \updo  &        &        &        \\
 H_{23} & G_{2} & F_{11} & F_{23} & F_{30} & \cdots & F_{k0} \\[-5pt]
\vdots \\
        &       &        & \updo  &        &        &        \\
 H_{2n_2} & G_{2} & F_{11} & F_{2n_2} & F_{30} & \cdots & F_{k0} \\[-5pt]
        & \updo  &        & \updo  & \updo  &        &        \\
 H_{31} & G_{3} & F_{11} & F_{21} & F_{31} & \cdots & F_{k0} \\[-5pt]
        &       &        &        & \updo  &        &        \\
 H_{32} & G_{3} & F_{11} & F_{21} & F_{32} & \cdots & F_{k0} \\[-5pt]
\vdots \\
        & \updo &        &        & \updo  &        &        \\
 H_{k1} & G_{k} & F_{11} & F_{21} & F_{31} & \cdots & F_{k1} \\[-5pt]
        &       &        &        &        &        & \updo  \\
 H_{k2} & G_{k} & F_{11} & F_{21} & F_{31} & \cdots & F_{k2} \\[-5pt]
\vdots \\
        &       &        &        &        &        & \updo  \\
 H_{kn_k} & G_{k} & F_{11} & F_{21} & F_{31} & \cdots & F_{kn_k} \\[-5pt]
        & \updo &        &        &        &        & \updo  \\
 H_{00} & G_{0} & F_{10} & F_{20} & F_{30} & \cdots & F_{k0} \\[-5pt]
        &       &        &        &        &        &        \\
\end{array}
\]
\end{table}

This is much easier to construct formally using Theorem~\ref{theoremA}: we just need to exhibit mutual left \emph{1-variable} adjoints for each of the $m$ functors obtained from $H_{00}$ by fixing all but one of the variables. Now fixing every variable except $a_{ij}$ in the functor
\[H_{11} = G_0(F_{10}, \ldots, F_{k0})\]
we construct a mutual left adjoint using
\begin{itemize}
 \item the mutual left adjoint for $F_{i0}$ with all but the $j$th variable fixed, and
\item the mutual left adjoint for $G_0$ with all but the $i$th variable fixed.
\end{itemize}
These compose to give the adjoint required.  We can depict this schematically as follows.  We depict the latter as
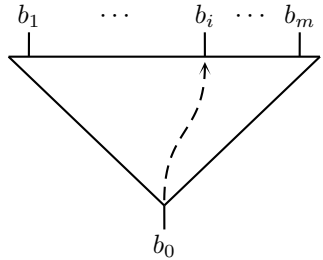
\[
\psset{linewidth=0.8pt,unit=0.9mm,labelsep=1pt,nodesep=0pt}
\pspicture(0,-10)(60,30)

\rput(0,28){\rnode{y1}{$b_1$}}  
\rput(13,28){$\cdots$}  
\rput(26,28){\rnode{yi}{$b_i$}}  
\rput(33,28){$\cdots$}  
\rput(40,28){\rnode{ym}{$b_m$}}  

\rput(0,22){\rnode{y1b}{$$}}  
\rput(26,22){\rnode{yib}{$$}}  
\rput(40,22){\rnode{ymb}{$$}}  

\rput(-3,22){\rnode{ystart}{$$}}  
\rput(43,22){\rnode{yend}{$$}}  

\rput(20,0){\rnode{y0b}{$$}}  %\rput(11,9){$g$}
\rput(20,-6){\rnode{y0}{$b_0$}}

\psset{nodesepA=2pt,nodesepB=0pt}
\ncline{-}{y1}{y1b} 
\ncline{-}{yi}{yib} 
\ncline{-}{ym}{ymb} 
\ncline{-}{y0}{y0b}

\psset{nodesep=0pt}

\ncline{-}{ystart}{yend} 
\ncline{-}{yend}{y0b} 
\ncline{-}{y0b}{ystart} 

\nccurve[linestyle=dashed,nodesep=2pt,angleA=90,angleB=-90,ncurvA=1,ncurvB=1]{->}{y0b}{yib}\naput{{\scriptsize $$}}

\endpspicture
\]
so then composing this with the former looks like the diagram below, where $F_{i0}$ and $G_0$ are the multifunctors pointing downwards, and the 1-variable left adjoint is indicated as the dotted arrow pointing upwards.
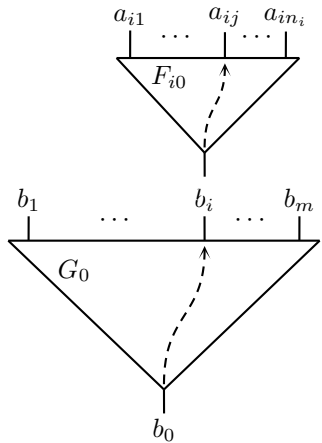
\[
\psset{linewidth=0.8pt,unit=0.9mm,labelsep=1pt,nodesep=0pt,dash=3.5pt 2.5pt}
\pspicture(15,-5)(60,70)

\rput(50,20){\pspicture(0,-10)(60,30)

\rput(0,28){\rnode{y1}{$b_1$}}  
\rput(13,25){$\cdots$}  
\rput(26,28){\rnode{yi}{$b_i$}}  
\rput(33,25){$\cdots$}  
\rput(40,28){\rnode{ym}{$b_m$}}  

\rput(0,22){\rnode{y1b}{$$}}  
\rput(26,22){\rnode{yib}{$$}}  
\rput(40,22){\rnode{ymb}{$$}}  

\rput(-3,22){\rnode{ystart}{$$}}  
\rput(43,22){\rnode{yend}{$$}}  

\rput(20,0){\rnode{y0b}{$$}}  %\rput(11,9){$g$}
\rput(20,-6){\rnode{y0}{$b_0$}}

\psset{nodesepA=2pt,nodesepB=0pt}
\ncline{-}{y1}{y1b} 
\ncline{-}{yi}{yib} 
\ncline{-}{ym}{ymb} 
\ncline{-}{y0}{y0b}

\psset{nodesep=0pt}

\ncline{-}{ystart}{yend} 
\ncline{-}{yend}{y0b} 
\ncline{-}{y0b}{ystart} \nbput[npos=0.7]{{$G_{0}$}}

\nccurve[linestyle=dashed,nodesep=2pt,angleA=90,angleB=-90,ncurvA=1,ncurvB=1]{->}{y0b}{yib}\naput{{\scriptsize $$}}

\endpspicture}

\rput(56,55){\pspicture(0,-10)(60,30)

\rput(9,20){\rnode{y1}{$a_{i1}$}}  
\rput(16,17){$\cdots$}  
\rput(23,20){\rnode{yi}{$a_{ij}$}}  
\rput(28,17){$\cdots$}  
\rput(32,20){\rnode{ym}{$a_{in_i}$}}  

\rput(9,14){\rnode{y1b}{$$}}  
\rput(23,14){\rnode{yib}{$$}}  
\rput(32,14){\rnode{ymb}{$$}}  

\rput(7,14){\rnode{ystart}{$$}}  
\rput(34,14){\rnode{yend}{$$}}  

\rput(20,0){\rnode{y0b}{$$}}  %\rput(11,9){$g$}
\rput(20,-6){\rnode{y0}{\textcolor{white}{$b_0$}}}

\psset{nodesepA=2pt,nodesepB=0pt}
\ncline{-}{y1}{y1b} 
\ncline{-}{yi}{yib} 
\ncline{-}{ym}{ymb} 
\ncline{-}{y0}{y0b}

\psset{nodesep=0pt}

\ncline{-}{ystart}{yend} 
\ncline{-}{yend}{y0b} 
\ncline{-}{y0b}{ystart} \nbput[npos=0.65]{{$F_{i0}$}}

\nccurve[linestyle=dashed,nodesep=2pt,angleA=90,angleB=-90,ncurvA=1,ncurvB=1]{->}{y0b}{yib}\naput{{\scriptsize $$}}

\endpspicture}

\endpspicture
\]
That is, starting from the functor 
\[F_{i0} \: A_{i1} \times \cdots \times A_{in_i} \tra A_{i0}\opb \]
we fix all but the $j$th variable and have a mutual left adjoint, that is right adjoint for
\[F_{i0}\opb(a_{11}, \ldots, a_{i,j-1}, \uscore, a_{i,j+1}, \ldots, a_{in_i}) \: A_{ij}\opb \tra A_{i0}\]
which is 
\[F_{ij}(a_{i,j+1}, \ldots, a_{i0}, \uscore, a_i, \ldots, a_{i,j-1}) \: A_{i0} \tra A_{ij}\opb.\]
Also consider
\[G_0\opb ( b_1, \ldots, b_{i-1}, \uscore, b_{i+1}, \ldots, b_k) \: B_i\opb \tra B_0\]
with its right adjoint
\[G_i(b_{i+1}, \ldots, b_k, \uscore, b_1, \ldots, b_{i-1}) \: B_0 \tra B_i\opb = A_{i0}.\]
Now we simply compose the functors
\[B_0 \tra A_{i0} \tra A_{ij}\opb\]
setting each
\[b_q = F_{q0}(a_{q1}, \ldots, a_{qn_q}).\]
Using the previous shorthand this is the composite
\[F_{ij}\big( \ \underline{\hh{3em}} \ , G_i (F_{i+1, 0}, \ldots, F_{k0}, \bullet, F_{10}, \ldots, F_{i0}), \ \underline{\hh{3em}} \ \big).\]
This completes the construction of composition.  Identities are given by identity adjunctions, which obviously satisfy unit conditions.  Associativity follows from associativity of composition of $n$-variable functors (with one another) and of 1-variable adjunctions (with one another).  
\end{prf}

\subsubsection*{Special cases}

\begin{enumerate}
 \item If any $n_i=0$ this amounts to fixing the $i$th variable of $G$.  If all but 1 of the $n_i$ is 0 then we have fixed every variable except one, and if we do this for each $n_i$ in turn we have effectively characterised the composite multivariable adjunction by producing the necessary 1-variable adjunctions as in Theorem~\ref{theoremA}.

\item If we compose with the identity adjunction (as a 1-ary adjunction) for all but one of the $i$'s, we have effectively composed in just one position.

\item If we take every $n_i=1$ or $k=1$ this says we can compose an $n$-variable adjunction with 1-variable adjunctions (pre- or post-) to get a new $n$-variable adjunction; this example is mentioned for composing 2-ary with 1-ary adjunctions in \cite{rie2}.

\end{enumerate}

\subsection{Multivariable mates}

We now give a multivariable version of the calculus of mates.  As for the adjunctions, we start with the 2-variable case and proceed inductively.

\begin{proposition}\label{theoremC}
 Suppose we have two 2-variable left adjunctions of functors
\[\begin{array}{ccc}
   F \: A \times B \tra C\opb & \mbox{and} & F' \: A' \times B' \tra {C'}\opb \\
G \: B \times C \tra A \opb && G' \: B' \times C' \tra {A'}\opb \\
H \: C \times A \tra B\opb && H' \: C' \times A' \tra {B'}\opb,
  \end{array}\]
together with functors
\[\begin{array}{c}
   S\: A \tra A' \\
T \: B \tra B' \\
U \: C \tra C'
  \end{array}\]
and a natural transformation

\[\psset{unit=0.08cm,labelsep=0pt,nodesep=3pt}
\pspicture(0,-3)(30,22)

% a1 a2 up right
% b1 b2

\rput(0,20){\rnode{a1}{$A\times B$}} % top left
\rput(30,20){\rnode{a2}{$A' \times B'$}} % top right

\rput(0,0){\rnode{b1}{$C\opb$}}   % bottom left
\rput(30,0){\rnode{b2}{${C'}\opb$}}  % bottom right

\psset{nodesep=3pt,labelsep=2pt,arrows=->}
\ncline{a1}{a2}\naput{{\scriptsize $S \times T$}} % top
\ncline{b1}{b2}\nbput{{\scriptsize $U\opb$}} % bottom
\ncline{a1}{b1}\nbput{{\scriptsize $F$}} % left
\ncline{a2}{b2}\naput{{\scriptsize $F'$}} % right

\psset{labelsep=1.5pt}
\pnode(12,7){c1}
\pnode(18,13){c2}
\ncline[doubleline=true]{c1}{c2} \nbput[npos=0.4]{{\scriptsize $\alpha$}}

\endpspicture\]
with components
\[\alpha_{a,b} \: F'(Sa, Tb) \tra UF(a,b).\]
Then for each $b \in B$ we have a natural transformation
\[\psset{unit=0.08cm,labelsep=0pt,nodesep=3pt}
\pspicture(0,-3)(30,22)

% a1 a2 up right
% b1 b2

\rput(0,20){\rnode{a1}{$A$}} % top left
\rput(30,20){\rnode{a2}{${A'}$}} % top right

\rput(0,0){\rnode{b1}{$C\opb$}}   % bottom left
\rput(30,0){\rnode{b2}{${C'}\opb$}}  % bottom right

\psset{nodesep=3pt,labelsep=2pt,arrows=->}
\ncline{a1}{a2}\naput{{\scriptsize $S$}} % top
\ncline{b1}{b2}\nbput{{\scriptsize $U\opb$}} % bottom
\ncline{a1}{b1}\nbput{{\scriptsize $F(\uscore, b)$}} % left
\ncline{a2}{b2}\naput{{\scriptsize $F'(\uscore, Tb)$}} % right

\psset{labelsep=1.5pt}
\pnode(12,7){c1}
\pnode(18,13){c2}
\ncline[doubleline=true]{c1}{c2} \nbput[npos=0.4]{{\scriptsize $\alpha_{\uscore, b}$}}

\endpspicture\]
with mate

\[\psset{unit=0.08cm,labelsep=0pt,nodesep=3pt}
\pspicture(0,-3)(30,22)

% a1 a2 up right medium
% b1 b2

\rput(0,20){\rnode{a1}{$C$}} % top left
\rput(30,20){\rnode{a2}{$C'$}} % top right

\rput(0,0){\rnode{b1}{$A\opb$}}   % bottom left
\rput(30,0){\rnode{b2}{${A'}\opb.$}}  % bottom right

\psset{nodesep=3pt,labelsep=2pt,arrows=->}
\ncline{a1}{a2}\naput{{\scriptsize $U$}} % top
\ncline{b1}{b2}\nbput{{\scriptsize $S\opb$}} % bottom
\ncline{a1}{b1}\nbput{{\scriptsize $G(b,\uscore)$}} % left
\ncline{a2}{b2}\naput{{\scriptsize ${G}'(Tb,\uscore)$}} % right

\psset{labelsep=-1pt}
\pnode(12,7){c1}
\pnode(18,13){c2}
\ncline[doubleline=true]{c1}{c2} \nbput[npos=0.5,labelsep=2pt]{{\scriptsize $\overline{{\alpha}_{\uscore,b}}$}}

\endpspicture\]
Then in fact the components $(\overline{{\alpha}_{\uscore,b}})_c$ are the components of a natural transformation

\[\psset{unit=0.08cm,labelsep=0pt,nodesep=3pt}
\pspicture(0,-3)(30,22)

% a1 a2 up right
% b1 b2

\rput(0,20){\rnode{a1}{$B\times C$}} % top left
\rput(30,20){\rnode{a2}{$B' \times C'$}} % top right

\rput(0,0){\rnode{b1}{$A\opb$}}   % bottom left
\rput(30,0){\rnode{b2}{${A'}\opb.$}}  % bottom right

\psset{nodesep=3pt,labelsep=2pt,arrows=->}
\ncline{a1}{a2}\naput{{\scriptsize $T \times U$}} % top
\ncline{b1}{b2}\nbput{{\scriptsize $S\opb$}} % bottom
\ncline{a1}{b1}\nbput{{\scriptsize $G$}} % left
\ncline{a2}{b2}\naput{{\scriptsize $G'$}} % right

\psset{labelsep=1.5pt}
\pnode(12,7){c1}
\pnode(18,13){c2}
\ncline[doubleline=true]{c1}{c2} \nbput[npos=0.4]{{\scriptsize $\bar{\alpha}$}}

\endpspicture\]

Dually if we start with 2-variable \emph{right} adjunctions then the result holds with all the \nts\ pointing in the opposite direction as below.

\[\psset{unit=0.08cm,labelsep=0pt,nodesep=3pt}
\pspicture(0,-3)(30,22)

% a1 a2 up right
% b1 b2

\rput(0,20){\rnode{a1}{$A\times B$}} % top left
\rput(30,20){\rnode{a2}{$A' \times B'$}} % top right

\rput(0,0){\rnode{b1}{$C\opb$}}   % bottom left
\rput(30,0){\rnode{b2}{${C'}\opb$}}  % bottom right

\psset{nodesep=3pt,labelsep=2pt,arrows=->}
\ncline{a1}{a2}\naput{{\scriptsize $S \times T$}} % top
\ncline{b1}{b2}\nbput{{\scriptsize $U\opb$}} % bottom
\ncline{a1}{b1}\nbput{{\scriptsize $F$}} % left
\ncline{a2}{b2}\naput{{\scriptsize $F'$}} % right

\psset{labelsep=1.5pt}
\pnode(18,13){c1}
\pnode(12,7){c2}
\ncline[doubleline=true]{c1}{c2} \nbput[npos=0.4]{{\scriptsize $\alpha$}}

\endpspicture\]

\end{proposition}

\begin{prf}
 We just need to check that the components $\bar{\alpha}_{b,c}= (\overline{{\alpha}_{\uscore,b}})_c$ are natural in $b$; \emph{a priori} they are natural in $c$.  We use the fact that $\alpha$ is natural in $b$.  As with Theorem~\ref{theoremA} the proof is possible by a 1-dimensional diagram chase, but we provide a 2-categorical proof as it is quite aesthetically pleasing.

Now the \nt\ $\overline{{\alpha}_{\uscore,b}}$ is given by the following composite 
% a1 a2 a3  a4
%       b3  b4
% c1    c3  c4  c5  c6

%%%%%%%% 
\[
\psset{unit=0.1cm,labelsep=2pt,nodesep=2pt}
\pspicture(-5,-20)(100,35)
\rput(0,30){\rnode{a1}{$C$}}
\rput(18,30){\rnode{a2}{$B \times C$}}
\rput(36,30){\rnode{a3}{$A\opb$}}
\rput(62,30){\rnode{a4}{${A'}\opb$}}

\rput(36,15){\rnode{b3}{$A\opb \ttimes B \opb$}}
\rput(62,15){\rnode{b4}{${A'}\opb \ttimes {B'}\opb$}}

\rput(11,-5){\rnode{c1}{$B \times C$}}
\rput(36,0){\rnode{c3}{$C$}}
\rput(62,0){\rnode{c4}{$C'$}}
\rput(80,0){\rnode{c5}{$B'\times C'$}}
\rput(98,0){\rnode{c6}{${{A'}\opb}$}}

% horizontals
\ncline{->}{a1}{a2} \naput{{\scriptsize $b \ttimes 1$}}
\ncline{->}{a2}{a3} \naput{{\scriptsize $G$}}
\ncline{->}{a3}{a4} \naput{{\scriptsize $S\opb$}}

\ncline{->}{b3}{b4} \naput[labelsep=1pt]{{\scriptsize $S\opb \ttimes T\opb$}}

\ncline{->}{c3}{c4} \nbput{{\scriptsize $U$}}
\ncline{->}{c4}{c5} \nbput{{\scriptsize $Tb \ttimes 1$}}
\ncline{->}{c5}{c6} \nbput{{\scriptsize $G'$}}

%verticals

\psset{labelsep=1pt}{
\ncline{->}{a3}{b3} \naput{{\scriptsize $1 \ttimes b\opb$}}
\ncline{->}{b3}{c3} \naput{{\scriptsize $F\opb$}}
\ncline{->}{a4}{b4} \nbput{{\scriptsize $1 \ttimes Tb\opb$}}
\ncline{->}{b4}{c4} \nbput{{\scriptsize ${F'}\opb$}}}

%diags
\ncline{->}{a1}{c3} \nbput{{\scriptsize $1$}}
\ncline{->}{a4}{c6} \naput{{\scriptsize $1$}}

%weirds
\nccurve[angleA=270,angleB=120]{->}{a1}{c1} \nbput{{\scriptsize $b\ttimes 1$}}
\nccurve[angleA=340,angleB=210]{->}{c1}{c5} \nbput{{\scriptsize $T \times U$}}

% 2-cells

\psset{doubleline=true}
\rput[c](21,19){\pcline{->}(4,4)(0,0) \naput[labelsep=1pt]{{\scriptsize 	     $\varepsilon_{b}$}}}
\rput[c](73,7){\pcline{->}(4,4)(0,0) \naput[labelsep=-1pt,npos=0.3]{{\scriptsize     	     $\eta'_{Tb}$}}}
\rput[c](46,5){\pcline{->}(4,4)(0,0) \naput[labelsep=-1pt,npos=0.3]	         	{{\scriptsize $\alpha\opb$}}}

\endpspicture
\]
taking care over the direction as the target is in ${A'}\opb$. Again we use the fact that a morphism $b_1 \tmap{F} b_2$ in $B$ corresponds to a \nt\ 
\[
\psset{unit=0.1cm,labelsep=2pt,nodesep=3pt}
\pspicture(20,20)

\rput(0,10){\rnode{a1}{$1$}}  % named node, with something placed there
\rput(20,10){\rnode{a2}{$B$}}  % named node, with something placed there

\pcline[doubleline=true]{->}(10,13)(10,7) \naput{{\scriptsize $f$}}

\ncarc[arcangle=45]{->}{a1}{a2}\naput{{\scriptsize $b_1$}}
\ncarc[arcangle=-45]{->}{a1}{a2}\nbput{{\scriptsize $b_2$}}

\endpspicture
\]
thus to check that the components 

\[
\psset{unit=0.1cm,labelsep=1pt,nodesep=2pt}
\pspicture(30,20)

%   b1
% a1  a2
%   c1

\rput(0,10){\rnode{a1}{$C$}}  % named node, with something placed there
\rput(30,10){\rnode{a2}{${A'}\opb$}}  % named node, with something placed there

\rput(15,17){\rnode{b1}{$B \ttimes C$}}
\rput(15,3){\rnode{c1}{$B \ttimes C$}}

\nccurve[angleA=50,angleB=185,ncurvA=0.1,ncurvB=0.9,nodesepB=1pt]{->}{a1}{b1}\naput{{\scriptsize $b \ttimes 1$}}
\nccurve[angleA=-5,angleB=130,ncurvB=0.1,ncurvA=0.9,nodesepB=1pt]{->}{b1}{a2}\naput{{\scriptsize $S\opb G$}}
\nccurve[angleA=-50,angleB=-185,ncurvA=0.1,ncurvB=0.9,nodesepB=1pt]{->}{a1}{c1}\nbput{{\scriptsize $b \ttimes 1$}}
\nccurve[angleA=5,angleB=-130,ncurvB=0.1,ncurvA=0.9,nodesepB=1pt]{->}{c1}{a2}\nbput{{\scriptsize $G'(T\ttimes U)$}}

\psset{doubleline=true}
\rput[c](15,7){\pcline{->}(0,6)(0,0) \naput[labelsep=1pt]{{\scriptsize $\bar{\alpha}_{b,\uscore}$}}}

\endpspicture
\]
are natural in $b$ we show that for all $b_1 \tmap{f} b_2$

\[
\psset{unit=0.1cm,labelsep=1pt,nodesep=2pt}
\pspicture(70,23)

\rput[b](0,0){
\pspicture(30,20)

%   b1
% a1  a2
%   c1

\rput(0,10){\rnode{a1}{$C$}}  % named node, with something placed there
\rput(40,10){\rnode{a2}{${A'}\opb$}}  % named node, with something placed there

\rput(20,19){\rnode{b1}{$B \ttimes C$}}
\rput(20,1){\rnode{c1}{$B \ttimes C$}}

%%%%%%%%%%%%%%%%%%%%% top left curves

\rput(17,19){\rnode{b}{}} %fake node

\ncarc[arcangle=38,nodesepB=4pt]{->}{a1}{b}\naput[npos=0.5,labelsep=1pt]{{\scriptsize $b_1 \ttimes 1$}}
\ncarc[arcangle=-38,nodesepB=4pt]{->}{a1}{b}\nbput[npos=0.4,labelsep=1pt]{{\scriptsize $b_2 \ttimes 1$}}

% top right curve
\ncline[nodesepB=1pt]{->}{b1}{a2}\naput[npos=0.4]{{\scriptsize $S\opb G$}}

% bottom left curve

\rput(17,1){\rnode{c}{}}

\ncline[nodesepB=4pt]{->}{a1}{c}\nbput[labelsep=0pt]{{\scriptsize $b_2 \ttimes 1$}}

% bottom right curve

\ncline[nodesepB=1pt]{->}{c1}{a2}\nbput{{\scriptsize $G'(T\ttimes U)$}}

% 2-cells

\psset{doubleline=true}
\rput[c](20,7){\pcline{->}(0,6)(0,0) \naput[labelsep=1pt]{{\scriptsize $\bar{\alpha}_{b,\uscore}$}}}

\rput[c](5.5,11.5){\pcline{->}(0,3.8)(2.5,0.5) \naput[labelsep=1pt]{{\scriptsize $f \ttimes 1$}}}

\endpspicture}

%%%%%%%%%%%%%%
\rput(37,10){$=$}
%%%%%%%%%%%%%%

\rput[b](63,0){
\pspicture(30,20)

%   b1
% a1  a2
%   c1

\rput(0,10){\rnode{a1}{$C$}}  % named node, with something placed there
\rput(40,10){\rnode{a2}{${A'}\opb.$}}  % named node, with something placed there

\rput(20,19){\rnode{b1}{$B \ttimes C$}}
\rput(21,1){\rnode{c1}{$B \ttimes C$}}

% top left curve

%\rput(17,18){\rnode{b}{}} %fake node

\ncline[nodesepB=1pt]{->}{a1}{b1}\naput{{\scriptsize $b_1 \ttimes 1$}}

% top right curve

\ncline[nodesepB=1pt]{->}{b1}{a2}\naput[npos=0.4]{{\scriptsize $S\opb G$}}

% bottom left curves

\rput(17,1){\rnode{c}{}}

\ncarc[arcangle=38,nodesepB=4pt]{->}{a1}{c}\naput[labelsep=0pt,npos=0.4]{{\scriptsize $b_1 \ttimes 1$}}
\ncarc[arcangle=-38,nodesepB=4pt]{->}{a1}{c}\nbput[labelsep=0pt]{{\scriptsize $b_2 \ttimes 1$}}

% bottom right curve

\ncline[nodesepB=1pt]{->}{c1}{a2}\nbput{{\scriptsize $G'(T\ttimes U)$}}

\psset{doubleline=true}
\rput[c](21,7){\pcline{->}(0,6)(0,0) \naput[labelsep=1pt]{{\scriptsize $\bar{\alpha}_{b,\uscore}$}}}

\rput[c](5.7,4.8){\pcline{->}(2.7,3.5)(0,0) \naput[labelsep=0pt]{{\scriptsize $f \ttimes 1$}}}

\endpspicture}

\endpspicture
\]

 We have 
% a1 a2 a3  a4
%       b3  b4
% c1    c3  c4  c5  c6
%%%% step 1
\[
\psset{unit=0.1cm,labelsep=2pt,nodesep=2pt}
\pspicture(-5,-10)(100,40)
\rput(-4,30){\rnode{a1}{$C$}}
\rput(16,30){\rnode{a2}{$B \times C$}}
\rput(36,30){\rnode{a3}{$A\opb$}}
\rput(62,30){\rnode{a4}{${A'}\opb$}}

\rput(36,15){\rnode{b3}{$A\opb \ttimes B \opb$}}
\rput(62,15){\rnode{b4}{${A'}\opb \ttimes {B'}\opb$}}

\rput(11,-5){\rnode{c1}{$B \times C$}}
\rput(36,0){\rnode{c3}{$C$}}
\rput(62,0){\rnode{c4}{$C'$}}
\rput(82,0){\rnode{c5}{$B'\times C'$}}
\rput(102,0){\rnode{c6}{${{A'}\opb}$}}

% horizontals
%\ncline{->}{a1}{a2} \naput{{\scriptsize $b \ttimes 1$}}

\ncarc[arcangle=40]{->}{a1}{a2}\naput{{\scriptsize $b_1 \ttimes 1$}}
\ncarc[arcangle=-40]{->}{a1}{a2}\nbput{{\scriptsize $b_2 \ttimes 1$}}
{\psset{doubleline=true}
\rput[c](3,28){\pcline{->}(0,4)(0,0) \naput[labelsep=1pt]{{\scriptsize $f\ttimes 1$}}}}

% \ncarc[arcangle=40]{->}{a1}{a2}\naput{{\scriptsize $$}}
% \ncarc[arcangle=-40]{->}{a1}{a2}\nbput{{\scriptsize $$}}

\ncline{->}{a2}{a3} \naput{{\scriptsize $G$}}
\ncline{->}{a3}{a4} \naput{{\scriptsize $S\opb$}}

\ncline{->}{b3}{b4} \naput[labelsep=1pt]{{\scriptsize $S\opb \ttimes T\opb$}}

\ncline{->}{c3}{c4} \nbput{{\scriptsize $U$}}
\ncline{->}{c4}{c5} \nbput{{\scriptsize $Tb_2 \ttimes 1$}}
\ncline{->}{c5}{c6} \nbput{{\scriptsize $G'$}}

%verticals

\psset{labelsep=1pt}{
\ncline{->}{a3}{b3} \naput{{\scriptsize $1 \ttimes b_2\opb$}}
\ncline{->}{b3}{c3} \naput{{\scriptsize $F\opb$}}
\ncline{->}{a4}{b4} \nbput{{\scriptsize $1 \ttimes Tb_2\opb$}}
\ncline{->}{b4}{c4} \nbput{{\scriptsize ${F'}\opb$}}}

%diags
\ncarc[arcangle=-20]{->}{a1}{c3} \nbput{{\scriptsize $1$}}
\ncline{->}{a4}{c6} \naput{{\scriptsize $1$}}

%weirds
\nccurve[angleA=270,angleB=120]{->}{a1}{c1} \nbput{{\scriptsize $b_2\ttimes 1$}}
\nccurve[angleA=340,angleB=210]{->}{c1}{c5} \nbput{{\scriptsize $T \times U$}}

% 2-cells

\psset{doubleline=true}
\rput[c](16,15){\pcline{->}(4,4)(0,0) \naput[labelsep=1pt]{{\scriptsize 	     $\varepsilon_{b_2}$}}}
\rput[c](73,7){\pcline{->}(4,4)(0,0) \naput[labelsep=-1pt,npos=0.3]{{\scriptsize $\eta'_{Tb_2}$}}}
\rput[c](46,5){\pcline{->}(4,4)(0,0) \naput[labelsep=-1pt,npos=0.3]	       {{\scriptsize $\alpha\opb$}}}

\endpspicture
\]
%%%% step 2
% a1 a2 a3  a4
%       b3  b4
% c1    c3  c4  c5  c6

\[
\psset{unit=0.1cm,labelsep=2pt,nodesep=2pt}
\pspicture(-5,-10)(100,35)

\rput(-12,6){$=$}

\rput(-4,30){\rnode{a1}{$C$}}
\rput(16,30){\rnode{a2}{$B \times C$}}
\rput(36,30){\rnode{a3}{$A\opb$}}
\rput(62,30){\rnode{a4}{${A'}\opb$}}

\rput(36,15){\rnode{b3}{$A\opb \ttimes B \opb$}}
\rput(62,15){\rnode{b4}{${A'}\opb \ttimes {B'}\opb$}}

\rput(11,-5){\rnode{c1}{$B \times C$}}
\rput(36,0){\rnode{c3}{$C$}}
\rput(62,0){\rnode{c4}{$C'$}}
\rput(82,0){\rnode{c5}{$B'\times C'$}}
\rput(102,0){\rnode{c6}{${{A'}\opb}$}}

% horizontals
\ncline{->}{a1}{a2} \naput{{\scriptsize $b_1 \ttimes 1$}}

% \ncarc[arcangle=40]{->}{a1}{a2}\naput{{\scriptsize $$}}
% \ncarc[arcangle=-40]{->}{a1}{a2}\nbput{{\scriptsize $$}}

\ncline{->}{a2}{a3} \naput{{\scriptsize $G$}}
\ncline{->}{a3}{a4} \naput{{\scriptsize $S\opb$}}

\ncline{->}{b3}{b4} \naput[labelsep=1pt]{{\scriptsize $S\opb \ttimes T\opb$}}

\ncline{->}{c3}{c4} \nbput{{\scriptsize $U$}}
\ncline{->}{c4}{c5} \nbput{{\scriptsize $Tb_2 \ttimes 1$}}
\ncline{->}{c5}{c6} \nbput{{\scriptsize $G'$}}

%verticals

\psset{labelsep=1pt}{

%\ncline{->}{a3}{b3} \naput{{\scriptsize $1 \ttimes b\opb$}}

\ncarc[arcangle=43]{->}{a3}{b3}\naput{{\scriptsize $1 \ttimes b_2\opb$}}
\ncarc[arcangle=-43]{->}{a3}{b3}\nbput{{\scriptsize $1 \ttimes b_1\opb$}}
{\psset{doubleline=true}
\rput[c](33.5,21){\pcline{->}(4,0)(0,0) \nbput[npos=0.35,labelsep=1pt]{{\scriptsize $1\ttimes \hh{-0.5pt} f\opb$}}}}

\ncline{->}{b3}{c3} \naput{{\scriptsize $F\opb$}}
\ncline{->}{a4}{b4} \nbput{{\scriptsize $1 \ttimes Tb_2\opb$}}
\ncline{->}{b4}{c4} \nbput{{\scriptsize ${F'}\opb$}}}

%diags
\ncarc[arcangle=-20]{->}{a1}{c3} \nbput{{\scriptsize $1$}}
\ncline{->}{a4}{c6} \naput{{\scriptsize $1$}}

%weirds
\nccurve[angleA=270,angleB=120]{->}{a1}{c1} \nbput{{\scriptsize $b_2\ttimes 1$}}
\nccurve[angleA=340,angleB=210]{->}{c1}{c5} \nbput{{\scriptsize $T \times U$}}

% 2-cells

\psset{doubleline=true}
\rput[c](16,15){\pcline{->}(4,4)(0,0) \naput[labelsep=1pt]{{\scriptsize 	     $\varepsilon_{b_1}$}}}
\rput[c](73,7){\pcline{->}(4,4)(0,0) \naput[labelsep=-1pt,npos=0.3]{{\scriptsize $\eta'_{Tb_2}$}}}
\rput[c](46,5){\pcline{->}(4,4)(0,0) \naput[labelsep=-1pt,npos=0.3]	       {{\scriptsize $\alpha\opb$}}}

\endpspicture
\]
% a1 a2 a3  a4
%       b3  b4
% c1    c3  c4  c5  c6
% step 3
\[
\psset{unit=0.1cm,labelsep=2pt,nodesep=2pt}
\pspicture(-5,-10)(100,40)

\rput(-12,6){$=$}

\rput(-4,30){\rnode{a1}{$C$}}
\rput(16,30){\rnode{a2}{$B \times C$}}
\rput(36,30){\rnode{a3}{$A\opb$}}
\rput(62,30){\rnode{a4}{${A'}\opb$}}

\rput(36,15){\rnode{b3}{$A\opb \ttimes B \opb$}}
\rput(62,15){\rnode{b4}{${A'}\opb \ttimes {B'}\opb$}}

\rput(11,-5){\rnode{c1}{$B \times C$}}
\rput(36,0){\rnode{c3}{$C$}}
\rput(62,0){\rnode{c4}{$C'$}}
\rput(82,0){\rnode{c5}{$B'\times C'$}}
\rput(102,0){\rnode{c6}{${{A'}\opb}$}}

% horizontals
\ncline{->}{a1}{a2} \naput{{\scriptsize $b_1 \ttimes 1$}}

% \ncarc[arcangle=40]{->}{a1}{a2}\naput{{\scriptsize $$}}
% \ncarc[arcangle=-40]{->}{a1}{a2}\nbput{{\scriptsize $$}}

\ncline{->}{a2}{a3} \naput{{\scriptsize $G$}}
\ncline{->}{a3}{a4} \naput{{\scriptsize $S\opb$}}

\ncline{->}{b3}{b4} \naput[labelsep=1pt]{{\scriptsize $S\opb \ttimes T\opb$}}

\ncline{->}{c3}{c4} \nbput{{\scriptsize $U$}}
\ncline{->}{c4}{c5} \nbput{{\scriptsize $Tb_2 \ttimes 1$}}
\ncline{->}{c5}{c6} \nbput{{\scriptsize $G'$}}

%verticals

\psset{labelsep=1pt}{

\ncline{->}{a3}{b3} \naput{{\scriptsize $1 \ttimes b_1\opb$}}

\ncline{->}{b3}{c3} \naput{{\scriptsize $F\opb$}}

%\ncline{->}{a4}{b4} \nbput{{\scriptsize $1 \ttimes Tb_2\opb$}}

\ncarc[arcangle=47]{->}{a4}{b4}\naput{{\scriptsize $1 \ttimes Tb_2\opb$}}
\ncarc[arcangle=-47]{->}{a4}{b4}\nbput{{\scriptsize $1 \ttimes Tb_1\opb$}}
{\psset{doubleline=true}
\rput[c](59.5,21){\pcline{->}(4,0)(0,0) \nbput[npos=0.35,labelsep=1pt]{{\scriptsize $1\ttimes T \hh{-1.0pt} f\opb$}}}}

\ncline{->}{b4}{c4} \nbput{{\scriptsize ${F'}\opb$}}}

%diags
\ncarc[arcangle=-20]{->}{a1}{c3} \nbput{{\scriptsize $1$}}
\ncarc[arcangle=20]{->}{a4}{c6} \naput{{\scriptsize $1$}}

%weirds
\nccurve[angleA=270,angleB=120]{->}{a1}{c1} \nbput{{\scriptsize $b_2\ttimes 1$}}
\nccurve[angleA=340,angleB=210]{->}{c1}{c5} \nbput{{\scriptsize $T \times U$}}

% 2-cells

\psset{doubleline=true}
\rput[c](16,15){\pcline{->}(4,4)(0,0) \naput[labelsep=1pt]{{\scriptsize 	     $\varepsilon_{b_1}$}}}
\rput[c](73,7){\pcline{->}(4,4)(0,0) \naput[labelsep=-1pt,npos=0.3]{{\scriptsize $\eta'_{Tb_2}$}}}
\rput[c](46,5){\pcline{->}(4,4)(0,0) \naput[labelsep=-1pt,npos=0.3]	       {{\scriptsize $\alpha\opb$}}}

\endpspicture
\]

% a1 a2 a3  a4
%       b3  b4
% c1    c3  c4  c5  c6
% step 4
\[
\psset{unit=0.1cm,labelsep=2pt,nodesep=2pt}
\pspicture(-5,-10)(100,35)

\rput(-12,6){$=$}

\rput(-4,30){\rnode{a1}{$C$}}
\rput(16,30){\rnode{a2}{$B \times C$}}
\rput(36,30){\rnode{a3}{$A\opb$}}
\rput(62,30){\rnode{a4}{${A'}\opb$}}

\rput(36,15){\rnode{b3}{$A\opb \ttimes B \opb$}}
\rput(62,15){\rnode{b4}{${A'}\opb \ttimes {B'}\opb$}}

\rput(11,-5){\rnode{c1}{$B \times C$}}
\rput(36,0){\rnode{c3}{$C$}}
\rput(62,0){\rnode{c4}{$C'$}}
\rput(82,0){\rnode{c5}{$B'\times C'$}}
\rput(102,0){\rnode{c6}{${{A'}\opb}$}}

% horizontals
\ncline{->}{a1}{a2} \naput{{\scriptsize $b_1 \ttimes 1$}}

% \ncarc[arcangle=40]{->}{a1}{a2}\naput{{\scriptsize $$}}
% \ncarc[arcangle=-40]{->}{a1}{a2}\nbput{{\scriptsize $$}}

\ncline{->}{a2}{a3} \naput{{\scriptsize $G$}}
\ncline{->}{a3}{a4} \naput{{\scriptsize $S\opb$}}

\ncline{->}{b3}{b4} \naput[labelsep=1pt]{{\scriptsize $S\opb \ttimes T\opb$}}

\ncline{->}{c3}{c4} \nbput{{\scriptsize $U$}}

%\ncline{->}{c4}{c5} \nbput{{\scriptsize $Tb_2 \ttimes 1$}}

\ncarc[arcangle=34]{->}{c4}{c5}\naput[labelsep=0pt,npos=0.4]{{\scriptsize $Tb_1 \ttimes 1$}}
\ncarc[arcangle=-34]{->}{c4}{c5}\nbput[labelsep=0pt,npos=0.4]{{\scriptsize $Tb_2 \ttimes 1$}}
{\psset{doubleline=true}
\rput[c](68,-2){\pcline{->}(0,4)(0,0) \naput[npos=0.35,labelsep=1pt]{{\scriptsize $T \hh{-1.0pt} f\ttimes \hh{-0.5pt} 1$}}}}

\ncline{->}{c5}{c6} \nbput{{\scriptsize $G'$}}

%verticals

\psset{labelsep=1pt}{

\ncline{->}{a3}{b3} \naput{{\scriptsize $1 \ttimes b_1\opb$}}

\ncline{->}{b3}{c3} \naput{{\scriptsize $F\opb$}}

\ncline{->}{a4}{b4} \nbput{{\scriptsize $1 \ttimes Tb_1\opb$}}

\ncline{->}{b4}{c4} \nbput{{\scriptsize ${F'}\opb$}}}

%diags
\ncarc[arcangle=-20]{->}{a1}{c3} \nbput{{\scriptsize $1$}}
\ncarc[arcangle=20]{->}{a4}{c6} \naput{{\scriptsize $1$}}

%weirds
\nccurve[angleA=270,angleB=120]{->}{a1}{c1} \nbput{{\scriptsize $b_2\ttimes 1$}}

%fake coord
\rput(83,-2){\rnode{cfake}{$$}}
\nccurve[angleA=340,angleB=220,ncurvB=0.6]{->}{c1}{cfake} \nbput{{\scriptsize $T \times U$}}

% 2-cells

\psset{doubleline=true}
\rput[c](16,15){\pcline{->}(4,4)(0,0) \naput[labelsep=1pt]{{\scriptsize 	     $\varepsilon_{b_1}$}}}
\rput[c](77,11){\pcline{->}(4,4)(0,0) \naput[labelsep=-1pt,npos=0.3]{{\scriptsize $\eta'_{Tb_1}$}}}
\rput[c](46,5){\pcline{->}(4,4)(0,0) \naput[labelsep=-1pt,npos=0.3]	       {{\scriptsize $\alpha\opb$}}}

\endpspicture
\]

% a1 a2 a3  a4
%       b3  b4
% c1    c3  c4  c5  c6
% step 5
\[
\psset{unit=0.1cm,labelsep=2pt,nodesep=2pt}
\pspicture(-5,-10)(100,35)

\rput(-12,6){$=$}

\rput(-4,30){\rnode{a1}{$C$}}
\rput(16,30){\rnode{a2}{$B \times C$}}
\rput(36,30){\rnode{a3}{$A\opb$}}
\rput(62,30){\rnode{a4}{${A'}\opb$}}

\rput(36,15){\rnode{b3}{$A\opb \ttimes B \opb$}}
\rput(62,15){\rnode{b4}{${A'}\opb \ttimes {B'}\opb$}}

\rput(11,-5){\rnode{c1}{$B \times C$}}
\rput(36,0){\rnode{c3}{$C$}}
\rput(62,0){\rnode{c4}{$C'$}}
\rput(82,0){\rnode{c5}{$B'\times C'$}}
\rput(102,0){\rnode{c6}{${{A'}\opb}.$}}

% horizontals
\ncline{->}{a1}{a2} \naput{{\scriptsize $b_1 \ttimes 1$}}

% \ncarc[arcangle=40]{->}{a1}{a2}\naput{{\scriptsize $$}}
% \ncarc[arcangle=-40]{->}{a1}{a2}\nbput{{\scriptsize $$}}

\ncline{->}{a2}{a3} \naput{{\scriptsize $G$}}
\ncline{->}{a3}{a4} \naput{{\scriptsize $S\opb$}}

\ncline{->}{b3}{b4} \naput[labelsep=1pt]{{\scriptsize $S\opb \ttimes T\opb$}}

\ncline{->}{c3}{c4} \nbput{{\scriptsize $U$}}

\ncline{->}{c4}{c5} \nbput{{\scriptsize $Tb_1 \ttimes 1$}}

\ncline{->}{c5}{c6} \nbput{{\scriptsize $G'$}}

%verticals

\psset{labelsep=1pt}{

\ncline{->}{a3}{b3} \naput{{\scriptsize $1 \ttimes b_1\opb$}}

\ncline{->}{b3}{c3} \naput{{\scriptsize $F\opb$}}

\ncline{->}{a4}{b4} \nbput{{\scriptsize $1 \ttimes Tb_1\opb$}}

\ncline{->}{b4}{c4} \nbput{{\scriptsize ${F'}\opb$}}}

%diags
\ncline{->}{a1}{c3} \nbput[npos=0.7]{{\scriptsize $1$}}
\ncarc[arcangle=20]{->}{a4}{c6} \naput{{\scriptsize $1$}}

%weirds
%\nccurve[angleA=270,angleB=120]{->}{a1}{c1} \nbput{{\scriptsize $b_2\ttimes 1$}}

\ncarc[arcangle=20]{->}{a1}{c1}\naput[labelsep=0pt,npos=0.5]{{\scriptsize $b_1 \ttimes 1$}}
\ncarc[arcangle=-30]{->}{a1}{c1}\nbput[labelsep=0pt,npos=0.4]{{\scriptsize $b_2 \ttimes 1$}}
{\psset{doubleline=true}
\rput[c](0,11){\pcline{->}(4,3)(0,0) \naput[npos=0.8,labelsep=1pt]{{\scriptsize $ f \ttimes 1$}}}}

%fake coord
\rput(83,-2){\rnode{cfake}{$$}}
\nccurve[angleA=340,angleB=220,ncurvB=0.6]{->}{c1}{cfake} \nbput{{\scriptsize $T \times U$}}

% 2-cells

\psset{doubleline=true}
\rput[c](20,19){\pcline{->}(4,4)(0,0) \naput[labelsep=1pt]{{\scriptsize 	     $\varepsilon_{b_1}$}}}
\rput[c](77,11){\pcline{->}(4,4)(0,0) \naput[labelsep=-1pt,npos=0.3]{{\scriptsize $\eta'_{Tb_1}$}}}
\rput[c](46,5){\pcline{->}(4,4)(0,0) \naput[labelsep=-1pt,npos=0.3]	       {{\scriptsize $\alpha\opb$}}}

\endpspicture
\]
\end{prf}

\begin{remark}
 Note that in the definition of the mate of $\alpha$ we could start by fixing the first variable instead of the second variable and then follow the analogous process to produce a natural transformation $\hat{\alpha}$ as below:
\[\psset{unit=0.08cm,labelsep=0pt,nodesep=3pt}
\pspicture(0,-3)(30,22)

% a1 a2 up right
% b1 b2

\rput(0,20){\rnode{a1}{$C\times A$}} % top left
\rput(30,20){\rnode{a2}{$C' \times A'$}} % top right

\rput(0,0){\rnode{b1}{$B\opb$}}   % bottom left
\rput(30,0){\rnode{b2}{${B'}\opb.$}}  % bottom right

\psset{nodesep=3pt,labelsep=2pt,arrows=->}
\ncline{a1}{a2}\naput{{\scriptsize $U \times S$}} % top
\ncline{b1}{b2}\nbput{{\scriptsize $T\opb$}} % bottom
\ncline{a1}{b1}\nbput{{\scriptsize $H$}} % left
\ncline{a2}{b2}\naput{{\scriptsize $H'$}} % right

\psset{labelsep=1.5pt}
\pnode(12,7){c1}
\pnode(18,13){c2}
\ncline[doubleline=true]{c1}{c2} \nbput[npos=0.4]{{\scriptsize $\hat{\alpha}$}}

\endpspicture\]

Note that $\bar{\hat{\alpha}} = \alpha = \hat{\bar{\alpha}}$ by the usual mates correspondence; the following result deals with a less trivial combination of these processes. This can be thought of as the 2-variable version of the mates correspondence.

\end{remark}

\begin{proposition}\label{twovarmate}

Given 2-variable adjunctions and a natural transformation $\alpha$ as above, $\bar{\bar{\alpha}} = \hat{\alpha}$. 

\end{proposition}

The proof of this result is analogous to the 1-variable case, which follows from the triangle identities for the adjunctions in question.  Therefore we start by making explicit the 2-variable version of the triangle identities, which must now involve three instances of units/counits.  In the following proofs we adopt notational shorthand as below, for simplification, clarity and to save space.

\begin{enumerate}

\item All objects have been omitted.  The source categories can always be determined from the functors shown, and whenever a variable in $A$ is required, it is understood to be $a$; likewise for $c \in C$.  For example:
\begin{itemize}
 \item $TH$ means $TH(c,a)$, and
\item $G(H,1)$ means $G(H(c,a),c))$. 
\end{itemize}

\item As in Remark~\ref{epsilonremark} we write all units and counits for all adjunctions as $\ve$; the source and target functors uniquely determine which adjunction is being used, and the object at which the component is being taken.  

\end{enumerate}

\begin{lemma}[{\bfseries Generalised triangle identity}]
For a 2-variable adjunction, the following triangles commute, along with all cyclic variants. 
% a1 a2
%  a3 

\[
\psset{unit=0.13cm,labelsep=2pt,nodesep=3pt}
\pspicture(40,20)

%%%%%%%%%% top

\rput(0,17){\rnode{a1}{$H(F,1)$}}  % top left
\rput(20,17){\rnode{a2}{$1$}}  % top right
\rput(10,5){\rnode{a3}{$H(F,G(1,F))$}}  % bottom

\ncline{->}{a1}{a2} \naput{{\scriptsize $\ve$}} % top
\ncline{->}{a1}{a3} \nbput{{\scriptsize $H(1,\ve)$}} % left
\ncline{->}{a3}{a2} \nbput{{\scriptsize $\ve$}} % right

\endpspicture
%
% a1 a2
%  a3 
%
\pspicture(20,20)

%%%%%%%%%% top

\rput(0,17){\rnode{a1}{$H(1,G)$}}  % top left
\rput(20,17){\rnode{a2}{$1$}}  % top right
\rput(10,5){\rnode{a3}{$H(F(G,1),G)$}}  % bottom

\ncline{->}{a1}{a2} \naput{{\scriptsize $\ve$}} % top
\ncline{->}{a1}{a3} \nbput{{\scriptsize $H(\ve,1)$}} % left
\ncline{->}{a3}{a2} \nbput{{\scriptsize $\ve$}} % right

\endpspicture
\]

\end{lemma}

\begin{prf}
 This follows from the ``cycle of isomorphisms'' as in Example~\ref{cycliso}.  
\end{prf}

\begin{prfof}{Proposition \ref{twovarmate}}
It suffices to show that these two natural transformations have the same component at $(c,a) \in C \times A$.  This is shown in the following (large) commutative diagram in which the top edge is the component $\bar{\bar{\alpha}}_{c,a}$ and the bottom edge $\hat{\alpha}_{c,a}$. 

Regions (3) and (4) are naturality squares, (5) and (6) are functorality of $H'$, (2) and (7) are generalised triangle identities and (1) commutes by dinaturality of $\ve$ as follows.  The counit $\ve$ in question has components
\[H'(c,G'(b,c)) \tra b\]
and is natural in $b$ but dinatural in $c$.  Writing out the dinaturality condition for the morphism
\[F'(S,TH) \tmap{\alpha} UF(1,H) \tmap{U\ve} U\]
yields region (1) as required.

% \begin{figure}[htb]
% \caption{\label{bigdiag} Diagram for proof of Proposition \ref{twovarmate}} 

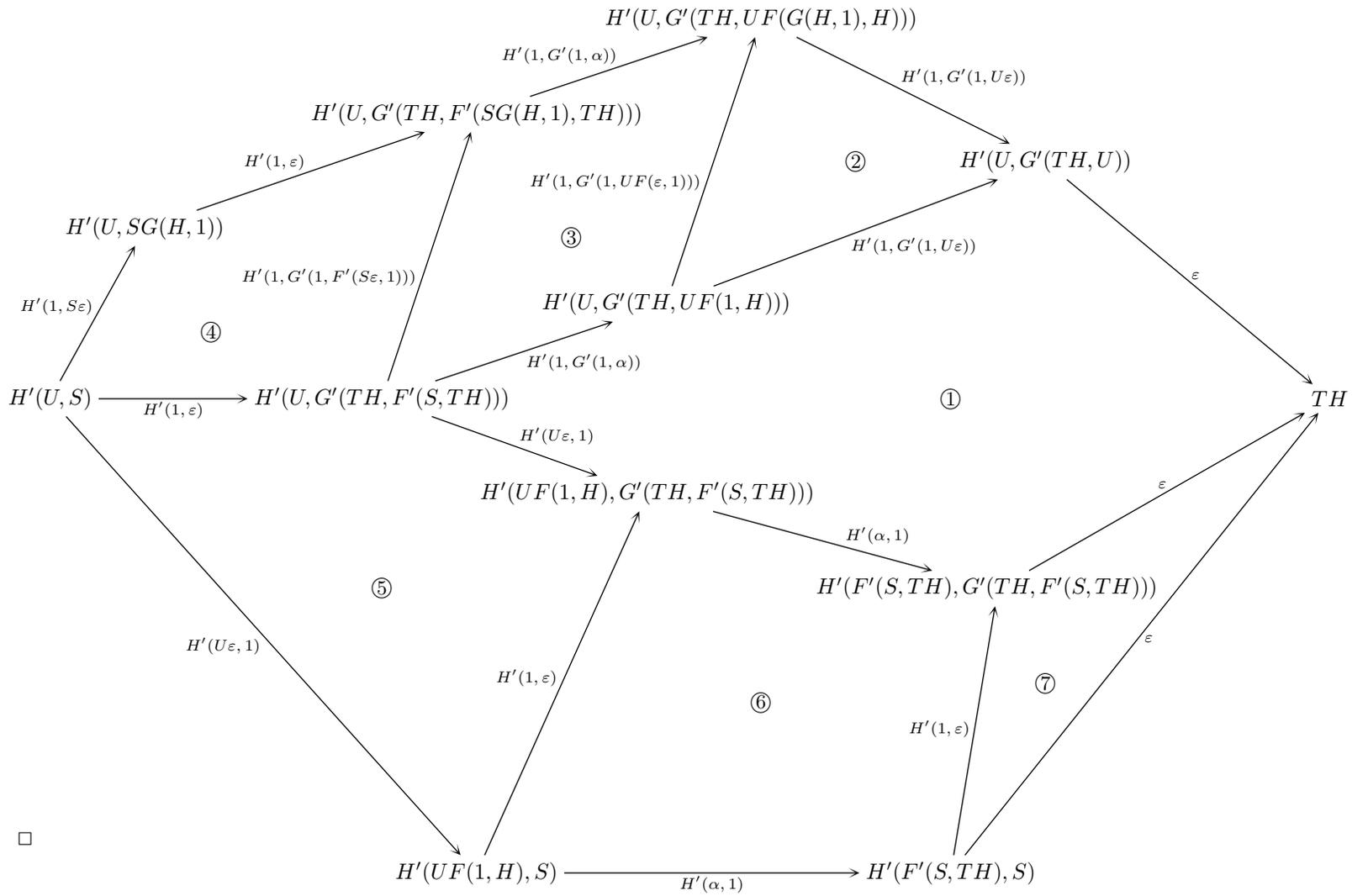
\psset{unit=1.5mm,labelsep=1pt,nodesep=3pt}
\pspicture(-10,-20)(150,100)
\rput{90}(20,50){

\pspicture(150,100)
\rput(5,50){\rnode{a1}{$H'(U,S)$}}  % named node, with something placed there
\rput(15,68){\rnode{a2}{$H'(U,SG(H,1))$}}
\rput(50,80){\rnode{a3}{$H'(U,G'(TH,F'(SG(H,1),TH)))$}}
\rput(80,90){\rnode{a4}{$H'(U,G'(TH,UF(G(H,1),H)))$}}
\rput(110,75){\rnode{a5}{$H'(U,G'(TH,U))$}}
\rput(140,50){\rnode{a6}{$TH$}}
\rput(50,0){\rnode{a7}{$H'(UF(1,H),S)$}}
\rput(100,0){\rnode{a8}{$H'(F'(S,TH),S)$}}

\rput(40,50){\rnode{b1}{$H'(U,G'(TH,F'(S,TH)))$}}
\rput(70,60){\rnode{b2}{$H'(U,G'(TH,UF(1,H)))$}}
\rput(68,40){\rnode{b3}{$H'(UF(1,H),G'(TH,F'(S,TH)))$}}
\rput(105,30){\rnode{b4}{$H'(F'(S,TH),G'(TH,F'(S,TH)))$\ \ \ }}

\ncline{->}{a1}{a2} \naput{{\scriptsize $H'(1,S\varepsilon$)}}
\ncline{->}{a2}{a3} \naput[npos=0.4]{{\scriptsize $H'(1,\ve)$}}
\ncline{->}{a3}{a4} \naput[npos=0.3]{{\scriptsize $H'(1,G'(1,\alpha))$}}
\ncline{->}{a4}{a5} \naput[npos=0.6]{{\scriptsize $H'(1,G'(1,U\ve))$}}
\ncline{->}{a5}{a6} \naput{{\scriptsize $\ve$}}

\ncline{->}{a1}{a7} \nbput{{\scriptsize $H'(U\ve,1)$}}
\ncline{->}{a7}{a8} \nbput{{\scriptsize $H'(\alpha,1)$}}
\ncline{->}{a8}{a6} \nbput{{\scriptsize $\ve$}}

\ncline{->}{b1}{b2} \nbput[npos=0.7]{{\scriptsize $H'(1,G'(1,\alpha))$}}
\ncline{->}{b2}{a5} \nbput[npos=0.6]{{\scriptsize $H'(1,G'(1,U\ve))$}}
\ncline{->}{b1}{b3} \naput[npos=0.65]{{\scriptsize $H'(U\ve,1)$}}
\ncline{->}{b3}{b4} \naput[npos=0.7]{{\scriptsize $H'(\alpha,1)$}}
\ncline{->}{b4}{a6} \naput{{\scriptsize $\ve$}}

\ncline{->}{a1}{b1} \nbput{{\scriptsize $H'(1,\ve)$}}
\ncline{->}{b1}{a3} \naput[npos=0.4]{{\scriptsize $H'(1,G'(1,F'(S\ve,1)))$}}
\ncline{->}{b2}{a4} \naput[npos=0.4]{{\scriptsize $H'(1,G'(1,UF(\ve,1)))$}}
\ncline{->}{a7}{b3} \naput{{\scriptsize $H'(1,\ve)$}}
\ncline{->}{a8}{b4} \naput{{\scriptsize $H'(1,\ve)$}}

\rput(100,50){
\rput(0,0){$\bigcirc$}
\rput(0,0){\small{{1}}}}

\rput(90,75){
\rput(0,0){$\bigcirc$}
\rput(0,0){\small{2}}}

\rput(60,67){
\rput(0,0){$\bigcirc$}
\rput(0,0){\small{3}}}

\rput(22,57){
\rput(0,0){$\bigcirc$}
\rput(0,0){\small{4}}}

\rput(40,30){
\rput(0,0){$\bigcirc$}
\rput(0,0){\small{5}}}

\rput(80,18){
\rput(0,0){$\bigcirc$}
\rput(0,0){\small{6}}}

\rput(110,20){
\rput(0,0){$\bigcirc$}
\rput(0,0){\small{7}}}

\endpspicture}
\endpspicture

%\end{figure}

\end{prfof}

\vspace{2em}

% 
% \begin{proposition}
%  This mate correspondence respects horizontal composition.
% \end{proposition}
% 
% 
% 
% \begin{prf}
%  Follows from the analogous result for ordinary (1-variable) mates.
% \end{prf}
% 
% Note that the mates correspondence also respects vertical composition in a certain sense that is more difficult to make precise; we will do this in the next section **xref**.  
% 

Now by allowing $B$ to be a finite product of categories, we get a notion of $n$-variable mates with respect to an $n$-variable adjunction.

\begin{theorem}\label{nmatethm}
 Suppose we have functors
\[\begin{array}{lcl}
 &  F_0 \: A_1 \times \cdots \times A_n \tra A_0\opb,\\
&F_0' \: A_1' \times \cdots \times A_n' \tra {A_0'}\opb &\\[4pt]
\end{array}\]
equipped with $n$-variable left adjoints, and for all $0 \leq i \leq n$ a functor
\[S_i \: A_i \tra A_i'.\]
Then a \nt\ 
\[\psset{unit=0.08cm,labelsep=0pt,nodesep=3pt}
\pspicture(70,22)

% a1 a2 up right
% b1 b2

\rput(0,20){\rnode{a1}{$A_1 \times \cdots \times A_n$}} % top left
\rput(70,20){\rnode{a2}{$A_1' \times \cdots \times A_n'$}} % top right

\rput(0,0){\rnode{b1}{$A_{0}\opb$}}   % bottom left
\rput(70,0){\rnode{b2}{${A_{0}'}\opb$}}  % bottom right

\psset{nodesep=3pt,labelsep=2pt,arrows=->}
\ncline{a1}{a2}\naput{{\scriptsize $S_1 \times \cdots \times S_n$}} % top
\ncline{b1}{b2}\nbput{{\scriptsize $S_{0}\opb$}} % bottom
\ncline{a1}{b1}\nbput{{\scriptsize $F_0$}} % left
\ncline{a2}{b2}\naput{{\scriptsize $F_0'$}} % right

\psset{labelsep=1.5pt}
\pnode(32,7){c1}
\pnode(38,13){c2}
\ncline[doubleline=true]{c1}{c2} \nbput[npos=0.4]{{\scriptsize $\alpha$}}

\endpspicture\]
has for all $1 \leq i \leq n$ a mate
\[\psset{unit=0.08cm,labelsep=0pt,nodesep=3pt}
\pspicture(70,22)

% a1 a2 up right
% b1 b2

\rput(0,20){\rnode{a1}{$A_{i+1} \times \cdots \times A_{i-1}$}} % top left
\rput(70,20){\rnode{a2}{$A_{i+1}' \times \cdots \times A_{i-1}'$}} % top right

\rput(0,0){\rnode{b1}{$A_i\opb$}}   % bottom left
\rput(70,0){\rnode{b2}{${A_i'}\opb$}}  % bottom right

\psset{nodesep=3pt,labelsep=2pt,arrows=->}
\ncline{a1}{a2}\naput{{\scriptsize $S_{i+1} \times \cdots \times S_{i-1}$}} % top
\ncline{b1}{b2}\nbput{{\scriptsize $S_{i}\opb$}} % bottom
\ncline{a1}{b1}\nbput{{\scriptsize $F_i$}} % left
\ncline{a2}{b2}\naput{{\scriptsize $F_i'$}} % right

\psset{labelsep=1.5pt}
\pnode(32,7){c1}
\pnode(38,13){c2}
\ncline[doubleline=true]{c1}{c2} \nbput[npos=0.4]{{\scriptsize $\alpha_i$}}

\endpspicture
\]
given as in Proposition~\ref{theoremC} with
\[\begin{array}{rcl}
   A &=& A_{i} \\
B &=& A_1 \times \cdots \times A_{i-1} \times A_{i+1} \times \cdots \times A_n \\
C &=& A_0.
  \end{array}\]

\end{theorem}

We now give the $n$-variable version of the mates correspondence, which follows from the 2-variable case (Proposition~\ref{twovarmate}).  First we need to fix our notation carefully.

\subsubsection*{Notation for $n$-variable mates.}

 Suppose we have functors
\[\begin{array}{lcl}
 &  F_0 \: A_1 \times \cdots \times A_n \tra A_0\opb,\\
&F_0' \: A_1' \times \cdots \times A_n' \tra {A_0'}\opb &\\[4pt]
\end{array}\]
equipped with $n$-variable left adjoints, and for all $0 \leq i \leq n$ a functor
\[S_i \: A_i \tra A_i'.\]

Then for any $0 \leq i \leq n$, given a natural transformation $\alpha$ as below 
\[\psset{unit=0.08cm,labelsep=0pt,nodesep=3pt}
\pspicture(70,22)

% a1 a2 up right
% b1 b2

\rput(0,20){\rnode{a1}{$A_{i+1} \times \cdots \times A_{i-1}$}} % top left
\rput(70,20){\rnode{a2}{$A_{i+1}' \times \cdots \times A_{i-1}'$}} % top right

\rput(0,0){\rnode{b1}{$A_i\opb$}}   % bottom left
\rput(70,0){\rnode{b2}{${A_i'}\opb$}}  % bottom right

\psset{nodesep=3pt,labelsep=2pt,arrows=->}
\ncline{a1}{a2}\naput{{\scriptsize $S_{i+1} \times \cdots \times S_{i-1}$}} % top
\ncline{b1}{b2}\nbput{{\scriptsize $S_{i}\opb$}} % bottom
\ncline{a1}{b1}\nbput{{\scriptsize $F_i$}} % left
\ncline{a2}{b2}\naput{{\scriptsize $F_i'$}} % right

\psset{labelsep=1.5pt}
\pnode(32,7){c1}
\pnode(38,13){c2}
\ncline[doubleline=true]{c1}{c2} \nbput[npos=0.4]{{\scriptsize $\alpha$}}

\endpspicture
\]
and any $j \neq i$ we denote by $\alpha_{ij}$ the mate
\[\psset{unit=0.08cm,labelsep=0pt,nodesep=3pt}
\pspicture(70,22)

% a1 a2 up right
% b1 b2

\rput(0,20){\rnode{a1}{$A_{j+1} \times \cdots \times A_{j-1}$}} % top left
\rput(70,20){\rnode{a2}{$A_{j+1}' \times \cdots \times A_{j-1}'$}} % top right

\rput(0,0){\rnode{b1}{$A_j\opb$}}   % bottom left
\rput(70,0){\rnode{b2}{${A_j'}\opb$}}  % bottom right

\psset{nodesep=3pt,labelsep=2pt,arrows=->}
\ncline{a1}{a2}\naput{{\scriptsize $S_{j+1} \times \cdots \times S_{j-1}$}} % top
\ncline{b1}{b2}\nbput{{\scriptsize $S_{j}\opb$}} % bottom
\ncline{a1}{b1}\nbput{{\scriptsize $F_j$}} % left
\ncline{a2}{b2}\naput{{\scriptsize $F_j'$}} % right

\psset{labelsep=1.5pt}
\pnode(32,7){c1}
\pnode(38,13){c2}
\ncline[doubleline=true]{c1}{c2} \nbput[npos=0.4]{{\scriptsize $\alpha_{ij}$}}

\endpspicture
\]
produced by Theorem~\ref{nmatethm}.  Note that in this notation, the mate called $\alpha_i$ in the theorem would be called $\alpha_{0i}$.

\begin{theorem}[{\bfseries The $n$-variable mates correspondence}]\label{nvarmates}
 
Given a pair of $n$-variable adjunctions, any distinct $i,j,k$ and a natural transformation $\alpha$ as above, we have
\[(\alpha_{ij})_{jk} = \alpha_{ik}.\]

\end{theorem}

\begin{prf}
Restricting to the functors $F_i, F_j, F_k$ and fixing all variables except those in $A_i, A_j, A_k$ we get a 2-variable adjunction.  The result is then simply an instance of Proposition~\ref{twovarmate} since it suffices to check it componentwise.
\end{prf}

% Some omitted explanation
%
% The component of $\alpha_{ij}$ at $a_{j+1}, \cdots a_{j-1}$ is the same as that obtained by first fixing all variables in $\alpha$ except $a_j$ and $a_k$, then taking the mate under this 2-variable adjunction, and then looking at its component at $(a_i, a_k)$.  

\begin{cor}
 Given a pair of $n$-variable adjunctions as above and a natural transformation 
\[\psset{unit=0.08cm,labelsep=0pt,nodesep=3pt}
\pspicture(70,22)

% a1 a2 up right
% b1 b2

\rput(0,20){\rnode{a1}{$A_1 \times \cdots \times A_n$}} % top left
\rput(70,20){\rnode{a2}{$A_1' \times \cdots \times A_n'$}} % top right

\rput(0,0){\rnode{b1}{$A_{0}\opb$}}   % bottom left
\rput(70,0){\rnode{b2}{${A_{0}'}\opb$}}  % bottom right

\psset{nodesep=3pt,labelsep=2pt,arrows=->}
\ncline{a1}{a2}\naput{{\scriptsize $S_1 \times \cdots \times S_n$}} % top
\ncline{b1}{b2}\nbput{{\scriptsize $S_{0}\opb$}} % bottom
\ncline{a1}{b1}\nbput{{\scriptsize $F_0$}} % left
\ncline{a2}{b2}\naput{{\scriptsize $F_0'$}} % right

\psset{labelsep=1.5pt}
\pnode(32,7){c1}
\pnode(38,13){c2}
\ncline[doubleline=true]{c1}{c2} \nbput[npos=0.4]{{\scriptsize $\alpha$}}

\endpspicture\]
we have
\[  (( \cdots   ((\alpha_{01})_{12} ) \cdots )_{n-1,n})_{n0} = \alpha.\]
That is, taking mates $n+1$ times is the identity.
\end{cor}

% Recall that in the 1-variable case, the mate relationship gave a bijection between certain sets of natural transformations.  We now give the multivariable version of this result.  The idea is that not only do we get bijections between different sets of natural transformations, but that all the bijections commute.  We make this precise in the following results, first the 2-variable case, and then the general $n$-variable case by induction.  

Note that the $n$-variable mates correspondence respects horizontal and vertical composition.  For horizontal composition this follows immediately from the analogous result for 1-variable mates.  For vertical composition a little more effort is required, but mainly just to make precise the meaning of ``respects vertical composition''.  However this is only a matter of indices.  The idea is not hard: composition of multivariable adjunctions is defined by fixing variables and composing the resulting 1-variable adjunctions, and the mates correspondence follows likewise.

To put this result in a more precise framework we will show that we have the structure of a cyclic double multicategory.

\section{Cyclic double multicategories}\label{multicat}

In this section we give the definition of ``cyclic double multicategory'', the structure into which multivariable adjunctions and mates organise themselves.  The idea is to combine the notions of double category and cyclic multicategory so that in our motivating example the cyclic action  expresses the multivariable mates correspondence. 

Recall that a double category can be defined as a category object in \Cat; similarly a double multicategory is a category object in the category \mcat\ of multicategories, and a cyclic double multicategory is a category object in the category \cmcat\ of ``cyclic multicategories''.  (Note that this could be called a ``double cyclic multicategory'' but this might sound as if there are two cyclic actions.)

We build up to the definition step by step, with some examples.

\subsection{Plain multicategories}

We begin by recalling the definition of  plain (non-symmetric) multicategories.

\begin{mydefinition}
 Let $T$ be the free monoid monad on \Set.  Write \tspan\ for the bicategory in which
\begin{itemize}
 \item 0-cells are sets,
\item 1-cells $A \cra B$ are $T$-spans
%
%   a1 
% a3 a2
%
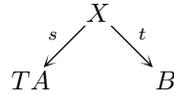
\psset{unit=0.1cm,labelsep=3pt,nodesep=1.5pt}
\pspicture(-5,9)(20,12)

%%%%%%%% top

\rput(9,9){\rnode{a1}{$X$}}  % top
\rput(0,0){\rnode{a3}{$TA$}}  % left
\rput(18,0){\rnode{a2}{$B$}}  % right

\ncline{->}{a1}{a3} \nbput{{\scriptsize $s$}} % top left
\ncline{->}{a1}{a2} \naput{{\scriptsize $t$}} %top right
\endpspicture

\vv{3em}

\item 2-cells are maps of $T$-spans.

\end{itemize}

\noi Composition is by pullback using the multiplication for $T$: the composite
\[A \cmap{X} B \cmap{Y} C\]
is given by the span
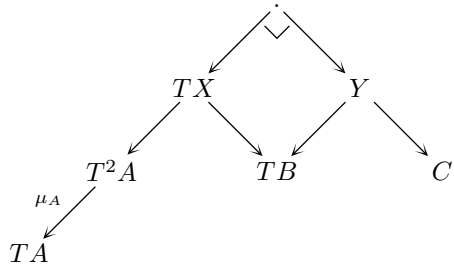
\[
\psset{unit=0.11cm,labelsep=3pt,nodesep=2pt}
\pspicture(-20,-10)(40,20)

%%%%%%%%%% top

\rput(10,20){\rnode{a1}{$.$}}  % top
\rput(0,10){\rnode{a3}{$TX$}}  % left
\rput(20,10){\rnode{a2}{$Y$}}  % right
\rput(10,0){\rnode{a4}{$TB$}}  % bottom
\rput(-10,0){\rnode{a5}{$T^2A$}} % bottom left
\rput(-20,-10){\rnode{a6}{$TA$}} % far bottom left
\rput(30,0){\rnode{a7}{$C$}} % bottom right

\ncline{->}{a1}{a3} \nbput{{\scriptsize $$}} % top left
\ncline{->}{a1}{a2} \naput{{\scriptsize $$}} %top right
\ncline{->}{a3}{a4} \nbput{{\scriptsize $$}} % bottom left
\ncline{->}{a2}{a4} \naput{{\scriptsize $$}} % bottom right

\ncline{->}{a3}{a5} \nbput{{\scriptsize $$}} % bottom left left
\ncline{->}{a5}{a6} \nbput{{\scriptsize $\mu_A$}} % bottom left left left
\ncline{->}{a2}{a7} \naput{{\scriptsize $$}}

\psline{-}(8.5,17.5)(10,16)(11.5,17.5)

\endpspicture
\]

A \demph{multicategory} $A$ is a monad in \tspan, thus
\begin{itemize}
 \item a 0-cell $A_0$,
\item a $T$-span%
%   a1 
% a3 a2
%
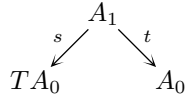
\psset{unit=0.1cm,labelsep=3pt,nodesep=1.5pt}
\pspicture(-5,9)(20,12)

%%%%%%%% top

\rput(9,9){\rnode{a1}{$A_1$}}  % top
\rput(0,0){\rnode{a3}{$TA_0$}}  % left
\rput(18,0){\rnode{a2}{$A_0$}}  % right

\ncline{->}{a1}{a3} \nbput{{\scriptsize $s$}} % top left
\ncline{->}{a1}{a2} \naput{{\scriptsize $t$}} %top right
\endpspicture

\vv{3em}
\end{itemize}
equipped with unit and multiplication 2-cells.  Explicitly, this gives
\begin{itemize}
 \item a set $A_0$ of objects,
\item for all $n \geq 0$ and objects $a_1, \ldots, a_n, a_0 \in A_0$ a set $A(a_1, \ldots, a_n; a_0)$ of $n$-ary ``multimaps''
\end{itemize}
equipped with
\begin{itemize}
 \item composition: for all sets of $k$ ordered strings $a_{i1}, \ldots, a_{im_i}, a_{i0}$ and $a_{00}$ in $A_0$ a function
\[A(a_{10}, \ldots, a_{k0}; a_{00}) \times A(\underline{a}_{1}; a_{10}) \times \cdots \times A(\underline{a}_{k}; a_{k0}) \tra A(\underline{a}_1, \ldots, \underline{a}_k; a_{00})\]
where we have written $\underline{a}_i$ for the string $a_{i1}, \ldots, a_{in_i}$, and
\item identities: for all $a \in A_0$ a function
\[1 \tmap{1_a} A(a;a)\]
\end{itemize}
satisfying the usual associativity and unit axioms.

\end{mydefinition}

Note that we can define composition at the $i$th input by composing with identities at every other input; this will be useful when giving the axioms for a cyclic multicategory and we denote it $\circ_i$.

\begin{myexamples}\label{mcatexamples} \ \\[-1em]
\begin{enumerate}
 \item Multifunctors: take objects to be categories and $k$-ary multimaps to be multifunctors, that is functors of the form
\[A_1 \times \cdots \times A_k \tra A_0.\]

\item Multicategories from monoidal categories: given any monoidal category $C$ there is a multicategory $M_C$ with the same objects, and with
\[M_C(x_1, \ldots, x_k; x_0) = C(x_1 \otimes \cdots \otimes x_k, x_0).\]

\item \label{profnonexample} Profunctors: we might try to use profunctors instead of functors in the above example, but this would form some sort of ``weak multicategory'' or ``multi-bicategory'' as profunctor composition is not strictly associative and unital.  However this is a pertinent case to consider.  A profunctor
\[A_1 \times \cdots \times A_k \pra A_0\]
is by definition a functor
\[A_1 \times \cdots \times A_k \times A_0\opb \tra \Set.\]
But this also gives rise to a profunctor
\[\hat{A}_i \times A_0\opb \pra A_i\opb\]
for each $i \neq 0$ where here $\hat{A}_i$ denotes the product 
\[A_1 \times \cdots A_{i-1} \times A_{i+1} \times \cdots \times A_k.\]
 Strictness aside, this is the sort of cyclic action we will be considering.  (In fact, the cyclic action in this example is strict although the composition is not.)

\end{enumerate}
\end{myexamples}

\subsection{Cyclic multicategories}

We now introduce the notion of a ``cyclic action'' on a multicategory.  Symmetric multicategories are multicategories with a symmetric group action that can be thought of as permuting the source elements of a given multimap.  Cyclic multicategories have a \emph{cyclic} group action that permutes the inputs \emph{and} outputs cyclically.  There is also a ``duality'' that is invoked each time an object moves between the input and output sides of a map under the cyclic action, as in the example with profunctors sketched above.  

Throughout, we work with $C_n$ the cyclic group of order $n$ considered as a subgroup of the symmetric group $S_n$ with canonical generator the cycle $\sigma_n = (123\cdots n)$.  We will often write this as $\sigma$ with its order being understood from the context.

\begin{mydefinition}\label{cyclicdef}
 A \demph{cyclic multicategory} $X$ is a multicategory equipped with
\begin{itemize}
 \item an involution on objects 
\[\begin{array}{rcl}
   X_0 & \tra & X_0\\
    x & \tmapsto & x^*
  \end{array}\]
\item for every $n \geq 0$ and ordered string $x_0, x_1, \ldots, x_n$ an isomorphism
\[\sigma = \sigma_{n+1} \: X(x_1, \cdots, x_n; x_0) \tmap{\sim} X(x_2, \ldots, x_n, x_0^*; x_1^*)\]
 \end{itemize}
such that the following axioms are satisfied.
\begin{enumerate}

\item Each isomorphism $\sigma_n$ is cyclic so that $(\sigma_n)^n = 1$.

\item The identity is preserved by $\sigma_2$, that is, the following diagram commutes
\[
\psset{unit=0.1cm,labelsep=1pt,nodesep=3pt,npos=0.4}
\pspicture(20,20)

% a1 a2
%    a3 

\rput(-3,13){\rnode{a1}{$1$}}  % top left
\rput(17,13){\rnode{a2}{$X(x;x)$}}  % top right
\rput(17,0){\rnode{a3}{$X(x^*;x^*)$}}  % bottom

\ncline[nodesepB=2pt]{->}{a1}{a2} \naput[labelsep=3pt,npos=0.6]{{\scriptsize $1_x$}} % top
\ncline{->}{a1}{a3} \nbput{{\scriptsize $1_{x^*}$}} % left
\ncline{->}{a2}{a3} \naput[labelsep=2pt]{{\scriptsize $\sigma_2$}} % right

\endpspicture
\]

\item Interaction between $\sigma$ and composition.

Let $c_i$ denote composition at the $i$th input only, that is
\[c_i \: X(y_1, \ldots, y_m; y_0) \times X(x_1, \ldots, x_n; y_i) \tra X(y_1, \ldots, y_{i-1}, x_1, \ldots, x_n, y_{i+1}, \ldots, y_m; y_0)\]

Then the following diagrams commute.

\begin{itemize}

\item For $i=1$, that is, for composition at the first input:

% a1 a2
% b1
% a3 a4

\[
\psset{unit=0.1cm,labelsep=3pt,nodesep=3pt}
\pspicture(0,-5)(75,45)

%%%%%%%%%% top

\rput(0,40){\rnode{a1}{$X(y_1, \ldots, y_m; y_0) \times X(x_1, \ldots, x_n; y_1)$}}  % top left
\rput(75,40){\rnode{a2}{$X(x_1, \ldots, x_{n}, y_2, \ldots, y_m; y_0)$}}  % top right
\rput(0,20){\rnode{b1}{$X(y_2, \ldots, y_m, y_0^*; y_1^*) \times X(x_2, \ldots, x_n, y_1^*; x_1^*) $}}  % middle left
\rput(0,0){\rnode{a3}{$X(x_2, \ldots, x_n, y_1^*; x_1^*) \times X(y_2, \ldots, y_m, y_0^*; y_1^*)$}}  % bottom left
\rput(75,0){\rnode{a4}{$X(x_2, \ldots, x_n,  y_{2}, \ldots, y_m, y_0^*; x_1^*).$}}  % bottom right

\ncline{->}{a1}{a2} \naput{{\scriptsize $c_1$}} % top
\ncline{->}{a3}{a4} \nbput{{\scriptsize $c_{-1}$}} % bottom
\ncline{->}{a1}{b1} \nbput{{\scriptsize $\sigma \times \sigma$}} % upper left
\ncline{->}{b1}{a3} \nbput{{\scriptsize $\iso$}} % lower left
\ncline{->}{a2}{a4} \naput{{\scriptsize $\sigma$}} % right

\endpspicture
\]

 \item If $i\neq 1$

% a1 a2
% a3 a4

\[
\psset{unit=0.1cm,labelsep=3pt,nodesep=3pt}
\pspicture(0,-5)(80,25)

%%%%%%%%%% top

\rput(0,20){\rnode{a1}{$X(y_1, \ldots, y_m; y_0) \times X(x_1, \ldots, x_n; y_i)$}}  % top left
\rput(80,20){\rnode{a2}{$X(y_1, \ldots, y_{i-1}, x_1, \ldots, x_n, y_{i+1}, \ldots, y_m; y_0)$}}  % top right
\rput(0,0){\rnode{a3}{$X(y_2, \ldots, y_m, y_0^*; y_1^*) \times X(x_1, \ldots, x_n; y_i)$}}  % bottom left
\rput(80,0){\rnode{a4}{$X(y_2, \ldots, y_{i-1}, x_1, \ldots, x_n, y_{i+1}, \ldots, y_m,y_0^*; y_1^*).$}}  % bottom right

\ncline{->}{a1}{a2} \naput{{\scriptsize $c_i$}} % top
\ncline{->}{a3}{a4} \nbput{{\scriptsize $c_{i-1}$}} % bottom
\ncline{->}{a1}{a3} \nbput{{\scriptsize $\sigma \times 1$}} % left
\ncline{->}{a2}{a4} \naput{{\scriptsize $\sigma$}} % right

\endpspicture
\]

\end{itemize}

In algebra
\[\sigma(g \circ_i f) = \left\{\begin{array}{ll}
                   (\sigma f) \circ_{-1} (\sigma g)  & i=1  \\
                   (\sigma g) \circ_{i-1} f  & 2 \leq i \leq n.
                 \end{array}\right.\]

\end{enumerate}
\end{mydefinition}

We can depict the axioms (3) pictorially as follows.  Depicting $f \in X(x_1, \ldots, x_n; x_0^*)$ as
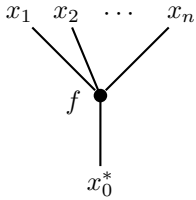
\[
\psset{linewidth=0.8pt,unit=0.9mm,labelsep=1pt,nodesep=2pt}
\pspicture(60,30)

\rput(3,22){\rnode{a1}{$x_1$}}  
\rput(10,22){\rnode{a2}{$x_2$}}  
\rput(27,22){\rnode{a3}{$x_n$}}  
\rput(15,-3){\rnode{a0}{$x_0^*$}}  
\rput(15,10){\rnode{a5}{$\pscircle*(0,0){1}
$}}  
\rput(11,9){$f$}

\rput(18,22){\rnode{a6}{$\cdots$}}

\ncline{-}{a1}{a5} \naput{{\scriptsize $$}}
\ncline{-}{a2}{a5} \naput{{\scriptsize $$}}
\ncline{-}{a3}{a5} \naput{{\scriptsize $$}}
\ncline{-}{a5}{a0} \naput{{\scriptsize $$}}

\endpspicture
\]
we depict $\sigma f$ as
\[
\psset{linewidth=0.8pt,unit=0.9mm,labelsep=1pt,nodesep=2pt}
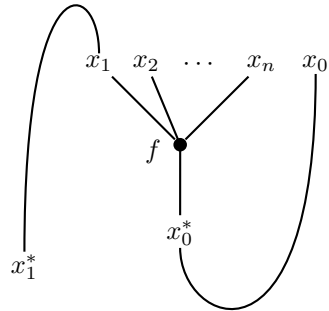
\pspicture(-10,-15)(60,30)

\rput(3,22){\rnode{x1}{$x_1$}}  
\rput(10,22){\rnode{x2}{$x_2$}}  
\rput(18,22){$\cdots$}  
\rput(27,22){\rnode{xn}{$x_n$}}  
\rput(15,10){\rnode{xc}{$\pscircle*(0,0){1}
$}}  \rput(11,9){$f$}
\rput(15,-3){\rnode{x0b}{$x_0^*$}}  

\ncline{-}{x1}{xc} \naput{{\scriptsize $$}}
\ncline{-}{x2}{xc} \naput{{\scriptsize $$}}
\ncline{-}{xn}{xc} \naput{{\scriptsize $$}}
\ncline{-}{xc}{x0b} \naput{{\scriptsize $$}}

\rput(-8,-8){\rnode{x1b}{$x_1^*$}}  
\rput(35,22){\rnode{x0}{$x_0$}}

\nccurve[nodesep=1pt,angleA=90,angleB=90,ncurvA=0.7,ncurvB=2.8]{x1}{x1b}
\nccurve[nodesep=1pt,angleA=-90,angleB=-90,ncurvA=0.7,ncurvB=2.8]{x0b}{x0}

\endpspicture
\]

Then the first axiom is depicted as shown below.

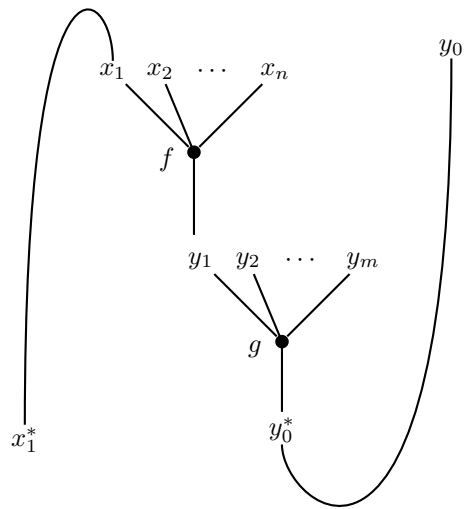
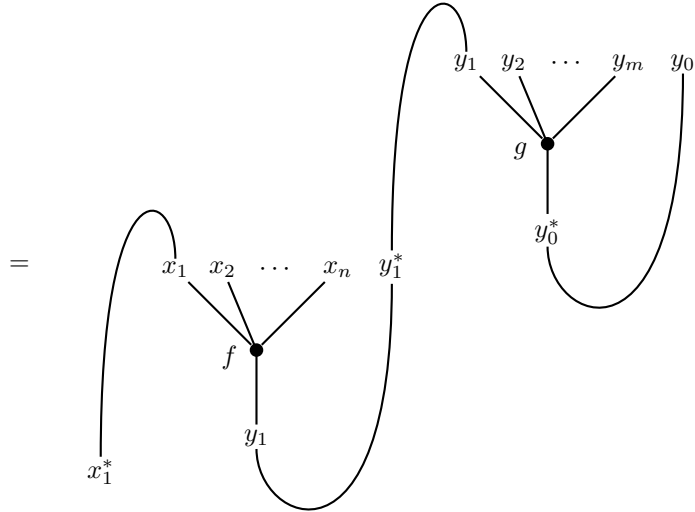
\[
\psset{linewidth=0.8pt,unit=0.9mm,labelsep=1pt,nodesep=2pt}
\pspicture(-40,-20)(50,50)

\rput(2,28)
{\pspicture(-10,-15)(60,30)

\rput(3,22){\rnode{x1}{$x_1$}}  
\rput(10,22){\rnode{x2}{$x_2$}}  
\rput(18,22){$\cdots$}  
\rput(27,22){\rnode{xn}{$x_n$}}  
\rput(15,10){\rnode{xc}{$\pscircle*(0,0){1}
$}}  \rput(11,9){$f$}
\rput(15,-3){\rnode{x0b}{$$}}  

\ncline{-}{x1}{xc} \naput{{\scriptsize $$}}
\ncline{-}{x2}{xc} \naput{{\scriptsize $$}}
\ncline{-}{xn}{xc} \naput{{\scriptsize $$}}
\ncline{-}{xc}{x0b} \naput{{\scriptsize $$}}

% \rput(-8,-8){\rnode{x1b}{$x_1\opb$}}  
% \rput(35,22){\rnode{x0}{$x_0$}}  

\endpspicture}

\rput(15,0)
{\pspicture(-10,-15)(60,30)

\rput(3,22){\rnode{y1}{$y_1$}}  
\rput(10,22){\rnode{x2}{$y_2$}}  
\rput(18,22){$\cdots$}  
\rput(27,22){\rnode{xn}{$y_m$}}  
\rput(15,10){\rnode{xc}{$\pscircle*(0,0){1}
$}}  \rput(11,9){$g$}
\rput(15,-3){\rnode{y0b}{$y_0^*$}}  

\ncline{-}{y1}{xc} \naput{{\scriptsize $$}}
\ncline{-}{x2}{xc} \naput{{\scriptsize $$}}
\ncline{-}{xn}{xc} \naput{{\scriptsize $$}}
\ncline{-}{xc}{y0b} \naput{{\scriptsize $$}}

\endpspicture}

\rput(-33,-12){\rnode{x1b}{$x_1^*$}}  
\rput(30,46){\rnode{y0}{$y_0$}}

\nccurve[nodesep=1pt,angleA=90,angleB=90,ncurvA=0.4,ncurvB=2.7]{x1}{x1b}
\nccurve[nodesep=1pt,angleA=-90,angleB=-90,ncurvA=0.3,ncurvB=2.8]{y0b}{y0}

\endpspicture
\]

\[
\psset{linewidth=0.8pt,unit=0.9mm,labelsep=1pt,nodesep=2pt}
\pspicture(-10,-25)(60,55)

\rput(-25,15){$=$}

\rput(20,0)
{\pspicture(-10,-15)(60,30)
\rput(3,22){\rnode{x1}{$x_1$}}  
\rput(10,22){\rnode{x2}{$x_2$}}  
\rput(18,22){$\cdots$}  
\rput(27,22){\rnode{xn}{$x_n$}}  
\rput(15,10){\rnode{xc}{$\pscircle*(0,0){1}
$}}  \rput(11,9){$f$}
\rput(15,-3){\rnode{x0b}{$y_1$}}  

\ncline{-}{x1}{xc} \naput{{\scriptsize $$}}
\ncline{-}{x2}{xc} \naput{{\scriptsize $$}}
\ncline{-}{xn}{xc} \naput{{\scriptsize $$}}
\ncline{-}{xc}{x0b} \naput{{\scriptsize $$}}

\rput(-8,-8){\rnode{x1b}{$x_1^*$}}  
\rput(35,22){\rnode{x0}{\textcolor{white}{$y_1^*$}}}

\nccurve[nodesep=1pt,angleA=90,angleB=90,ncurvA=0.7,ncurvB=2.8]{x1}{x1b}
\nccurve[nodesep=1pt,angleA=-90,angleB=-90,ncurvA=0.7,ncurvB=2.8]{x0b}{x0}
\endpspicture}

\rput(63,30.5)
{\pspicture(-10,-15)(60,30)
\rput(3,22){\rnode{x1}{$y_1$}}  
\rput(10,22){\rnode{x2}{$y_2$}}  
\rput(18,22){$\cdots$}  
\rput(27,22){\rnode{xn}{$y_m$}}  
\rput(15,10){\rnode{xc}{$\pscircle*(0,0){1}
$}}  \rput(11,9){$g$}
\rput(15,-3){\rnode{x0b}{$y_0^*$}}  

\ncline{-}{x1}{xc} \naput{{\scriptsize $$}}
\ncline{-}{x2}{xc} \naput{{\scriptsize $$}}
\ncline{-}{xn}{xc} \naput{{\scriptsize $$}}
\ncline{-}{xc}{x0b} \naput{{\scriptsize $$}}

\rput(-8,-8){\rnode{x1b}{$y_1^*$}}  
\rput(35,22){\rnode{x0}{$y_0$}}

\nccurve[nodesep=1pt,angleA=90,angleB=90,ncurvA=0.7,ncurvB=2.8]{x1}{x1b}
\nccurve[nodesep=1pt,angleA=-90,angleB=-90,ncurvA=0.7,ncurvB=2.8]{x0b}{x0}
\endpspicture}

\endpspicture
\]
The second is a little harder to depict as the diagrams are ambiguous, but in fact it says that the two different ways of interpreting the diagram below are the same.

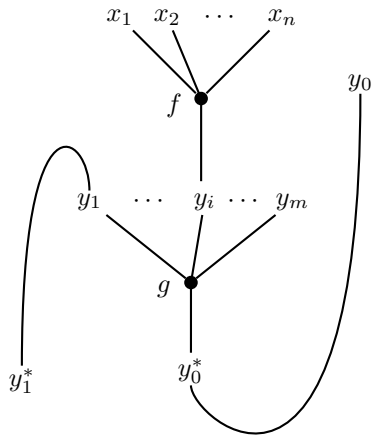
\[
\psset{linewidth=0.8pt,unit=0.9mm,labelsep=1pt,nodesep=2pt}
\pspicture(-40,-20)(50,50)

\rput(16.5,27.2)
{\pspicture(-10,-15)(60,30)

\rput(3,22){\rnode{x1}{$x_1$}}  
\rput(10,22){\rnode{x2}{$x_2$}}  
\rput(18,22){$\cdots$}  
\rput(27,22){\rnode{xn}{$x_n$}}  
\rput(15,10){\rnode{xc}{$\pscircle*(0,0){1}
$}}  \rput(11,9){$f$}
\rput(15,-3){\rnode{x0b}{$$}}  

\ncline{-}{x1}{xc} \naput{{\scriptsize $$}}
\ncline{-}{x2}{xc} \naput{{\scriptsize $$}}
\ncline{-}{xn}{xc} \naput{{\scriptsize $$}}
\ncline{-}{xc}{x0b} \naput{{\scriptsize $$}}

% \rput(-8,-8){\rnode{x1b}{$x_1\opb$}}  
% \rput(35,22){\rnode{x0}{$x_0$}}  

\endpspicture}

\rput(15,0)
{\pspicture(-10,-15)(60,30)

\rput(0,22){\rnode{y1}{$y_1$}}  
\rput(9,22){$\cdots$}  
\rput(17,22){\rnode{x2}{$y_i$}}  
\rput(23,22){$\cdots$}  
\rput(30,22){\rnode{xn}{$y_m$}}  
\rput(15,10){\rnode{xc}{$\pscircle*(0,0){1}
$}}  \rput(11,9){$g$}
\rput(15,-3){\rnode{y0b}{$y_0^*$}}  

\ncline{-}{y1}{xc} \naput{{\scriptsize $$}}
\ncline{-}{x2}{xc} \naput{{\scriptsize $$}}
\ncline{-}{xn}{xc} \naput{{\scriptsize $$}}
\ncline{-}{xc}{y0b} \naput{{\scriptsize $$}}

\endpspicture}

\rput(-20,-12){\rnode{x1b}{$y_1^*$}}  
\rput(30,32){\rnode{y0}{$y_0$}}

\nccurve[nodesep=1pt,angleA=90,angleB=90,ncurvA=0.6,ncurvB=3]{y1}{x1b}
\nccurve[nodesep=1pt,angleA=-90,angleB=-90,ncurvA=0.2,ncurvB=2.8]{y0b}{y0}

\endpspicture
\]

Note that these two axioms are equivalent to a single axiom involving the cyclic action and composition at every variable.

\begin{myexamples} First we give some slightly degenerate examples.
 
\begin{enumerate}
 \item If $X$ is a cyclic multicategory with only one object (so the involution must be the identity) then we have a notion of non-symmetric cyclic operad as given by Batanin and Berger in \cite{bb1}.  This is in contrast to the definition of (symmetric) cyclic operad in \cite{kv1}.

\item More generally, the involution can be the identity even for a non-trivial set of objects.  
\end{enumerate}

\end{myexamples}

\begin{myexample} 
We define a cyclic multicategory \MAdj\ as follows.  Take objects to be categories, and a multimap
\[A_1 , \ldots,  A_n \tra A_0\]
to be a functor
\[A_1 \times \cdots \times A_n \tmap{F_0} A_0\]
equipped with all $n$-variable left adjoints, $F_1, \ldots, F_n$.  The involution $()^\ast$ is then given by $()\opb$ and the cyclic action is  given by
\[\sigma \: F_i \tmapsto F_{i+1}\]
and the axioms are satisfied by construction.   We could also do this with $n$-variable right adjoints.  

\end{myexample}

We now characterise \MAdj\ using profunctors.  Recall our profunctor example that was not quite a true example (Example \ref{mcatexamples}.\ref{profnonexample}) as profunctor composition is not strictly associative or unital; nevertheless it has a strict cyclic action.

In fact the multivariable profunctors are $n$-variable adjunctions \emph{internal} to the bicategory \Prof\ of categories, profunctors and \nts; all profunctors have such adjoints, just as all (1-variable) profunctors have left and right adjoints.  Since functors are representable profunctors, $n$-variable adjunctions can be thought of as $n$-variable profunctors $F_0$  such that $F_0$ and all its $n$-variable adjoints in \Prof\ are representable, or, more precisely, equipped with representations as follows.

\begin{prop}\label{threepointsix}
Let $P : A_1 \times \cdots \times A_n \times A_0 \tra \Set$ be a profunctor equipped with a representation for each $P(a_1, \ldots, a_{i-1}, \uscore, a_{i+1}, \ldots, a_n, a_0)$.  That is, given
\[a_{i+1}, \ldots, a_{i-1}\]
an object $F_i(a_{i+1}, \ldots, a_{i-1}) \in A_i$ and an isomorphism
\begin{equation}\label{equ1}
P(a_1, \ldots, a_0) \cong A_i(F_i(a_{i+1}, \ldots, a_{i-1}), a_i)
\end{equation}
natural in $a_i$.  

Then the $F_i$ canonically extend to functors
\[A_{i+1} \times \cdots \times A_{i-1} \tmap{F_i} A_i\opb\]
forming an $n$-variable (left) adjunction.
\end{prop}

\begin{prf}
By standard results about parametrised representability, each $F_i$ extends to a functor
\[A_{i+1} \times \cdots \times A_{i-1} \tmap{F_i} A_i\opb\]
unique making the isomorphism (\ref{equ1}) above natural in all variables.  For the $n$-variable adjunction we then compose the isomorphisms
\[A_j(F_j(a_{j+1}, \ldots, a_{j-1}), a_j) \tmap{\cong} P(a_1, \ldots, a_0) \tmap{\cong} A_i(F_i(a_{i+1}, \ldots, a_{i-1}), a_i).  \]
 \end{prf}

\begin{prop}
 Such profunctors \emph{equipped} with such representations form a cyclic multicategory isomorphic to \MAdj.
\end{prop}

\begin{myremark}
 Note that the content of this result is that, composition of profunctors equipped with these choices of representation can be made strict.  In general profunctor composition is by coend, and the isomorphic possible choices of these coends results in associativity and unit axioms only holding up to isomorphism.   However, the representations effectively give us a coherent choice of coends for forming the composites in such a way that the axioms are guaranteed to hold strictly.  The only difficulty is notation.
\end{myremark}

\begin{proof}
 We must define composition.  We use notation analogous to Theorem~\ref{twopointseven} and start with the 1-variable case.

Given profunctors $A \pmap{P} B \pmap{Q} C$ so functors
\[\begin{array}{cl}
   A \times B\opb \tmap{P} \Set, & \mbox{and}\\
   B \times C\opb \tmap{Q} \Set
  \end{array}\]
with adjunctions
\[
\psset{labelsep=2pt}
\pspicture(40,8)
\rput(0,3){\rnode{A}{$A$}}
\rput(18,3){\rnode{B}{$B$}}
\rput(36,3){\rnode{C}{$C$}}

\psset{nodesep=2pt,arrows=->, labelsep=2pt}
\ncline[offset=5pt]{A}{B}\nbput{$\scriptstyle\top$}\naput{{\scriptsize $G_1$}}
\ncline[offset=5pt]{B}{A}\naput{{\scriptsize $F_1$}}

\ncline[offset=5pt]{B}{C}\nbput{$\scriptstyle\top$}\naput{{\scriptsize $G_2$}}
\ncline[offset=5pt]{C}{B}\naput{{\scriptsize $F_2$}}

\endpspicture
\]
with
\[\begin{array}{c}
P(a,b)  \cong  A(F_1b, a) \cong B(b,G_1a)  \\
Q(b,c) \cong  B(F_2c,b) \cong C(c, G_2b) 
\end{array}
\]
natural in $a,b,c$; here the first isomorphism in each row is from the given representations, and the second is from Proposition~\ref{threepointsix}.  We then define the composite
\[\begin{array}{rcll}
Q \circ P & = & C(c, G_2G_1a) \\
& \cong & A(F_1F_2c, a) & \mbox{by composition of adjunctions.}
\end{array}\]
As this composition has been defined by composing the functors $G_1, G_2$ and the functors $F_1, F_2$, it is strictly associative and unital.  The $n$-variable case follows similarly, as does the cyclic action and the isomorphism of multicategories.
\end{proof}

\begin{myremark} 
 The idea is that we consider the functor
\[\Cat \tra \Prof\]
that is the identity on objects and on morphisms sends a functor
\[A \tmap{G} B\]
to the profunctor
\[A \pra B\]
given by
\[\begin{array}{ccl}
   A \times B\opb & \tra & \Set \\
   (a,b) & \tmapsto & B(b,Ga).
  \end{array}\]
With the usual composition in \Prof\ this is only a pseudo-functor, giving us a ``sub-pseudo-multicategory'' of \Prof\ that is somehow ``equivalent'' to \MAdj. In order to get an honest multicategory we must specifiy data as above, giving us a strict multicategory biequivalent to the more natural arising pseudo-multicategory.
\end{myremark}

% Before introducing 2-cells to our multicategories, we include one last abstract result although, again, we will not be using it any further in this work.

% \begin{theorem}
%  The forgetful functor
% \[\CMcat \tra \Mcat\]
% is monadic.
% \end{theorem}
% 
% \begin{proof}
%  The free functor adds cyclic actions freely.  ?? think about involution a bit. It adds the involution first, and then cyclic copies of morphisms after.
% \end{proof}

\subsection{Cyclic double multicategories}

We are now ready to introduce the ``cubical'' 2-cells we need.  Recall that a double category can be defined very succinctly as a category object in the category of (small) categories.  We proceed analogously for the multi-versions.

\begin{mydefinition}
 A \demph{double multicategory} is a category object in the category \Mcat\ of multicategories.

A \demph{cyclic double multicategory} is a category object in the category \CMcat\ of cyclic multicategories.
\end{mydefinition}

Note that pullbacks in the category \CMcat\ are defined in the obvious way, so this definition makes sense.  As with double categories, it is desirable to give an elementary description.  A cyclic double multicategory $X$ has as underlying data a diagram
\[\psset{labelsep=2pt,nodesep=2pt}
\pspicture(15,5)
\rput(0,3){\rnode{a1}{B}}
\rput(15,3){\rnode{a2}{A}}
\psset{nodesep=3pt,arrows=->}
\ncline[offset=3pt]{a1}{a2}\naput{{\scriptsize $s$}}
\ncline[offset=-3pt]{a1}{a2}\nbput{{\scriptsize $t$}}
\endpspicture
\]
in \CMcat.  

% Thus we have the following underlying data.
% 
% \begin{itemize}
%  \item A cyclic multicategory $A$ with objects $x$ and multi-homsets
% \[A(x_1, \ldots, x_n; x_0).\]
% We call its objects we the 0-cells of $X$ and its multimaps the vertical 1-cells of $X$.
% 
% \item A cyclic multicategory $B$ with objects $f$ and multi-homsets
% \[B(f_1, \ldots, f_n; f_0).\]
% We call its objects the horizontal 1-cells of $X$ and its multimaps the 2-cells of $X$.
% \end{itemize}

Recall that the underlying data for a multicategory $A$ is in turn a diagram in sets of the following form
%
%   a1 
% a3 a2
%
\[
\psset{unit=0.1cm,labelsep=3pt,nodesep=1.5pt}
\pspicture(20,12)

%%%%%%%% top

\rput(9,9){\rnode{a1}{$A_1$}}  % top
\rput(0,0){\rnode{a3}{$TA_0$}}  % left
\rput(18,0){\rnode{a2}{$A_0$}}  % right

\ncline{->}{a1}{a3} \nbput{{\scriptsize $d$}} % top left
\ncline{->}{a1}{a2} \naput{{\scriptsize $c$}} %top right
\endpspicture\]
where $T$ is the free monoid monad on \Set.  Thus for a category object in \Mcat\ we have a diagram of the following form in \Set:
\[
\psset{unit=0.1cm,labelsep=3pt,nodesep=1.5pt}
\pspicture(20,25)

%   b1
% b3  b2

%
%   a1 
% a3 a2

%%%%%%%% top

\rput(10,25){\rnode{b1}{$B_1$}}  % top
\rput(0,15){\rnode{b3}{$TB_0$}}  % left
\rput(20,15){\rnode{b2}{$B_0$}}  % right

\rput(10,10){\rnode{a1}{$A_1$}}  % top
\rput(0,0){\rnode{a3}{$TA_0$}}  % left
\rput(20,0){\rnode{a2}{$A_0$}}  % right

\psset{nodesep=2pt,arrows=->, labelsep=1pt}

\ncline{->}{a1}{a3} \naput[labelsep=0pt,npos=0.45]{{\scriptsize $d$}} % top left
\ncline{->}{a1}{a2} \nbput{{\scriptsize $c$}} %top right

\ncline{->}{b1}{b3} \nbput{{\scriptsize $d$}} % top left
\ncline{->}{b1}{b2} \naput{{\scriptsize $c$}} %top right

\ncline[offset=2pt]{b1}{a1}\naput{{\scriptsize $t$}}
\ncline[offset=-2pt]{b1}{a1}\nbput{{\scriptsize $s$}}

\ncline[offset=2pt]{b2}{a2}\naput{{\scriptsize $t$}}
\ncline[offset=-2pt]{b2}{a2}\nbput{{\scriptsize $s$}}

\ncline[offset=2pt]{b3}{a3}\naput{{\scriptsize $t$}}
\ncline[offset=-2pt]{b3}{a3}\nbput{{\scriptsize $s$}}

\endpspicture\]
where the sets correspond to data as follows:
\[\begin{array}{rcl}
   A_0 &=& \mbox{0-cells} \\
   A_1 &=& \mbox{vertical (multi) 1-cells} \\
   B_0 &=& \mbox{horizontal (plain) 1-cells} \\
   B_1 &=& \mbox{2-cells}.
  \end{array}\]
Commuting conditions tell us that 2-cells might be depicted as: 
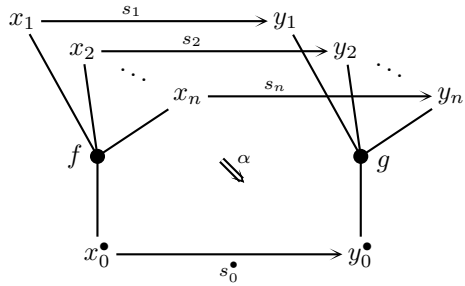
\[
\psset{linewidth=0.8pt,unit=0.1cm,labelsep=1pt,nodesep=2pt}
\pspicture(0,-10)(60,30)
\rput(0,28){\rnode{a1}{$x_1$}}  
\rput(8,24){\rnode{a2}{$x_2$}}  
\rput(22,18){\rnode{a3}{$x_n$}}  
\rput(10,-3){\rnode{a0}{$x_0\opb$}}  
\rput(10,10){\rnode{a5}{$\pscircle*(0,0){1}
$}}  
\rput(7,10){$f$}

\rput{-26}(15,20.7){\rnode{a6}{$\cdots$}}  
\rput{-26}(49,21.3){\rnode{a6}{$\cdots$}}

\rput(35,28){\rnode{b1}{$y_1$}}  
\rput(43,24){\rnode{b2}{$y_2$}}  
\rput(57,18){\rnode{b3}{$y_n$}}  
\rput(45,-3){\rnode{b0}{$y_0\opb$}}  
\rput(45,10){\rnode{b5}{$\pscircle*(0,0 ){1}$}}  
\rput(48,9){$g$}

\ncline{-}{a1}{a5} \naput{{\scriptsize $$}}
\ncline{-}{a2}{a5} \naput{{\scriptsize $$}}
\ncline{-}{a3}{a5} \naput{{\scriptsize $$}}
\ncline{-}{a5}{a0} \naput{{\scriptsize $$}}

\ncline{-}{b1}{b5} \naput{{\scriptsize $$}}
\ncline{-}{b2}{b5} \naput{{\scriptsize $$}}
\ncline{-}{b3}{b5} \naput{{\scriptsize $$}}
\ncline{-}{b5}{b0} \naput{{\scriptsize $$}}

\psset{npos=0.4}

\ncline{->}{a1}{b1} \naput{{\scriptsize $s_1$}}
\ncline{->}{a2}{b2} \naput{{\scriptsize $s_2$}}
\ncline{->}{a3}{b3} \naput[npos=0.3]{{\scriptsize $s_n$}}
\ncline{->}{a0}{b0} \nbput[labelsep=2pt,npos=0.5]{{\scriptsize $s_0\opb$}}

{\psset{doubleline=true}
\rput[c](26,6){\pcline{->}(0,4)(4,0) \naput[labelsep=1pt]{{\scriptsize $\alpha$}}}}

\endpspicture
\]

Inside the structure of a (cyclic) double multicategory we have two categories given by 
\begin{itemize}
 \item 0-cells, horizontal 1-cells and horizontal composition, and
\item vertical 1-cells, 2-cells and horizontal composition
\end{itemize}
and two (cyclic) multicategories with objects and multimaps given by
\begin{itemize}
 \item 0-cells, vertical 1-cells and vertical multi-composition, and
\item horizontal 1-cells, 2-cells and vertical multi-composition.
\end{itemize}
Furthermore these must all be compatible, in the following sense.  In addition to the underlying diagram 
\[\psset{labelsep=2pt,nodesep=2pt}
\pspicture(15,5)
\rput(0,3){\rnode{a1}{B}}
\rput(15,3){\rnode{a2}{A}}
\psset{nodesep=3pt,arrows=->}
\ncline[offset=3pt]{a1}{a2}\naput{{\scriptsize $s$}}
\ncline[offset=-3pt]{a1}{a2}\nbput{{\scriptsize $t$}}
\endpspicture
\]
in \CMcat\ we must have an identity map 
\[I \: A \tra B\] 
and a composition map 
\[\gamma \: B \times_A B \tra B\]
and  $s,t,I,\gamma$ must all be maps of cyclic multicategories, that is, they must respect (co)domains, composition, involution and cyclic actions in passing from $B$ to $A$.  Note that $s/t$ give ``horizontal source and target'', $I$ gives ``horizontal identities'' and $\gamma$ ``horizontal composition''.  Respecting (co)domains and composition is analogous to the axioms for a double category, just with multimaps instead of 1-ary maps where appropriate; notably this gives us interchange between horizontal and vertical composition.  

Respecting involution and cyclic actions gives us the following information.

\begin{enumerate}
 \item Horizontal source and target respect involution: we have an involution on $A_0$ (0-cells) and an involution on $B_0$ (horizontal 1-cells), both written $(\ )^*$ under which
\[ x \tmap{f} y \hh{2em} \tmapsto \hh{2em} x^* \tmap{f^*} y^*.\]

\item Horizontal identities respect involution: for any 0-cell $x \in A_0$ we have a horizontal 1-cell identity $I_x \: x \tra x \in B_0$.  This assignation must satisfy the following equality of horizontal 1-cells:
\[I_{x^*} = (I_x)^*.\]

\item Horizontal composition respects involution: given composable horizontal 1-cells in $B_0$
\[x \map{f} y \map{g} z \]
we must have the following equality of horizontal 1-cells:
\[(gf)^* = g^* f^*.\]

\item Horizontal source and target respect cyclic action: given a 2-cell $\alpha \in B_1$ we have the following equalities of vertical 1-cells
\[\begin{array}{cl}
   s(\sigma \alpha) = \sigma(s\alpha), & \mbox{and} \\
   t(\sigma \alpha) = \sigma(t\alpha).
  \end{array}\] 

\item Horizontal identities respect cyclic action: given a vertical 1-cell $f \in A_1$ we have a horizontal 2-cell identity $I_f \in B_1$.  This assignation must satisfy the following equality of 2-cells:
\[\sigma(I_f) = I_{\sigma f}.\]

\item Horizontal composition respects cyclic action: given horizontally composable 2-cells $\alpha, \beta \in B_1$ we have the following equality of 2-cells:
\[\sigma(\beta \ast \alpha) = \sigma\beta \ast \sigma\alpha\]
where as usual we write horizontal composition of 2-cells as $\beta \ast \alpha$.

\end{enumerate}

\subsection{Multivariable adjunctions}

In this section we show how to organise multivariable adjunctions and mates into a cyclic double multicategory.  In fact, just as for the 1-variable case, there are many choices of such a structure on this underlying data.  The difference is that now, because of the extra variables, there are also extra choices but many of them are rather unnatural so there is more danger of confusion.  The other source of confusion is that the standard notation used in the 1-variable case does not generalise very easily to express all the possible choices in the multivariable case.  We begin by giving the most obvious choices of structure.

\begin{theorem}\label{matesvert}\label{mainthm}

There is a cyclic double multicategory \MADJ\ extending the cyclic multicategory \MAdj\ of multivariable (left) adjunctions, given as follows.
\begin{itemize}
 \item 0-cells are categories.
\item Horizontal 1-cells are functors.
\item A vertical 1-cell $A_1 , \ldots,  A_n \map{F} A_0\opb$ is a functor $F$ equipped with $n$-variable left adjoints.

\item 2-cells are natural transformations
\[\psset{unit=0.08cm,labelsep=0pt,nodesep=3pt}
\pspicture(0,-2)(60,24)

% a1 a2 up right
% b1 b2

\rput(0,20){\rnode{a1}{$A_1 \times \cdots \times A_n $}} % top left
\rput(60,20){\rnode{a2}{$B_1 \times \cdots \times B_n$}} % top right

\rput(0,0){\rnode{b1}{$A_0\opb$}}   % bottom left
\rput(60,0){\rnode{b2}{${B_0}\opb$}}  % bottom right

\psset{nodesep=3pt,labelsep=2pt,arrows=->}
\ncline{a1}{a2}\naput{{\scriptsize $S_1 \times \cdots \times S_n$}} % top
\ncline{b1}{b2}\nbput{{\scriptsize $S_0\opb$}} % bottom
\ncline{a1}{b1}\nbput{{\scriptsize $F$}} % left
\ncline{a2}{b2}\naput{{\scriptsize $G$}} % right

\psset{labelsep=1.5pt}
\pnode(30,7){c1}
\pnode(36,13){c2}
\ncline[doubleline=true]{c1}{c2} \nbput[npos=0.4]{{\scriptsize $\alpha$}}

\endpspicture\]
(note direction).  Here, despite the direction of the \nt, the horizontal source of $\alpha$  \emph{as a 2-cell of \MADJ} is $F$ and the horizontal target is $G$; the vertical source is $S_1, \ldots, S_n$ and the vertical target is $S_0\opb$.

\end{itemize}
 
The cyclic action on 2-cells is given by the multivariable mates correspondence.

\end{theorem}

\begin{prf}
It only remains to prove that the cyclic composition axioms hold for 2-cells; these are the axioms given in Definition~\ref{cyclicdef}, applied to the multicategory whose objects are horizontal 1-cells and whose multimaps are 2-cells.  For the first axiom, it suffices to consider the following 2-cells.

\[\psset{unit=0.08cm,labelsep=0pt,nodesep=3pt}
\pspicture(20,55)

\rput(0,45){\pspicture(20,22)

% a1 a2 up right
% b1 b2

\rput(0,20){\rnode{a1}{$A\times B$}} % top left
\rput(30,20){\rnode{a2}{$A' \times B'$}} % top right

\rput(0,0){\rnode{b1}{$C\opb$}}   % bottom left
\rput(30,0){\rnode{b2}{${C'}\opb$}}  % bottom right

\psset{nodesep=3pt,labelsep=2pt,arrows=->}
\ncline{a1}{a2}\naput{{\scriptsize $S\times T$}} % top
\ncline{b1}{b2}\nbput{{\scriptsize $U\opb$}} % bottom
\ncline{a1}{b1}\nbput{{\scriptsize $F_0$}} % left
\ncline{a2}{b2}\naput{{\scriptsize ${F_0}'$}} % right

\psset{labelsep=1pt}
\pnode(18,13){c2}
\pnode(12,7){c1}
\ncline[doubleline=true,arrowinset=0.6,arrowlength=0.8,arrowsize=0.5pt 2.1]{c1}{c2} \naput[npos=0.4]{{\scriptsize $\alpha$}}

\endpspicture}

\rput(0,15){
\pspicture(20,24)

% a1 a2 up right
% b1 b2

\rput(0,20){\rnode{a1}{$C\opb\times D$}} % top left
\rput(30,20){\rnode{a2}{${C'}\opb \times C'$}} % top right

\rput(0,0){\rnode{b1}{$E\opb$}}   % bottom left
\rput(30,0){\rnode{b2}{${E'}\opb$}}  % bottom right

\psset{nodesep=3pt,labelsep=2pt,arrows=->}
\ncline{a1}{a2}\naput{{\scriptsize $U\opb\times V$}} % top
\ncline{b1}{b2}\nbput{{\scriptsize $W\opb$}} % bottom
\ncline{a1}{b1}\nbput{{\scriptsize $G_0$}} % left
\ncline{a2}{b2}\naput{{\scriptsize ${G_0}'$}} % right

\psset{labelsep=1pt}
\pnode(18,13){c2}
\pnode(12,7){c1}
\ncline[doubleline=true,arrowinset=0.6,arrowlength=0.8,arrowsize=0.5pt 2.1]{c1}{c2} \naput[npos=0.4]{{\scriptsize $\beta$}}
\endpspicture}

\endpspicture\]

This gives the general axiom by considering $B$ and $D$ to be products.  Using multicategorical notation, and our previous notation for multivariable mates, we need to show
\[(\beta \circ_1 \alpha)_{01} = \alpha_{01} \circ_2 \beta_{01}.\]
Now the component of $(\beta \circ_1 \alpha)_{01}$ at $(b,d,e)$ is obtained as follows:
\begin{enumerate}
 \item fix $b$ and $d$ in the composite $\beta \circ_1 \alpha$,
\item take the 1-variable mate,
\item evaluate at $e$.
\end{enumerate}
Now step (1) is the same as fixing $b$ in $\alpha$, $d$ in $\beta$ and then composing the squares vertically. So the axiom is an instance of 1-variable mates respecting vertical composition.  

For the second axiom it suffices to consider the following 2-cells.

\[\psset{unit=0.08cm,labelsep=0pt,nodesep=3pt}
\pspicture(20,55)

\rput(0,45){
\pspicture(20,22)

% a1 a2 up right
% b1 b2

\rput(0,20){\rnode{a1}{$A$}} % top left
\rput(30,20){\rnode{a2}{$A'$}} % top right

\rput(0,0){\rnode{b1}{$C\opb$}}   % bottom left
\rput(30,0){\rnode{b2}{${C'}\opb$}}  % bottom right

\psset{nodesep=3pt,labelsep=2pt,arrows=->}
\ncline{a1}{a2}\naput{{\scriptsize $S$}} % top
\ncline{b1}{b2}\nbput{{\scriptsize $U\opb$}} % bottom
\ncline{a1}{b1}\nbput{{\scriptsize $F$}} % left
\ncline{a2}{b2}\naput{{\scriptsize ${F}'$}} % right

\psset{labelsep=1pt}
\pnode(18,13){c2}
\pnode(12,7){c1}
\ncline[doubleline=true,arrowinset=0.6,arrowlength=0.8,arrowsize=0.5pt 2.1]{c1}{c2} \naput[npos=0.4]{{\scriptsize $\alpha$}}

\endpspicture}

\rput(0,15){
\pspicture(20,24)

% a1 a2 up right
% b1 b2

\rput(0,20){\rnode{a1}{$B\times C\opb$}} % top left
\rput(30,20){\rnode{a2}{${B'} \times {C'}\opb$}} % top right

\rput(0,0){\rnode{b1}{$D\opb$}}   % bottom left
\rput(30,0){\rnode{b2}{${D'}\opb$}}  % bottom right

\psset{nodesep=3pt,labelsep=2pt,arrows=->}
\ncline{a1}{a2}\naput{{\scriptsize $T\times U\opb$}} % top
\ncline{b1}{b2}\nbput{{\scriptsize $V\opb$}} % bottom
\ncline{a1}{b1}\nbput{{\scriptsize $G_0$}} % left
\ncline{a2}{b2}\naput{{\scriptsize ${G_0}'$}} % right

\psset{labelsep=1pt}
\pnode(18,13){c2}
\pnode(12,7){c1}
\ncline[doubleline=true,arrowinset=0.6,arrowlength=0.8,arrowsize=0.5pt 2.1]{c1}{c2} \naput[npos=0.4]{{\scriptsize $\beta$}}
\endpspicture}

\endpspicture\]

This gives the general axiom by letting $A$ and $C$ be products.  

We need to show
\[(\beta \circ_2 \alpha)_{01} = \beta_{01} \circ_1 \alpha.\]
Note that the $\alpha$ on the right hand side is not a mate, as in the axiom given in Definition~\ref{cyclicdef}. 

The component of $(\beta \circ_2 \alpha)_{01}$ at $(a,d)$ is obtained as follows:
\begin{enumerate}
 \item fix $a$ in the composite $\beta \circ_2 \alpha$,
\item take the 1-variable mate, and
\item evaluate at $d$.
\end{enumerate}

Step (1) is the same as taking the following horizontal composite:

\[\psset{unit=0.08cm,labelsep=0pt,nodesep=3pt}
\pspicture(0,-5)(90,30)

% a1 a2 a3 up right
% b1 b2 b3

\rput(0,25){\rnode{a1}{$B$}} % top left
\rput(45,25){\rnode{a2}{${B'}$}} % top middle
\rput(90,25){\rnode{a3}{${B'}$}} % top right

\rput(0,0){\rnode{b1}{$D\opb$}}   % bottom left
\rput(45,0){\rnode{b2}{${D'}\opb$}}  % bottom middle
\rput(90,0){\rnode{b3}{${D'}\opb$}}  % bottom right

\psset{nodesep=3pt,labelsep=2pt,arrows=->}
\ncline{a1}{a2}\naput{{\scriptsize $T$}} % top
\ncline{a2}{a3}\naput{{\scriptsize $1$}} % top

\ncline{b1}{b2}\nbput{{\scriptsize $V\opb$}} % bottom
\ncline{b2}{b3}\nbput{{\scriptsize $1$}} % bottom

\ncline{a1}{b1}\nbput{{\scriptsize $G_0(\suscore, Fa)$}} % left
\ncline{a2}{b2}\nbput[npos=0.6]{{\scriptsize ${G_0'}(\suscore, UFa)$}} % right
\ncline{a3}{b3}\naput{{\scriptsize ${G_0'}(\suscore, F'Sa)$}} % right

{\psset{doubleline=true,arrowinset=0.6,arrowlength=0.5,arrowsize=0.5pt 2.1,nodesep=0pt}
\rput[c](16,12){\pcline{->}(0,0)(4,4) \naput[labelsep=-1pt]{{\scriptsize $\beta_{\suscore,Fa}$}}}}

{\psset{doubleline=true,arrowinset=0.6,arrowlength=0.5,arrowsize=0.5pt 2.1,nodesep=0pt}
\rput[c](63,12){\pcline{->}(0,0)(4,4) \nbput[labelsep=0pt]{{\scriptsize $G_0'(\suscore,\alpha_a)$}}}}

% 
% 
% \psset{labelsep=1.5pt}
% \pnode(18,13){c1}
% \pnode(12,7){c2}
% \ncline[doubleline=true,arrowinset=0.6,arrowlength=0.8,arrowsize=0.5pt 2.1]{c1}{c2} \nbput[npos=0.4]{{\scriptsize $\alpha_{\uscore, b}$}}
% 
% 
% \psset{labelsep=1.5pt}
% \pnode(18,13){c1}
% \pnode(12,7){c2}
% \ncline[doubleline=true,arrowinset=0.6,arrowlength=0.8,arrowsize=0.5pt 2.1]{c1}{c2} \nbput[npos=0.4]{{\scriptsize $\alpha_{\uscore, b}$}}
% 

\endpspicture\]
and the axiom then follows from the fact that 1-variable mates respect horizontal composition, together with the fact that the mate of $G_0'(\uscore,\alpha_a)$ is
$G_1'(\alpha_a,\uscore)$; this last fact can be shown using a diagram chase involving dinaturality of a counit, functoriality of $G_1$, and a triangle identity for the 1-variable adjunction in question. \end{prf}

\begin{myremark}
 The direction of the natural tranformation for 2-cells is crucial so that the multivariable mates correspondence can be applied.  There is a cyclic double multicategory involving multivariable \emph{right} adjunctions in which the 2-cells must be given by natural transformations pointing down, as in
\[\psset{unit=0.08cm,labelsep=0pt,nodesep=3pt}
\pspicture(0,-2)(60,24)

% a1 a2 down left
% b1 b2

\rput(0,20){\rnode{a1}{$A_1 \times \cdots \times A_n $}} % top left
\rput(60,20){\rnode{a2}{$B_1 \times \cdots \times B_n$}} % top right

\rput(0,0){\rnode{b1}{$A_0\opb$}}   % bottom left
\rput(60,0){\rnode{b2}{${B_0}\opb$}}  % bottom right

\psset{nodesep=3pt,labelsep=2pt,arrows=->}
\ncline{a1}{a2}\naput{{\scriptsize $S_1 \times \cdots \times S_n$}} % top
\ncline{b1}{b2}\nbput{{\scriptsize $S_0\opb$}} % bottom
\ncline{a1}{b1}\nbput{{\scriptsize $F$}} % left
\ncline{a2}{b2}\naput{{\scriptsize $G$}} % right

\psset{labelsep=1.5pt}
\pnode(36,13){c1}
\pnode(30,7){c2}
\ncline[doubleline=true]{c1}{c2} \nbput[npos=0.4]{{\scriptsize $\alpha$}}

\endpspicture\]

To be precise we write $\MADJ_L$ for the multivariable left adjunctions and $\MADJ_R$ for the multivariable right adjunctions.  We will need the latter construction in the next section.

\end{myremark}

\begin{theorem}
 There is an isomorphism of double multicategories
\[(\ )\opb \: \MADJ_L \lra \MADJ_R.\]
\end{theorem}

This isomorphism is analogous to the isomorphism of double categories
\[\Ladj \cong \Ladj_R.\]
We now discuss isomorphisms analogous to the isomorphism of double categories
\[\Ladj \cong \Radj.\]
Recall that these double categories have the same 0- and 1-cells, but the 2-cells are natural transformations living in squares involving the left adjoints, for \Ladj, and the right adjoints, for \Radj.  For the $n$-variable version we have instead of left and right adjoints, a cycle of $n+1$ possible mutual adjoints.  This gives us many possible variants of the double cyclic multicategory \MADJ.  

% Note that in both \Ladj\ and \Radj\ the vertical 1-cells point in the direction of the \emph{left} adjoint in the adjunction; the \L\ and \R\ refer to the side of the adjunction involved in the 2-cells only.  Yet another variant would have the vertical 1-cells pointing in the direction of the \emph{right} adjoint in the adjunction; this is analogous to passing from $\MADJ_L$ to $\MADJ_R$.  **repetition from before?**

% (There an unfortunate contravariance here meaning that with mutual left adjunctions the vertical 1-cells point in the direction of the right adjoint, and vice versa.  This convention is chosen to avoid writing $\opb$ on every single source category of the multifunctor in the adjunction. **fix this**)

For the multivariable case the situation is further complicated by the fact that we have a choice of 2-cell convention for \emph{each arity $n$}, and these can all be chosen independently.  These choices are the $w_n$ in the following theorem.

\begin{theorem}
Suppose we have fixed for each $n \in \N$ an integer $w_n$ with
 $ 0 \leq w_n \leq n$. Write this infinite sequence of natural numbers as $\mathbf{w}$. Then we have a cyclic double multicategory $\MADJ_\bw$ with the same 0- and 1-cells as \MADJ\ (with multivariable left adjunctions, say) but where for each $n$ an $n$-ary 2-cell is as shown below
\[\psset{unit=0.08cm,labelsep=0pt,nodesep=3pt}
\pspicture(0,-3)(80,25)

% a1 a2 up right
% b1 b2

\rput(0,20){\rnode{a1}{$A_{w_n + 1} \times \cdots \times A_{w_n -1} $}} % top left
\rput(80,20){\rnode{a2}{$B_{w_n + 1} \times \cdots \times B_{w_n - 1}$}} % top right

\rput(0,0){\rnode{b1}{$A_{w_n}\opb$}}   % bottom left
\rput(80,0){\rnode{b2}{$B_{w_n}\opb$}}  % bottom right

\psset{nodesep=3pt,labelsep=2pt,arrows=->}
\ncline{a1}{a2}\naput{{\scriptsize $S_{w_n + 1} \times \cdots \times S_{w_n - 1}$}} % top
\ncline{b1}{b2}\nbput{{\scriptsize $S_{w_n}\opb$}} % bottom
\ncline{a1}{b1}\nbput{{\scriptsize $F_{w_n}$}} % left
\ncline{a2}{b2}\naput{{\scriptsize $G_{w_n}$}} % right

\psset{labelsep=1.5pt}
\pnode(39,7){c1}
\pnode(45,13){c2}
\ncline[doubleline=true]{c1}{c2} \nbput[npos=0.4]{{\scriptsize $\alpha$}}

\endpspicture\]
(note direction).  We emphasise that the horizontal source is still $F_0$ and the horizontal target is $G_0$; the vertical source is $S_1, \ldots, S_n$ and the vertical target is $S_0\opb$.  If each $w_n = 0$ we get the original version of \MADJ.

Composition proceeds via the mates correspondence.  

Then for all \bw\ there is an isomorphism of cyclic double multicategories
\[ \MADJ \cong \MADJ_{\bw}\]
which is the multivariable generalisation of the double category isomorphism
\[ \Ladj \cong \Radj.\]

\end{theorem}

\section{Application to algebraic monoidal model categories}\label{modelcat}

One aim of this work is to study an algebraic version of Hovey's notion of monoidal model category \cite{hov1}.  In such a model category we have hom and tensor structures that must interact well with the given model structure.  One such interaction requirement is that the 2-variable adjunction for hom and tensor should be a morphism of the underlying algebraic weak factorisation systems of the model category.  An important consequence of the defining axioms is that the total derived functors of the 2-variable adjunction given by the tensor and hom define a closed monoidal structure on the homotopy category of the model category. 

% The essential point in the proof is that the necessary coherence natural transformations, themselves maps between composite multivariable adjunctions expressing associativity and so forth, also derive. 

A model category has, among other things, two weak factoristion systems. In an \emph{algebraic} model category \cite{rie1} these are \emph{algebraic} weak factorisation systems \cite{gt1}. In this case, elements in the left and right classes of the weak factorisation systems specifying the model structure become coalgebras and algebras for the comonads and monads of the algebraic weak factorisation systems. An algebraic model category with a closed monoidal structure is a \emph{monoidal algebraic model category} \cite{rie2} just when the tensor/hom/cotensor 2-variable adjunction is a ``2-variable adjunction of algebraic weak factorisation systems''. This notion makes use of the definition of parametrised mates and motivates much of the present work.

As in \cite{rie1}, we abbreviate ``algebraic weak factorisation system'' to ``awfs''. First we recall the definition of awfs and of a standard (1-variable) adjunction of awfs. Throughout this section, given a category $A$ we write $A^\2$ for the category whose objects are morphisms of $A$, and whose morphisms are commuting squares. That is, $A^\2$ is the category $\mathbf{Cat}(\2, A)$ where $\2$ denotes the category containing a single non-trivial arrow. We have domain and codomain projections $\dom,\cod \colon A^\2 \tra A$.

A functorial factorisation on a category $A$ is given by a pair of functors $L,R \colon A^\2 \tra A^\2$ with $\dom L = \dom$, $\cod R = \cod$, and $\cod L = \dom R$. We call this last functor $E$, so we can write the factorisation of a morphism $f$ as below.
\[
\psset{unit=0.1cm,labelsep=3pt,nodesep=3pt}
\pspicture(20,17)

%%%%%%%%%% top

\rput(0,15){\rnode{a1}{$a$}}  % top left
\rput(20,15){\rnode{a2}{$b$}}  % top right
\rput(10,3){\rnode{a3}{$Ef$}}  % bottom

\ncline{->}{a1}{a2} \naput{{\scriptsize $f$}} % top
\ncline{->}{a1}{a3} \nbput{{\scriptsize $Lf$}} % left
\ncline{->}{a3}{a2} \nbput[labelsep=1pt]{{\scriptsize $Rf$}} % right

\endpspicture
\]

An awfs on a category $A$ is given by a functorial factorisation together with extra structure making

 \begin{itemize} 
\item $L$ a comonad on $A^\2$, and 
\item $R$ a monad on $A^\2$, such that
\item the canonical map $LR \tra RL$ given by multiplication and comultiplication is a distributive law. 
 \end{itemize} 
 The idea is that the $L$-coalgebras are the left maps (equipped with structure specifying their liftings) and the $R$-algebras are the right maps. 

% 
% 
% %%%%%%%%%%%%%%%%%%%%%%%%%%%%%%
% 
% As in **cite Emily 2** we abbreviate ``algebraic weak factorisation system'' to ``awfs''.
% 
% First we recall the definition of awfs and of a standard (1-variable) morphism of awfs.  Throughout this section, given a category $A$ we write $A^\2$ for the category whose objects are morphisms of $A$, and whose morphisms are commuting squares.  That is, $A^\2$ is the category $\Cat(\2, A)$ where $\2$ denotes the category containing a single non-trivial arrow.
% 
% Recall that an awfs on a category $A$ is given by a pair $(L,R)$ where
% \begin{itemize}
%  \item $L$ is a comonad on $A^\2$, and
% \item $R$ is a monad on $A^\2$,
% \end{itemize}
% satisfying certain axioms.  The idea is that the $L$-coalgebras are the left maps (equipped with structure specifying their liftings) and the $R$-algebras are the right maps.  Any map $a \tmap{f} b$ then factors as $Lf$ followed by $Rf$, and it is convenient to write the object in the ``middle'' as $Ef$---in fact there is a functor
% \[E: A^\2 \tra A\]
% arising canonically from the awfs.  Thus we have the following factorisation of $f$.
% \[
% \psset{unit=0.1cm,labelsep=3pt,nodesep=3pt}
% \pspicture(20,17)
% 
% %%%%%%%%%% top
% 
% \rput(0,15){\rnode{a1}{$a$}}  % top left
% \rput(20,15){\rnode{a2}{$b$}}  % top right
% \rput(10,3){\rnode{a3}{$Ef$}}  % bottom
% 
% \ncline{->}{a1}{a2} \naput{{\scriptsize $f$}} % top
% \ncline{->}{a1}{a3} \nbput{{\scriptsize $Lf$}} % left
% \ncline{->}{a3}{a2} \nbput[labelsep=1pt]{{\scriptsize $Rf$}} % right
% 
% \endpspicture
% \]
% 

\begin{mydefinition}
 A adjunction of awfs
\[\begin{array}{ccc}
 (L_1, R_1) & \tra & (L_2, R_2) \\
\mbox{on \ } A_1 && \mbox{on \ } A_2
\end{array}\]
consists of the following.

\begin{itemize}
 \item An adjunction \[
\psset{labelsep=2pt}
\pspicture(20,8)
\rput(0,3){\rnode{A}{$A_1$}}
\rput(20,3){\rnode{B}{$A_2$}}
\psset{nodesep=2pt,arrows=->, labelsep=2pt}
\ncline[offset=5pt]{A}{B}\nbput{$\scriptstyle\perp$}\naput{{\scriptsize $F$}}
\ncline[offset=5pt]{B}{A}\naput{{\scriptsize $G$}}
\endpspicture
\]

\item Natural transformations $\lambda$ and $\rho$ making
\begin{enumerate}
 \item $(F^\2, \lambda)$ into a colax comonad map $L_1 \tra L_2$, and
\item $(G^\2, \rho)$ into a lax monad map $R_2\tra R_1$
\end{enumerate}
where
\[\begin{array}{cl}
   \lambda = (1, \alpha), & \mbox{and}\\
   \rho = (\overline{\alpha}, 1).
  \end{array}\]
\end{itemize}

\end{mydefinition}

Here $\overline{\alpha}$ denotes the mate of $\alpha$, about which some further comments are called for.  \emph{A priori} the natural transformations $\lambda$ and $\rho$ are as shown below
\[\psset{unit=0.08cm,labelsep=0pt,nodesep=3pt}
\pspicture(0,-3)(40,22)

% a1 a2 down left
% b1 b2

\rput(0,20){\rnode{a1}{$A_1^\2$}} % top left
\rput(20,20){\rnode{a2}{$A_1^\2$}} % top right

\rput(0,0){\rnode{b1}{$A_2^\2$}}   % bottom left
\rput(20,0){\rnode{b2}{$A_2^\2$}}  % bottom right

\psset{nodesep=3pt,labelsep=2pt,arrows=->}
\ncline{a1}{a2}\naput{{\scriptsize $L_1$}} % top
\ncline{b1}{b2}\nbput{{\scriptsize $L_2$}} % bottom
\ncline{a1}{b1}\nbput{{\scriptsize ${F}^\2$}} % left
\ncline{a2}{b2}\naput{{\scriptsize ${F}^\2$}} % right

\psset{labelsep=1.5pt}
\pnode(13,13){c1}
\pnode(7,7){c2}
\ncline[doubleline=true]{c1}{c2} \naput[npos=0.4]{{\scriptsize $\lambda$}}
\endpspicture
\pspicture(0,-3)(30,22)

% a1 a2 down left
% b1 b2

\rput(0,20){\rnode{a1}{$A_1^\2$}} % top left
\rput(20,20){\rnode{a2}{$A_1^\2$}} % top right

\rput(0,0){\rnode{b1}{$A_2^\2$}}   % bottom left
\rput(20,0){\rnode{b2}{$A_2^\2$}}  % bottom right

\psset{nodesep=3pt,labelsep=2pt,arrows=->}
\ncline{a1}{a2}\naput{{\scriptsize $R_1$}} % top
\ncline{b1}{b2}\nbput{{\scriptsize $R_2$}} % bottom
\ncline{b1}{a1}\naput{{\scriptsize ${G}^\2$}} % left
\ncline{b2}{a2}\nbput{{\scriptsize ${G}^\2$}} % right

\psset{labelsep=1.5pt}
\pnode(7,13){c1}
\pnode(13,7){c2}
\ncline[doubleline=true]{c1}{c2} \nbput[npos=0.4]{{\scriptsize $\rho$}}
\endpspicture
\]
but it turns out that such $\lambda$ and $\rho$ are completely determined by respective natural transformations as below

\[\psset{unit=0.08cm,labelsep=0pt,nodesep=3pt}
\pspicture(0,-3)(40,22)

% a1 a2 down left
% b1 b2

\rput(0,20){\rnode{a1}{$A_1^\2$}} % top left
\rput(20,20){\rnode{a2}{$A_1$}} % top right

\rput(0,0){\rnode{b1}{$A_2^\2$}}   % bottom left
\rput(20,0){\rnode{b2}{$A_2$}}  % bottom right

\psset{nodesep=3pt,labelsep=2pt,arrows=->}
\ncline{a1}{a2}\naput{{\scriptsize $E_1$}} % top
\ncline{b1}{b2}\nbput{{\scriptsize $E_2$}} % bottom
\ncline{a1}{b1}\nbput{{\scriptsize ${F}^\2$}} % left
\ncline{a2}{b2}\naput{{\scriptsize ${F}$}} % right

\psset{labelsep=1.5pt}
\pnode(13,13){c1}
\pnode(7,7){c2}
\ncline[doubleline=true]{c1}{c2} \naput[npos=0.4]{{\scriptsize $$}}
\endpspicture
\pspicture(0,-3)(30,22)

% a1 a2 down left
% b1 b2

\rput(0,20){\rnode{a1}{$A_1^\2$}} % top left
\rput(20,20){\rnode{a2}{$A_1^\2$}} % top right

\rput(0,0){\rnode{b1}{$A_2^\2$}}   % bottom left
\rput(20,0){\rnode{b2}{$A_2^\2.$}}  % bottom right

\psset{nodesep=3pt,labelsep=2pt,arrows=->}
\ncline{a1}{a2}\naput{{\scriptsize $E_1$}} % top
\ncline{b1}{b2}\nbput{{\scriptsize $E_2$}} % bottom
\ncline{b1}{a1}\naput{{\scriptsize ${G}^\2$}} % left
\ncline{b2}{a2}\nbput{{\scriptsize ${G}$}} % right

\psset{labelsep=1.5pt}
\pnode(7,13){c1}
\pnode(13,7){c2}
\ncline[doubleline=true]{c1}{c2} \nbput[npos=0.4]{{\scriptsize $$}}
\endpspicture
\]
It is these that are required to be mates $\alpha$ and $\overline{\alpha}$ respectively, under the adjunctions $F^\2 \ladj G^\2$ and $F \ladj G$.  (Note that $(\ )^\2$ is actually the 2-functor $\Cat(\2, \uscore)$ so preserves adjunctions.)

It turns out that the appropriate generalisation for the $n$-variable case involves generalising the functor $(\ )^\2$ as well, as follows.

\begin{mydefinition} 
 Let $F \: A_1 \times \cdots \times A_n \tra A_0$ be an $n$-variable functor.  We define a functor
\[\hat{F} \: A_1^\2 \times \cdots \times A_n^\2 \tra A_0^\2\]
as follows.  Consider morphisms
\[a_{i0} \tmap{f_i} a_{i1} \in A_i\]
for each $1 \leq i \leq n$.  We need to define a morphism $\hat{F}(f_1, \ldots, f_n)$ in $A_0$.  Consider the commuting hypercube in $A_1 \times \cdots \times A_n$ built from $f_i$'s as follows.
\begin{itemize}
 \item Vertices are given by $(a_{1k_1}, \ldots, a_{nk_n})$ where each $k_i=0$ or 1 (thus, the $i$th term is either the source or target of $f_i$).

\item Edges are given by $(1, \ldots, 1, f_i, 1, \ldots, 1)$ for some $1 \leq i \leq n$.
\end{itemize}
This hypercube commutes since each path from
\[(a_{10}, \ldots, a_{n0}) \mbox{ \ to \ } (a_{11}, \ldots, a_{n1})\]
composes to $(f_1, \ldots, f_n)$.  

We apply $F$ to this diagram and take the ``obstruction'' map induced by the colimit over the diagram
\[ \big\{ (a_{10}, \ldots, a_{n0}) \psset{unit=0.1cm,nodesep=0pt} \pspicture(27,4)
\pcline{->}(2.5,1.1)(24.5,1.1)\naput[labelsep=2pt]{\ensuremath{\scriptstyle{(1, \ldots, 1, f_i, 1, \ldots, 1)}}} \endpspicture
 (a_{10}, \ldots, a_{i1}, \ldots, a_{n0}) \big\}_{0 \leq i \leq n} \]
and call this map $\hat{F}(f_1, \ldots, f_n)$ in ${A_0}^\2$; its domain is the above colimit and its codomain is $(a_{11}, \ldots, a_{n1})$.

  The action on morphisms is then induced in the obvious way.  In fact $\hat{(\ )}$ is a pseudo-functor so preserves adjunctions.  Furthermore, a straightforward but notationally involved proof shows that $\hat{(\ )}$ preserves $n$-variable adjunctions, as we first learned from Dominic Verity.
\end{mydefinition}

\begin{myremark}
Given an awfs $(L,R)$ on a category $A$, we get a dual awfs $(R\opb, L\opb)$ on $A\opb$.  Note that
\begin{itemize}
 \item $L$ is a comonad on $A$, so $L\opb$ is a monad on $A\opb$, and
\item $R$ is a monad on $A$, so $R\opb$ is a comonad on $A\opb$.
\end{itemize}
Also, given awfs $(L_1,R_1)$ on $A_1$ and $(L_2, R_2)$ on $A_2$ we get an awfs \[(L_1 \times L_2, R_1 \times R_2)\] on $A_1 \times A_2$. 
\end{myremark}

\begin{mydefinition}\label{morphawfs}
Suppose we have for each $0 \leq i \leq n$ a category $A_i$ equipped with an awfs $(L_i, R_i)$.  Then an \demph{$n$-variable adjunction of awfs}
\[A_1 \times \cdots \times A_n \tra A_0\opb\]
is given by the following.
\begin{itemize}
 \item A functor $F_0: A_1 \times \cdots \times A_n \tra A_0\opb$ equipped with $n$-variable \emph{right} adjoints $F_1, \ldots, F_n$, and

\item For each $i$ a natural transformation $\lambda_i$ as shown below

\[\psset{unit=0.08cm,labelsep=0pt,nodesep=3pt}
\pspicture(0,-3)(70,22)

% a1 a2 down left
% b1 b2

\rput(0,20){\rnode{a1}{$A_{i+1}^\2 \times \cdots \times A_{i-1}^\2$}} % top left
\rput(70,20){\rnode{a2}{$A_{i+1}^\2 \times \cdots \times A_{i-1}^\2$}} % top right

\rput(0,0){\rnode{b1}{${A_i\opb}^\2$}}   % bottom left
\rput(70,0){\rnode{b2}{${A_i\opb}^\2$}}  % bottom right

\psset{nodesep=3pt,labelsep=2pt,arrows=->}
\ncline{a1}{a2}\naput{{\scriptsize $L_{i+1} \times \cdots \times L_{i-1}$}} % top
\ncline{b1}{b2}\nbput{{\scriptsize $R_i\opb$}} % bottom
\ncline{a1}{b1}\nbput{{\scriptsize $\hat{F}_i$}} % left
\ncline{a2}{b2}\naput{{\scriptsize $\hat{F}_i$}} % right

\psset{labelsep=1.5pt}
\pnode(38,13){c1}
\pnode(32,7){c2}
\ncline[doubleline=true]{c1}{c2} \naput[npos=0.4]{{\scriptsize $\lambda_i$}}

\endpspicture\]
making $(\hat{F}_i, \lambda_i)$ into a colax comonad map
\[L_{i+1} \times \cdots \times L_{i-1} \tra R_i\opb.\]

\item As in the 1-variable case, such a $\lambda_i$ is completely determined by a natural transformation

\[\psset{unit=0.08cm,labelsep=0pt,nodesep=3pt}
\pspicture(0,-3)(70,22)

% a1 a2 down left
% b1 b2

\rput(0,20){\rnode{a1}{$A_{i+1}^\2 \times \cdots \times A_{i-1}^\2$}} % top left
\rput(70,20){\rnode{a2}{$A_{i+1}\times \cdots \times A_{i-1}$}} % top right

\rput(0,0){\rnode{b1}{${A_i\opb}^\2$}}   % bottom left
\rput(70,0){\rnode{b2}{$A_i\opb$}}  % bottom right

\psset{nodesep=3pt,labelsep=2pt,arrows=->}
\ncline{a1}{a2}\naput{{\scriptsize $E_{i+1} \times \cdots \times E_{i-1}$}} % top
\ncline{b1}{b2}\nbput{{\scriptsize $E_i\opb$}} % bottom
\ncline{a1}{b1}\nbput{{\scriptsize $\hat{F}_i$}} % left
\ncline{a2}{b2}\naput{{\scriptsize ${F}_i$}} % right

\psset{labelsep=1.5pt}
\pnode(38,13){c1}
\pnode(32,7){c2}
\ncline[doubleline=true]{c1}{c2} \naput[npos=0.4]{{\scriptsize $\alpha_i$}}

\endpspicture\]
and we require the $\alpha_i$ to be parametrised mates.

\end{itemize}

\end{mydefinition}

% \emph{A priori} only 2-variable adjunctions of awfs are required to make the definition of monoidal algebraic model category, using the motivating example of hom and tensor (Section~\ref{homtensor}).  However, in practice the higher arity versions are useful in the following way.  Given a $\cV$-category $\cC$ tensored over $\cV$ we may wish to know if the functor
% \[V \times \uscore \: \cC \tra \cC\]
% is actually a $\cV$-functor.  For this it is equivalent to study the functors
% \[\begin{array}{ccl}
%    \cV \times \cV \times \cC & \tra & \cC \\
%    (V, W, C)  & \tmapsto & V \otimes (W \otimes C) \\
%    (V, W, C) & \tmapsto & (V \otimes W) \otimes C
%   \end{array}\]
% and ask for ``associativity'', that is, for these functors to be equal.  In that case the functor is part of a 3-variable adjunction.  
% 
% 

\begin{myexample} 
An algebraic, or perhaps constructive, encoding of the classical result that the simplicial hom-space from a simplicial set $A$ to a Kan complex $X$ is again a Kan complex is that the tensor-hom 2-variable adjunction is a 2-variable adjunction of awfs. This example is prototypical, so we explain it further. The sets of maps 
\[I = \{ \partial\Delta^n \tra \Delta^n \mid n \geq 0\}\] 
and 
\[J = \{ \Lambda^n_k \tra \Delta^n \mid n \geq 1, 0 \leq k \leq n\}\]
 generate two awfs $(C,F_t)$ and $(C_t, F)$ on \sSet\ by Garner's algebraic small object argument \cite{gar2}. A simplicial set $X$ is a \textbf{Kan complex} if the unique map $X \tra \Delta^0$ satisfies the right lifting property with respect to $J$.

The sets $I$ and $J$ determine the cofibrations and fibrations in Quillen's model structure on \sSet, which is a monoidal algebraic model category. The key technical step in the proof of this fact is that the 2-variable morphism 
\psset{labelsep=2pt}
\[\begin{array}{ccc} 
\sSet^\2 \times \sSet^\2 & \ltmap{-\hat{\times}-} & \sSet^\2 \\ 
(C,F_t) \times (C_t,F) & \ltmapsto & (C_t,F)
  \end{array}\]
induced from the cartesian product is part of a 2-variable adjunction of awfs.

The modern proof of the non-algebraic version of this result makes use of the closure properties of left classes of weak factorisation systems and is non-constructive; see \cite{gj1}. This argument does not suffice to prove the algebraic statement. However, the classical constructive proof does suffice:  the proof given in  \cite[Theorem 6.9]{may3} explicitly constructs the required lifts of $\mathbf{Hom}(A,X) \tra \Delta^0$ against $J$, supposing that similar lifts for $X \to \Delta^0$ are given. By the main result of \cite{rie2}, this argument shows that the 2-variable right adjoint 

\[\begin{array}{ccc} 
(\mathbf{sSet}^\2)^\bullet \times \mathbf{sSet}^\2 &\ltmap{\hat{\mathbf{Hom}}} & \mathbf{sSet}^\2 \\ 
(F_t, C) \times (C_t,F) &\ltmapsto & (C_t,F)
\end{array}\] 
defines a 2-variable adjunction of awfs. By our main theorem (Theorem~\ref{mainthm}) this is equivalent to the desired statement. See \cite{rie2} for more details.
\end{myexample}

An important corollary of our main theorem in this context is the following result.

\begin{thm} Multivariable adjunctions of awfs compose to yield new multivariable adjunctions of awfs.
\end{thm}
\begin{prf} 
Multivariable colax comonad morphisms compose multicategorically.  Using the notation of Definition~\ref{morphawfs}, the composite is obtained by composing the $F_i$ and the $\lambda_i$ in the obvious way.  

Now, by the relationship between the $\lambda_i$ and the $\alpha_i$, the composite of the $\lambda_i$ is determined by the multicategorical composite of the $\alpha_i$.  So we check that these composites satisfy the mate condition required by the definition.  This follows from Theorem~\ref{matesvert}.
\end{prf}

% Assuming that sources and targets match up, the constituent multifunctors compose in the way detailed in Section \ref{sec:multi}. By a well-known result of Appelgate, colax comonad morphisms correspond to lifted functors between the associated categories of coalgebras \cite{johnstoneadjoint}; the latter clearly compose as well. It remains only to show that the parametrised mates of the pasted composite of the natural transformations of Definition \ref{defn:nvaradjawfs} is the pasted composite of the parametrised mates. But this is precisely what we have proven.
% \end{proof}

While only 2-variable adjunctions of awfs are required to make the definition of a monoidal algebraic model category, the higher arity versions are useful in the following way. Enriched categories, functors, adjunctions, and 2-variable adjunctions over a closed symmetric monoidal category $\cV$ can be encoded by an \emph{a priori} unenriched tensor/hom/cotensor 2-variable adjunction together with coherence isomorphisms.  These are isomorphisms between various composite 2-, 3- and 4-variable functors \cite{shul1}. There are many equivalent ways to encode this data having to do with choices of left and right adjoints. Our main result allows a seamless translation between these equivalent formulations. Related considerations arise in homotopy theory where these arguments may be used to prove that the total derived functor of a $\cV$-functor between $\cV$-model categories admits a canonical enrichment over the homotopy category of $\cV$.

%\bibliography{../bib/bib1208}

%\begin{thebibliography}{10}

\ed

%%%%%%%%%%%%%%%%%%%%%%%%%%%%%%%%%%%%%%%%%%%%%%%%%%%%%%%%%%%%%%%%%%%%%%%%%%%%%%%%%%%%%%%%%%%%%%%%%%%%%%%%%%%%%%%%%%%%%%%%%%%%%%%%%%%%%%%%%%%%%%%%%%%%%%%%%%%%%%%%%%%%%%%%%%%%%%%%%%%%%%%%%%%%%%%%%%%%%%%%%%%%%%%%%%%%%%%%%%%%%%%%%%%%

\[
\psset{unit=0.1cm,labelsep=2pt,nodesep=2pt}
\pspicture(,)
\pnode(,){}             % named node, nothing placed there
\rput(,){\rnode{}{$$}}  % named node, with something placed there
\rput(,){$$}            % unnamed node, something placed there
\nput{0}{}{}    % put things at an existing node: angle, nodename, stuff
\ncline{->}{}{} \naput{{\scriptsize $$}}
\ncline{}{} \naput{{\scriptsize $$}}
\psline*[par]{arrows}(x0,y0)(x1,y1)(xn,yn)
\ncarc[arcangle=35]{->}{}{}
\nccurve[nodesep=,angleA=,angleB=,ncurvA=,ncurvB=]{LA}{LB}\naput{{\scriptsize $$}}
\endpspicture
\]

just an arrow

\[
\psset{unit=0.1cm,labelsep=2pt,nodesep=2pt}
\pspicture(,)
\rput(,){\rnode{}{$$}}  % named node, with something placed there
\ncline{->}{}{} \naput{{\scriptsize $$}}
\endpspicture
\]

2-cell for placing

\psset{doubleline=true}
\rput{0}(21,17){\pcline{->}(4,4)(0,0) \naput[labelsep=1pt]{{\scriptsize $\varepsilon_{b_1}$}}}